\documentclass[moor,nonblindrev]{informs3}
\OneAndAHalfSpacedXII

\usepackage[utf8]{inputenc}
\usepackage[american]{babel}
\usepackage{bm,mathrsfs}
\usepackage{booktabs,tabularx,multirow}
\usepackage{enumitem}
\usepackage{float}
\floatstyle{ruled}
\newfloat{algorithm}{tbp}{loa}
\floatname{algorithm}{Algorithm}
\usepackage{placeins}
\usepackage{microtype}
\usepackage{etoolbox}

\definecolor{myblue}{RGB}{0,20,114}
\usepackage[hidelinks,colorlinks=true,urlcolor=myblue,linkcolor=myblue,citecolor=myblue,hypertexnames=false]{hyperref}
\def\EMAIL#1{\href{mailto:#1}{#1}}

\def\theARTICLETOP{}
\def\theLRHFirstLine{}
\def\theLRHSecondLine{}
\def\theRRHFirstLine{}
\def\theRRHSecondLine{}

\usepackage{natbib}
\bibpunct[, ]{(}{)}{,}{a}{}{,}
\def\bibfont{\small}

\TheoremsNumberedThrough
\EquationsNumberedThrough

\let\classproof\proof
\let\classendproof\endproof
\renewenvironment{proof}[1][Proof]
  {\classproof{#1.}}
  {\halmos\classendproof}

\setlist[enumerate]{itemsep=0pt,parsep=0pt,topsep=4pt}

\newcommand{\E}{\mathbb{E}}
\newcommand{\I}{\mathbb{I}}
\newcommand{\PP}{\mathbb{P}}

\allowdisplaybreaks[4]

\begin{document}
\RUNAUTHOR{Wang, Cao, Hu, Hu, and Xia}
\RUNTITLE{Long-Run CVaR Reinforcement Learning}
\TITLE{Long-Run Conditional Value at Risk Reinforcement Learning}
\ARTICLEAUTHORS{%
\AUTHOR{Qixin Wang}
\AFF{School of Management, Fudan University, Shanghai, 200433, China, \EMAIL{23110690014@m.fudan.edu.cn}}
\AUTHOR{Hao Cao}
\AFF{School of Management, Fudan University, Shanghai, 200433, China, \EMAIL{hcao21@m.fudan.edu.cn}}
\AUTHOR{Jian-Qiang Hu}
\AFF{School of Management, Fudan University, Shanghai, 200433, China, \EMAIL{hujq@fudan.edu.cn}}
\AUTHOR{Mingjie Hu}
\AFF{School of Management, Fudan University, Shanghai, 200433, China, \EMAIL{23110690009@m.fudan.edu.cn}}
\AUTHOR{Li Xia}
\AFF{School of Business, Sun Yat-Sen University, Guangzhou, 510275, China, \EMAIL{xiali5@sysu.edu.cn}}
}
\ABSTRACT{%
We propose a reinforcement learning method for minimizing the long-run conditional value at risk (CVaR) in Markov decision processes with finite state and action spaces. The method uses a single sample trajectory to update a quantile estimate, a $Q$-function, and a randomized policy, without requiring knowledge of transition probabilities or cost distributions. Under appropriate conditions, we show that the estimates track their changing targets as the policy is updated and exploration decreases, even when the proportion of updates for some state and action pairs tends to zero. With a strict action gap, the policy converges to a locally optimal deterministic policy if it enters a suitable neighborhood sufficiently late. We also extend the method to an objective combining CVaR and mean cost. Numerical experiments illustrate its performance in machine replacement and energy storage and compare it with existing methods.
}
\KEYWORDS{Long-Run CVaR, Reinforcement Learning, Stochastic Approximation, Risk Sensitive, Markov Processes}
\maketitle

\section{Introduction}\label{sec_intro}
Reinforcement learning (RL) provides a useful approach for improving the performance of dynamic systems when their transition mechanisms and cost distributions are unavailable. Most classical methods evaluate a policy through an expected cumulative criterion. In many long-run decision problems, however, a policy with favorable mean performance may still generate large losses with small probability. For such problems, the tail behavior of the cost distribution is a natural part of the performance measure.

Conditional value at risk (CVaR) measures the expected loss in the upper tail of a distribution and admits a convenient optimization representation \citep{rockafellar2000optimization,brown2007large}. Existing studies have developed simulation and stochastic approximation (SA) procedures for estimating and optimizing value at risk (VaR) and CVaR in static stochastic systems \citep{uryasev2001conditional,hong2009simulating,tamar2015optimizing,hu2024quantile,hu2024simulation,wang2025}. Existing procedures may also involve output distributions that depend on changing decision parameters. In our setting, however, the target is the stationary CVaR of a controlled Markov chain, whereas observations are generated along a single trajectory under an evolving policy. The additional difficulty is to relate these observations to the stationary CVaR of the current policy.

For dynamic systems, much of the CVaR literature considers discounted cumulative cost, for which dynamic programming and value iteration methods are available \citep{pflug2016time,chow2015risk}. We instead study the CVaR of the stationary one stage cost distribution. This criterion evaluates the tail cost of an operating period observed after the controlled system has reached steady state. It is different from the CVaR of discounted return, cumulative cost, or a path dependent loss, and it does not distinguish isolated high cost periods from clusters having the same stationary marginal distribution. The criterion is useful when the operating objective is to control the frequency and severity of costly periods rather than the tail of an accumulated return. Model-based formulations have been studied through Bellman local optimality equations and policy iteration \citep{xia2023risk}. These results provide the fixed point relations used in our analysis. \citet{bonetti2023risk} optimize reward-based CVaR using policy gradient and trust region constructions, whereas we study tabular SA that explicitly tracks the policy dependent stationary VaR and the associated fixed threshold control $Q$-function along a single adaptive trajectory. The distinction concerns the recursive objects and their analysis, rather than a general requirement for known transition probabilities and cost distributions.

Risk sensitive RL has also been studied for exponential criteria, utility functions, dynamic risk measures, and the CVaR of accumulated costs \citep{borkar2010learning,prashanth2022risk,tamar2015optimizing,chow2018risk,jiang2018risk,chen2023provably,stanko2021cvar}. The tabular reward-extended differential (RED) CVaR Q learning method of \citet[Theorem C.1.1]{rojas2025burning} concerns the limiting one stage reward. With $R=-C$ and lower tail probability $1-\phi$, the continuous distribution criterion agrees, up to sign, with ours. With a frequency balance requirement, it establishes convergence to optimal CVaR quantities and a Bellman solution modulo a constant. Our results instead concern moving VaR and Q targets under decreasing exploration. Removing the requirement is accompanied by weighted contraction and component step size conditions, and our policy conclusion is local. 

The main difficulty in our setting is that the VaR entering the transformed one stage cost depends on the stationary distribution generated by the current policy. Because the policy changes during learning, the observed state and cost sequence is nonhomogeneous. The VaR estimate, the relative $Q$-function, and the policy therefore form three coupled recursive processes. To address this difficulty, we develop a tabular SA method based on the Bellman local optimality equation. The basic idea is to view the VaR and Q calculations as a root finding problem and a fixed point problem and to approximate their solutions through separate but coupled recursions. A third recursion moves the policy toward the minimizing actions of the current Q estimate. The recursions use distinct step sizes and only one sample trajectory, in the same spirit as multi-timescale SA \citep{Hu2022,borkar2008stochastic,kushner2003sa}. A componentwise projection bounds the Q estimates, while the policy is adjusted incrementally rather than replaced by a greedy policy.

Another difficulty is that changes in the policy also change its stationary VaR and the corresponding control $Q$-function. Classical asynchronous stochastic approximation already allows unequal update frequencies for a fixed contraction map, and tracking a changing parameter has been studied in average cost learning (see, e.g., \citealp{abounadi2001learning}). Here, we need to account for these changes while updating each $Q$ entry only when its state and action pair is visited and projecting the updated value onto a bounded interval. Some pairs may be visited infinitely often but account for a fraction of the observations that tends to zero.

We analyze the three recursions in turn. For VaR estimation, we use a Poisson equation to handle the dependence in observations generated by the changing policy. This establishes that the estimates remain bounded and track the stationary VaR of the current policy, without assuming that the policy converges. For the $Q$-function estimator, we measure the progress of each entry by the step sizes accumulated at its own updates. Under appropriate conditions, we bound the estimation error over successive update cycles, accounting for the random noise, the changes in the target, and the projection. This proves that the $Q$ estimates track the changing control $Q$-function even when relative update frequencies tend to zero. Finally, we show that if the policy iterates enter a suitable neighborhood of a deterministic policy satisfying a strict action-gap condition sufficiently late, they remain in that neighborhood and converge to it. This policy is locally optimal among stationary randomized policies. Table~\ref{tab:closest_methods} compares our method with \citet{xia2023risk, rojas2025burning} in more detail.

\begin{table}[t]
\centering
\caption{Comparison of stationary CVaR methods.}
\label{tab:closest_methods}
\small
\setlength{\tabcolsep}{2.5pt}
\renewcommand{\arraystretch}{1.08}
\begin{tabular}{@{}>{\raggedright\arraybackslash}p{2.1cm}>{\raggedright\arraybackslash}p{2.5cm}>{\raggedright\arraybackslash}p{3.5cm}>{\raggedright\arraybackslash}p{3.2cm}>{\raggedright\arraybackslash}p{3.2cm}@{}}
\toprule
Study & Long-run risk object & Information and recursion & Conditions used in the analysis & Main conclusion \\
\midrule
\citet{xia2023risk} & Upper tail CVaR of the stationary one stage cost & Known transition and cost laws; sensitivity analysis and policy iteration & Ergodicity and the regularity conditions of the model-based analysis & Bellman local optimality characterization; deterministic optimality and a global algorithm \\\addlinespace[3pt]
\citet{rojas2025burning}, Theorem C.1.1 & Lower tail CVaR of the limiting one stage reward & Online temporal difference recursions for Q, VaR, and CVaR & Communicating model; bounded rewards; uniqueness; balanced infinite updates & Optimal CVaR quantities and Q modulo a constant \\\addlinespace[3pt]
This paper & Upper tail CVaR of the stationary one stage cost & Single trajectory recursions for VaR, $Q$-function, and policy & Local density positivity; finite second moments; uniform minorization; weighted contraction & Moving target tracking; policy convergence only after local entry with a strict gap \\
\bottomrule
\end{tabular}
\end{table}

Our main contributions are summarized as follows.
\begin{enumerate}
    \item We propose a reinforcement learning method for stationary cost CVaR that combines quantile estimation, $Q$ learning, and policy improvement along a single sample trajectory without knowing the transition probabilities or cost distributions. We establish the asymptotic behavior of the quantile and $Q$-function estimators without prior policy convergence or balanced update frequencies, and prove convergence to a deterministic policy that is locally optimal among stationary randomized policies when a strict action gap holds and the iterates enter a suitable neighborhood sufficiently late. We also extend the method to an objective combining CVaR and mean cost.

    \item We compare the proposed method with RED through numerical experiments. The results show that our method yields smaller average CVaR gaps and fewer observed changes in the greedy policies than RED.
\end{enumerate}

The rest of this paper is organized as follows. Section~\ref{sec_prob_alg} introduces the stationary one stage CVaR Markov decision process and develops our algorithm. Section~\ref{sec_converge} establishes the asymptotic results and the local stability result and then gives the mean CVaR extension. Section~\ref{sec_experi} reports the numerical experiments. Section~\ref{sec_conclude} concludes the paper.

\section{Reinforcement Learning for Long-Run CVaR}\label{sec_prob_alg}
This section introduces the stationary one stage cost CVaR criterion and then presents the recursive learning procedure. We first state the model and its steady state interpretation, then recall the Bellman local optimality equation, and finally identify the quantile and fixed point problems approximated by the recursions. This separation distinguishes the model-based characterization from the sample-based updates. 

\subsection{Preliminaries}\label{sec_prob}
We consider an infinite horizon discrete time Markov decision process. The finite state space is $\mathcal S=\{s^{(1)},\ldots,s^{(|\mathcal S|)}\}$, and $\mathcal A(s)$ denotes the finite set of actions available at state $s$. Let $\mathcal A=\bigcup_{s\in\mathcal S}\mathcal A(s)$. We use $s^{(0)}$ to denote one selected state in $\mathcal S$ that serves as the reference state. For every feasible pair $(s,a)$, let $\mathsf K(B,s'\mid s,a)$ denote the joint probability that the one stage cost belongs to a Borel set $B$ and the next state is $s'$. Its cost and state marginals are denoted by $F(\cdot;s,a)$ and $p(\cdot\mid s,a)$, respectively, so that $\sum_{s'\in\mathcal S}\mathsf K((-\infty,x],s'\mid s,a) =F(x;s,a),~ \mathsf K(\mathbb R,s'\mid s,a) =p(s'\mid s,a)$. Both marginals and their joint law are unknown.

\begin{assumption}\label{A.cost}
For every feasible pair $(s,a)$, the cost marginal $F(\cdot;s,a)$ is absolutely continuous and has a finite second moment. The joint transition and cost law $\mathsf K(\cdot,\cdot\mid s,a)$ is time homogeneous.
\end{assumption}

\begin{remark}
Assumption~\ref{A.cost} specifies the continuous cost model and supplies the second moment needed for the martingale estimates. The joint law allows the cost and next state to be conditionally dependent. It is satisfied by Gaussian, Student $t$ with more than two degrees of freedom, logistic, and bounded continuous cost distributions. 
\end{remark}

A stationary randomized policy is a mapping $d$ with $d(s)\in\Delta(\mathcal A(s))$, where $d(s,a)$ is the probability of selecting action $a\in\mathcal A(s)$ at state $s$. Let $\mathcal D$ denote the collection of all such policies. For $a\in\mathcal A(s)$, $\delta(a)$ denotes the probability vector with entry one at $a$ and zero at every other feasible action. For $d\in\mathcal D$, define $P^d(s,s'):=\sum_{a\in\mathcal A(s)}d(s,a)p(s'\mid s,a)$. Under Assumption~\ref{A.mixing}, $P^d$ has a unique stationary distribution $\mu^d$, and $\pi^d(s,a):=\mu^d(s)d(s,a)$ is the stationary state and action distribution.

For a fixed level $\phi\in(0,1)$, the VaR of the cost at time $n$ under policy $d$ and initial state $s$ is $\operatorname{VaR}_\phi^d[C(s_n,a_n)\mid s_0=s]:=\inf\{x:\PP^d(C(s_n,a_n)\le x\mid s_0=s)\ge\phi\}$. The corresponding CVaR is defined by
\begin{align*}
\operatorname{CVaR}_\phi^d\left[C\left(s_n,a_n\right)\mid s_0=s\right] &=\frac{1}{1-\phi}\int_\phi^1\operatorname{VaR}_\varphi^d\left[C\left(s_n,a_n\right)\mid s_0=s\right]d\varphi.
\end{align*}
Because $\phi$ is fixed, its subscript is suppressed hereinafter. The long-run VaR and CVaR criteria under policy $d$ are defined as: $\operatorname{VaR}^d(s) :=\lim_{N\to\infty} (1/N) \sum_{n=0}^{N-1}\operatorname{VaR}^d[C(s_n,a_n)\mid s_0=s]$, $\operatorname{CVaR}^d(s) :=\lim_{N\to\infty} (1/N) \sum_{n=0}^{N-1}\operatorname{CVaR}^d[C(s_n,a_n)\mid s_0=s]$. In this paper, we consider $\min_{d\in\mathcal D}\operatorname{CVaR}^d(s)$. Let $C^d$ denote the stationary one stage cost, whose cumulative distribution function is $\PP(C^d\le x)=\sum_{s\in\mathcal S}\sum_{a\in\mathcal A(s)}\pi^d(s,a)F(x;s,a)$. The CVaR of $C^d$ measures the severity of a high cost period drawn from stationary operation instead of the CVaR of an accumulated cost. Lemma~\ref{lemma_exchange_VaR} reduces the long-run criteria to the VaR and CVaR of this stationary cost. The proof of Lemma~\ref{lemma_exchange_VaR} is given in Appendix~\ref{sec_preliminary_app}.

\begin{lemma}\label{lemma_exchange_VaR}
Under Assumptions~\ref{A.cost}, \ref{A.function}, and~\ref{A.mixing}, for every stationary policy $d$ and every $s\in\mathcal{S}$, $\operatorname{VaR}^d(s) =\operatorname{VaR}[C^d]$, $\operatorname{CVaR}^d(s) =\operatorname{CVaR}[C^d]$.
\end{lemma}

Accordingly, we omit the initial state from $\operatorname{VaR}^d$ and $\operatorname{CVaR}^d$ hereinafter. The lemma shows that the long-run criterion considered in this paper is exactly the quantile and tail expectation of the stationary one stage cost.

Thus, the objective is equivalent to \(\min_{d\in\mathcal D}\operatorname{CVaR}^d\). We use the Bellman local optimality equation from \citet{xia2023risk} as the starting point. For a deterministic candidate $d^*(s)=\delta(a^*(s))$, the equations have the following form:
\begin{align}
a^*\left(s\right) &=\arg\min_{a\in\mathcal A(s)}\left\{\tilde c\left(\operatorname{VaR}^{d^*},s,a\right)+\E_{s'\sim p\left(\cdot\mid s,a\right)}\left[V^{d^*}\left(s'\right)\right]\right\},\nonumber\\
V^{d^*}\left(s\right)+\operatorname{CVaR}^{d^*} &=\min_{a\in\mathcal A(s)}\left\{\tilde c\left(\operatorname{VaR}^{d^*},s,a\right)+\E_{s'\sim p\left(\cdot\mid s,a\right)}\left[V^{d^*}\left(s'\right)\right]\right\},\qquad s\in\mathcal S.\label{Eq.AROE}
\end{align}
where $\tilde c(v,s,a):=v+\E[(C(s,a)-v)^+]/(1-\phi)$, $[x]^+:=\max\{x,0\}$, and $V^{d^*}$ is the corresponding relative value function. When the transition probabilities and cost distributions are known, the terms in~\eqref{Eq.AROE} can be evaluated and policy improvement can be carried out directly. In the sample-based setting, however, the VaR parameter is determined by the stationary distribution of the current policy. As the policy changes, both the observation law and the one stage cost change. Thus, the central problem is to learn the stationary VaR and the corresponding relative $Q$-function while updating the policy.

\subsection{Algorithm}\label{sec_alg}
This subsection develops the three recursive components of our algorithm. We first identify the fixed point associated with a fixed VaR parameter, then establish the VaR recursion, and finally describe the incremental policy update. 

\textbf{Q estimation.} Fix a reference state $s^{(0)}\in\mathcal S$. For a scalar $v$ and a $Q$-function $Q$, define
\begin{align*}
\left(\mathcal T_vQ\right)\left(s,a\right)
&=\tilde c\left(v,s,a\right)+\sum_{s'\in\mathcal S}p\left(s'\mid s,a\right)\min_{a'\in\mathcal A(s')}Q\left(s',a'\right)-\min_{a'\in\mathcal A(s^{(0)})}Q\left(s^{(0)},a'\right).
\end{align*}
Under regular assumptions, $\mathcal T_v$ has a unique fixed point, denoted by $Q(v)$. Its existence is proved in Lemma~\ref{lemma_Q_continuity_app}. For a stationary policy $d$, we use the shorthand $Q^d:=Q(\operatorname{VaR}^d)$. For a stationary policy $d$, define $\Psi_d(v):=\sum_{s\in\mathcal S}\mu^d(s)\sum_{a\in\mathcal A(s)}d(s,a)\tilde c(v,s,a)$, and write $h_v(s):=\min_{a\in\mathcal A(s)}Q(v)(s,a)$ and $g(v):=h_v(s^{(0)})$. The fixed point equation gives $Q(v)(s,a)+g(v)=\tilde c(v,s,a)+\sum_y p(y\mid s,a)h_v(y)$. Multiplying by $\mu^d(s)d(s,a)$ and summing yields
\begin{align}\label{Eq.threshold_control_bridge}
\Psi_d(v)-g(v)
&=\sum_{s\in\mathcal S}\mu^d(s)\sum_{a\in\mathcal A(s)}d(s,a)\left[Q(v)(s,a)-h_v(s)\right]\ge0,
\end{align}
where the equality uses $\sum_s\mu^d(s)\sum_a d(s,a)p(y\mid s,a)=\mu^d(y)$. Thus, we have
\begin{align}\label{Eq.threshold_control_value}
g(v)&=\min_{\widetilde d\in\mathcal D}\Psi_{\widetilde d}(v),
&g\left(\operatorname{VaR}^d\right)&\le\Psi_d\left(\operatorname{VaR}^d\right)=\operatorname{CVaR}^d.
\end{align}
If every stationary state probability is positive, equality in the inequality holds exactly when $d(s)$ is supported on minimizing actions of its own $Q^d$ at every state. 

For a policy $d^*$ satisfying the Bellman local optimality equation, the fixed point satisfies
\begin{align*}
a^*\left(s\right)
&=\arg\min_{a\in\mathcal A(s)}Q^{d^*}\left(s,a\right),\\
Q^{d^*}\left(s,a\right)+\operatorname{CVaR}^{d^*}
&=\tilde c\left(\operatorname{VaR}^{d^*},s,a\right)+\E_{s'\sim p\left(\cdot\mid s,a\right)}\left[\min_{a'\in\mathcal A(s')}Q^{d^*}\left(s',a'\right)\right],
\end{align*}
where $\operatorname{CVaR}^{d^*}=\min_{a'\in\mathcal A(s^{(0)})}Q^{d^*}(s^{(0)},a')$. Thus, we can use an SA procedure for the fixed point of $\mathcal T_v$ with an additional estimation of $\operatorname{VaR}^{d^*}$ to construct the Q recursion.

\textbf{VaR estimation.} For a continuous cumulative distribution function $F_X$, the quantile defined by the generalized inverse satisfies $F_X(v)=\phi$. If independent observations are attainable, a Robbins and Monro recursion for this equation is $v_{n+1}=v_n+\alpha_n(\phi-\I\{X_n\le v_n\})$. In our setting, an independent observation from the stationary cost distribution is not available after each policy update. We therefore use the current one stage cost and show that the recursion tracks the VaR of the current policy.

\textbf{Policy update.} For each state $s$, let $\mathcal D_{\epsilon_n}(s):=\{x\in\mathbb R^{|\mathcal A(s)|}:\sum_{a\in\mathcal A(s)}x(a)=1,\ x(a)\ge\epsilon_n\}$, and let $\Pi_n[\cdot]$ denote the statewise Euclidean projection onto $\mathcal D_{\epsilon_n}(s)$. We project $d_n(s)$ to $\mathcal D_{\epsilon_n}(s)$ at each iteration. If several actions minimize $Q_n(s,\cdot)$, the action with the smallest index is selected. The lower bound $\epsilon_n$ keeps every action available during learning and tends to zero. Replacing the current policy directly by a greedy policy can cause an abrupt change in the stationary cost distribution. We instead use an incremental update. The update seeks a policy whose selected actions minimize its own threshold control $Q$-function. 

The complete procedure is Algorithm \ref{GOODALGORITHM}. Let $v_n$, $Q_n$, and $d_n$ denote the VaR, $Q$-function, and policy estimates at iteration $n$. Only the visited Q component is projected, and all three estimates are updated from the same observation. Let $\widetilde C(v,s,a):=v+[C(s,a)-v]^+/(1-\phi)$, so that $\E[\widetilde C(v,s,a)]=\tilde c(v,s,a)$. Fix $Q_{\max}>0$ and let $\Pi_{Q_{\max}}(x):=\min\{Q_{\max},\max\{-Q_{\max},x\}\}$ be the scalar projection onto $[-Q_{\max},Q_{\max}]$.

\begin{algorithm}[htbp]
\caption{The Long-Run CVaR RL.}
\label{GOODALGORITHM}
\begin{enumerate}[leftmargin=*]
\item {\textbf{Input} projection radius $Q_{\max}>0$; initial state $s_0\in\mathcal S$ and estimates $v_0\in\mathbb R$, $d_0$, and $Q_0$ with $|Q_0(s,a)|\le Q_{\max}$ and $d_0(s)\in\mathcal D_{\epsilon_0}(s)$; learning rate sequences $\{\alpha_n\}$, $\{\beta_n(s,a):s\in\mathcal S,~a\in\mathcal A(s)\}$, and $\{\gamma_n\}$; exploration sequence $\{\epsilon_n\}$; a reference state $s^{(0)}\in\mathcal S$.}
\item \textbf{Initialize} the iteration counter $n\gets0$.
\item \textbf{Iterate} until a stopping rule is satisfied. Select an action $a_n\sim d_n(s_n)$. Observe the next state $s_{n+1}\sim p(\cdot\mid s_n,a_n)$ and the random cost $C(s_n,a_n)$.
{\allowdisplaybreaks\small
\begin{align}
v_{n+1}
&=v_n+\alpha_n\left(\phi-\I\left\{C\left(s_n,a_n\right)\le v_n\right\}\right),\label{Eq.v}\\
Q_{n+1}\left(s,a\right)
&=\left\{
\begin{aligned}
\Pi_{Q_{\max}}\left(\left(1-\beta_n(s,a)\right)Q_n(s,a)+\beta_n(s,a)\left[\widetilde C\left(v_n,s,a\right)\right.\right.\\[-1mm]
\left.\left.\qquad+\min_{a'\in\mathcal A(s_{n+1})}Q_n\left(s_{n+1},a'\right)-\min_{a'\in\mathcal A(s^{(0)})}Q_n\left(s^{(0)},a'\right)\right]\right),
&\quad s=s_n,~a=a_n,\\
Q_n\left(s,a\right),
&\quad\text{otherwise},
\end{aligned}
\right.\label{Eq.Q}\\
d_{n+1}\left(s\right)
&=\Pi_n\left[d_n\left(s\right)
+\gamma_n\left(\delta\left(\arg\min_{a\in\mathcal A(s)}Q_{n+1}\left(s,a\right)\right)-d_n\left(s\right)\right)\right],\qquad s\in\mathcal S,\label{Eq.d}\\
n&\gets n+1.\nonumber
\end{align}}
\end{enumerate}
\end{algorithm}

The three recursions are separate but coupled. The VaR update estimates the stationary quantile of the current policy, the Q update approximates the fixed point associated with the current VaR estimate, and the policy update moves gradually toward the current minimizing action. The relation $\gamma_n=o(\alpha_n)$ makes the policy nearly fixed on bounded VaR time intervals. Each Q component is studied on its own accumulated gain sequence, allowing a nonminimizing action to be sampled less frequently as exploration decreases. The projection $\Pi_{Q_{\max}}$ prevents the Q estimates from leaving a sufficiently large interval, while the desired fixed points remain in its interior.

\section{Convergence Analysis}\label{sec_converge}
This section establishes the convergence of Algorithm~\ref{GOODALGORITHM}. The main difficulties are that the observations are generated by a changing Markov chain and each Q component is updated only when its state and action pair is visited. We first introduce the notations and assumptions, then prove VaR and Q tracking, and finally study the local behavior of the policy recursion. The proofs are given in Appendix~\ref{sec_convergence_app}.

Let $(\bar\Omega,\mathscr F,\PP)$ be the probability space induced by the algorithm, where $\bar\Omega$ is the collection of all sample trajectories, $\mathscr F$ is a sigma field of subsets of $\bar\Omega$, and $\PP$ is a probability measure on $\mathscr F$. Let $\mathscr F_n$ contain the initial values and all observations completed before action $a_n$ is selected, together with $(s_n,v_n,d_n,Q_n)$, and let $\mathscr G_n:=\sigma(\mathscr F_n,a_n)$. Thus, $s_n$, $v_n$, $d_n$, and $Q_n$ are $\mathscr F_n$ measurable, whereas $a_n$ is $\mathscr G_n$ measurable. Write $C_n=C(s_n,a_n)$ for the cost observed at iteration $n$. For every Borel set $B$ and state $s'$, $\PP(C_n\in B,\ s_{n+1}=s'\mid\mathscr G_n) =\mathsf K(B,s'\mid s_n,a_n)$.
Taking the two marginals gives $\PP(C_n\le x\mid\mathscr G_n)=F(x;s_n,a_n)$ and $\PP(s_{n+1}=s'\mid\mathscr G_n)=p(s'\mid s_n,a_n)$. For two positive sequences $\{x_n\}$ and $\{y_n\}$, we write $x_n=o(y_n)$ when $x_n/y_n$ converges to zero and $x_n=O(y_n)$ when $x_n/y_n$ is bounded for all sufficiently large $n$. For policies, define $\|d-d'\|:=\max_{s\in\mathcal S}\|d(s)-d'(s)\|_1$, and for $Q$-functions define $\|Q\|_\infty:=\max_{s\in\mathcal S,~a\in\mathcal A(s)}|Q(s,a)|$. Let $N_n(s,a):=\sum_{k=0}^n\I\{s_k=s,a_k=a\}$ be the number of visits to $(s,a)$ by iteration $n$. Let $v_\phi(s,a):=F^{-1}(\phi;s,a)$. Since the state and action spaces are finite, we define $v_{\min}:=\min_{s\in\mathcal S,~a\in\mathcal A(s)}v_\phi(s,a)$ and $v_{\max}:=\max_{s\in\mathcal S,~a\in\mathcal A(s)}v_\phi(s,a)$, which are finite. We also fix deterministic numbers $\underline v<v_{\min}$, $\overline v>v_{\max}$, and a constant $\eta_v>0$, such that $F(\underline v;s,a)\le\phi-\eta_v$ and $F(\overline v;s,a)\ge\phi+\eta_v$ for every $(s,a)$. For a stationary policy $d$, let $\mu^d$ be the stationary state distribution and define the stationary cost distribution by $F^d(v):=\sum_{s\in\mathcal S}\mu^d(s)\sum_{a\in\mathcal A(s)}d(s,a)F(v;s,a)$.

We now state the regularity assumptions on the cost distributions, input sequences, and controlled Markov chain used in the analysis.

\begin{assumption}\label{A.function}
For every feasible $(s,a)$, $F(\cdot;s,a)$ has a continuous density $f(\cdot;s,a)$. There are constants $r_v>0$, $M_f<\infty$, and $\varepsilon_f>0$ such that $f(v;s,a)\le M_f$ for all $v\in\mathbb R$ and $f(v;s,a)\ge\varepsilon_f$ whenever $\underline v-r_v\le v\le\overline v+r_v$.
\end{assumption}

\begin{assumption}\label{A.stepsize}
The learning and exploration sequences satisfy the following conditions:
\begin{enumerate}[label=\upshape(\alph*), ref=\theassumption (\alph*)]
\item\label{A.alpha} $\alpha_n>0$ is nonincreasing, $\alpha_{n+1}/\alpha_n$ converges to one, $\sum_n\alpha_n=\infty$, and $\sum_n\alpha_n^2<\infty$.
\item\label{A.gamma} $0<\gamma_n<1$ is nonincreasing, $\gamma_{n+1}/\gamma_n$ converges to one, $\sum_n\gamma_n=\infty$, and $\sum_n\gamma_n^2<\infty$.
\item\label{A.ostep} $\gamma_n=o(\alpha_n)$.
\item\label{A.beta} There is a nonincreasing function $f_\beta$ such that $0<\beta_n(s,a)=f_\beta(N_n(s,a))<1$, $f_\beta(n+1)/f_\beta(n)$ converges to one, $\sum_nf_\beta(n)=\infty$, and $\sum_nf_\beta^2(n)<\infty$.
\item\label{A.epsilonstep} $0<\epsilon_n\le1/\max_{s\in\mathcal S}|\mathcal A(s)|$ is decreasing, $\epsilon_n$ converges to zero, $\sum_n\epsilon_n=\infty$, and $(\epsilon_{n-1}-\epsilon_n)/(\gamma_n\epsilon_n)$ converges to zero.
\end{enumerate}
\end{assumption}

\begin{assumption}\label{A.mixing}
There are an integer $m_p\ge1$, a constant $\delta_p>0$, and a probability vector $\nu$ with $\min_s\nu(s)>0$ such that, for every admissible policy sequence and every $n$, $\PP(s_{n+m_p}=s'\mid\mathscr F_n) \ge\delta_p\nu(s'),~ s'\in\mathcal S$ almost surely (a.s.).
\end{assumption}

\begin{remark}
Assumption~\ref{A.function} imposes an upper density bound on the whole real line and a positive lower bound only on the specified neighborhood containing the stationary policy quantiles. Gaussian, Student $t$ with more than two degrees of freedom, and logistic distributions satisfy the condition for any fixed choice of the interval. A bounded-support distribution satisfies Assumption \ref{A.function} provided that its density is continuous on \(\mathbb R\) and strictly positive throughout the specified interval. Since the interval is compact and the number of state–action pairs is finite, a common positive lower bound then exists. Similar conditions are used in SA-based quantile estimation \citep{Hu2022,hu2024quantile}.

Assumption~\ref{A.stepsize} is a standard condition used in SA \citep{kushner2003sa}. The relation $\gamma_n=o(\alpha_n)$ makes the policy nearly fixed when the VaR recursion is viewed on its own time scale. The condition on $\epsilon_n$ requires the exploration bound to decrease slowly relative to one policy step.

Assumption~\ref{A.mixing} is a common finite step minorization condition. It allows structural zeros in one step transitions. For each frozen stationary policy, the common lower bound yields uniform geometric mixing. Similar uniform mixing conditions are used in SA with controlled Markov noise and in average cost RL \citep{borkar2008stochastic,kushner2003sa,abounadi2001learning,yang2024relative}. Assumption~\ref{A.mixing} also controls the change in the stationary distribution as shown in Lemma \ref{lemma_stationary_lipschitz}.

\begin{lemma}\label{lemma_stationary_lipschitz}
Under Assumption~\ref{A.mixing}, every stationary policy $d$ has a unique stationary distribution $\mu^d$, and, for all stationary policies $d$ and $d'$, $\|\mu^d-\mu^{d'}\|_1 \le L_\mu\|d-d'\|$, where we define $L_\mu:=m_p/\delta_p$.
\end{lemma}
\end{remark}

\subsection{Asymptotic Behavior of the VaR Estimator}\label{sec_asymptotic_var}
The VaR recursion uses the current one stage cost rather than an independent observation from the stationary cost distribution. Its error consists of a martingale observation term and a transient term generated by the controlled Markov chain. A Poisson equation is used to control the transient term, leaving a stable root finding recursion whose root moves with the policy. Because $\gamma_n=o(\alpha_n)$, the VaR estimate can track this moving root. 

We begin with the following stability property needed for the VaR analysis. 

\begin{lemma}\label{lemma_v_bound}
Under Assumptions~\ref{A.cost} and \ref{A.alpha}, we have $\E[\sup_{n\ge0}v_n^2]<\infty$. Consequently, $\sup_n\E[v_n^2]<\infty$ and $\sup_n|v_n|<\infty$ a.s.
\end{lemma}

Theorem \ref{th_v} then gives the conclusion for the current policy $d_n$ without assuming that the policy sequence has already converged.

\begin{theorem}\label{th_v}
Under Assumptions~\ref{A.cost}--\ref{A.mixing}, $|v_n-\operatorname{VaR}^{d_n}|\to0$ a.s.
\end{theorem}

\subsection{Asymptotic Behavior of the $Q$-function Estimator}\label{sec_Qfunc}

The Q recursion is asynchronous, and a nonminimizing action may be visited much less frequently than the current minimizing action. Common assumptions may therefore be unsuitable under decreasing exploration. We instead require each component to accumulate an infinite amount of its own Q gain while the VaR estimate changes negligibly. 

For each scalar $v$, recall that $Q(v)$ is the fixed point of $\mathcal T_v$ and that $Q^d=Q(\operatorname{VaR}^d)$. Let $w(s,a)>0$ be fixed weights for all feasible state and action pairs, normalized so that $\min_{s\in\mathcal S,~a\in\mathcal A(s)}w(s,a)=1$, and write $\|Q\|_w:=\max_{s\in\mathcal S,~a\in\mathcal A(s)}|Q(s,a)|/w(s,a)$ and $w_{\max}:=\max_{s\in\mathcal S,~a\in\mathcal A(s)}w(s,a)$. We make the following assumptions.

\begin{assumption}\label{A.Qlearning}
We assume:
\begin{enumerate}[label=\upshape(\alph*), ref=\theassumption (\alph*)]
\item\label{A.Qcontraction} There is $\rho_Q\in(0,1)$ such that, for every $v\in[v_{\min}-r_v,v_{\max}+r_v]$ and all $Q$-functions $Q,Q'$, $\|\mathcal T_vQ-\mathcal T_vQ'\|_w \le\rho_Q\|Q-Q'\|_w$.
\item\label{A.Qprojection} The initial estimate satisfies $|Q_0(s,a)|\le Q_{\max}$ for every $(s,a)$, and
\begin{align*}
Q_{\max}
&>1+\frac{w_{\max}}{1-\rho_Q}
\frac{\displaystyle \max_{s\in\mathcal S,~a\in\mathcal A(s)}\E\left[|C(s,a)|\right]+(2-\phi)\max\left\{|v_{\min}-r_v|,|v_{\max}+r_v|\right\}}{1-\phi}.
\end{align*}
\item\label{A.Qgain} For every feasible $(s,a)$, $\sum_n\beta_n(s,a)\I\{s_n=s,a_n=a\}=\infty$ a.s., and for every $T>0$,
\begin{align*}
\sup\left\{\sum_{k=n}^{m}\alpha_k:m\ge n,~\sum_{k=n}^{m}\beta_k(s,a)\I\{s_k=s,a_k=a\}\le T\right\}
&\longrightarrow0\quad\text{a.s.}
\end{align*}
\end{enumerate}
\end{assumption}

\begin{remark}
Assumption~\ref{A.Qcontraction} allows different scales for the Q components. As shown in Proposition \ref{prop_weight_certificate}, it imposes a restriction on the reference state and on the transition probabilities.
\begin{proposition}\label{prop_weight_certificate}
If there are positive state weights $w(s)$ and a constant $\rho_Q\in(0,1)$ such that, for every $s\in\mathcal S$ and $a\in\mathcal A(s)$, $\sum_{y\ne s^{(0)}}p(y\mid s,a)w(y)+[1-p(s^{(0)}\mid s,a)]w(s^{(0)}) \le\rho_Qw(s)$, then Assumption~\ref{A.Qcontraction} holds with $w(s,a)=w(s)$. Under Assumption~\ref{A.mixing}, the existence of such a certificate, the existence of positive state and action weights satisfying~\ref{A.Qcontraction}, and the condition
\begin{align}\label{Eq.reference_occupancy_necessary}
\mu^d\left(s^{(0)}\right)&>\frac12,\qquad d\in\mathcal D,
\end{align}
are equivalent.
\end{proposition}

Assumption~\ref{A.Qprojection} imposes an explicit lower bound on the projection radius. In particular, since $\mathcal T_v0(s,a)=\tilde c(v,s,a)$ and $w(s,a)\ge1$, we have $\|\mathcal T_v0\|_w \le (\max_{s\in\mathcal S,~a\in\mathcal A(s)}\E[|C(s,a)|]+(2-\phi)|v|)/(1-\phi)$.
The contraction fixed point therefore satisfies $\|Q(v)\|_w\le\|\mathcal T_v0\|_w/(1-\rho_Q)$, and hence $\|Q(v)\|_\infty\le w_{\max}\|Q(v)\|_w$. The displayed choice of $Q_{\max}$ places all fixed points required by the proof at least one unit inside the projection interval. 

Assumption~\ref{A.Qgain} requires every Q component to receive an infinite total gain and imposes a uniform pathwise condition: whenever one component accumulates at most a fixed amount of Q gain after iteration $n$, the corresponding accumulated VaR gain converges to zero. As shown in Proposition~\ref{prop_gain_certificate}, under a polynomial specification, Assumption~\ref{A.Qgain} follows from the exploration lower bound and Assumption~\ref{A.mixing}.
\begin{proposition}\label{prop_gain_certificate}
If Assumption~\ref{A.mixing} holds, and $f_\beta(j)=\beta_0/(j+1)^{p_\beta}$, $\alpha_n=\alpha_0/(n+1)^{p_\alpha}$, $\gamma_n=\gamma_0/(n+1)^{p_\gamma}$, and $\epsilon_n=\epsilon_0/(n+1)^{p_\epsilon}$, where $\alpha_0>0$, $0<\beta_0<1$, $0<\gamma_0<1$, $0<\epsilon_0\le1/\max_{s\in\mathcal S}|\mathcal A(s)|$, and $1/2 <p_\beta<1,~ 0<p_\epsilon<1,~ p_\beta+p_\epsilon<p_\alpha<p_\gamma<1$, then Assumption~\ref{A.Qgain} holds.
\end{proposition}
\end{remark}

Finally, we show that the $Q$ estimator tracks its changing target, even when some state--action pairs are updated much less frequently than others. Theorem \ref{th_Q} states this result. The proof controls the effects of observation noise and changes in VaR.

\begin{theorem}\label{th_Q}
Under Assumptions~\ref{A.cost}--\ref{A.Qlearning}, $\|Q_n-Q^{d_n}\|_\infty\to0$ a.s.
\end{theorem}

\subsection{Local Stability of the Policy Recursion}\label{sec_main}
Theorems \ref{th_v}--\ref{th_Q} establish tracking of the two quantities used in the policy update, but they do not imply that a policy generated from an arbitrary initial condition enters the neighborhood of a locally optimal policy. We therefore study the local behavior of the policy recursion around a deterministic policy with a strict action gap. Continuity of the policy VaR and $Q(v)$ gives a neighborhood in which the same action remains the unique minimizer. Once a sufficiently late policy iterate enters this neighborhood, the tracking results keep the estimated minimizing actions correct and the policy recursion moves toward the deterministic policy. To address this issue, we first make the following assumption.

\begin{assumption}\label{A.policylimit}
There is a deterministic policy $d^*$ with $d^*(s)=\delta(a^*(s))$ such that\newline
$\Delta_Q :=\min_{\substack{s\in\mathcal S,~a\in\mathcal A(s)\\a\ne a^*(s)}}[Q^{d^*}(s,a)-Q^{d^*}(s,a^*(s))]>0$.
\end{assumption}

\begin{remark}
Assumption~\ref{A.policylimit} is a local strict action gap condition. Since $d\mapsto Q^d$ is continuous, it keeps the minimizing action unchanged in a sufficiently small policy neighborhood. States with only one feasible action do not enter the minimum in Assumption~\ref{A.policylimit}. If every state has only one action, the policy is fixed and the policy optimization conclusions are immediate, while this assumption concerns the remaining case. Under this assumption, Proposition \ref{prop_performance_difference} derives the corresponding performance difference bound from the fixed point equation, stationarity, the CVaR representation, and the Lipschitz property of the policy VaR.
\begin{proposition}\label{prop_performance_difference}
Under Assumptions~\ref{A.cost}, \ref{A.function},~\ref{A.mixing},~\ref{A.Qcontraction}, and \ref{A.policylimit}, there is a finite constant $C_{\rm p}$ such that, for every stationary policy $d$,
\begin{align*}
\operatorname{CVaR}^d-\operatorname{CVaR}^{d^*}
&\ge\sum_{s\in\mathcal S}\mu^d(s)\sum_{a\in\mathcal A(s)}d(s,a)\left[Q^{d^*}(s,a)-\min_{b\in\mathcal A(s)}Q^{d^*}(s,b)\right]-C_{\rm p}\left\|d-d^*\right\|^2.
\end{align*}
Consequently, $d^*$ is locally optimal among stationary randomized policies.
\end{proposition}
\end{remark}

Theorem \ref{th_opt} states the local stability result. It shows that every sufficiently late entry into a smaller strict gap neighborhood is preserved and leads to convergence a.s. 

\begin{theorem}\label{th_opt}
Under Assumptions~\ref{A.cost}--\ref{A.policylimit}, there exists a constant $r>0$ independent of the sample path such that, almost surely, there is a finite $N_0(\omega)$ such that, for every $N\ge N_0(\omega)$, if $\|d_N(\omega)-d^*\|\le r/4$, then $d_n(\omega)$ converges to $d^*$. Moreover, $d^*$ is locally optimal among stationary randomized policies.
\end{theorem}

Once a sufficiently late iterate enters the strict action gap neighborhood, the policy recursion converges to the corresponding locally optimal policy. The entry condition is separate from the tracking results. In an implementation, a policy obtained from an approximate model, an offline preliminary calculation, or several initial policies may be used to provide a local starting point. 

Theorem~\ref{prop_post_rate} gives the order $O(n^{-p_\epsilon})$ for the policy recursion after the minimizing actions have been identified. It does not control the random action identification time and is not a finite sample rate for the complete learning procedure.

\begin{theorem}\label{prop_post_rate}
Let $\gamma_n=\gamma_0/(n+1)^{p_\gamma}$ and $\epsilon_n=\epsilon_0/(n+1)^{p_\epsilon}$, where $\gamma_0,\epsilon_0>0$ and $0<p_\epsilon<p_\gamma<1$. Under the conditions of Theorem~\ref{th_opt}, after the minimizing actions have been identified, for every state $s$, $\|d_n(s)-d^*(s)\|_1 =O(\epsilon_n)$ a.s.
\end{theorem}

\subsection{Extension to Mean CVaR Optimization}\label{sec.mean-cvar}
The same analysis can be used when a mean cost term is added to CVaR, but the objective function and its minimizing actions will change. For a fixed $\lambda \ge 0$, we consider the mean--CVaR optimization problem $\min_{d \in \mathcal{D}} J_\lambda(d)$, where $J_\lambda(d) := \operatorname{CVaR}^{d} + \lambda \bar{C}^{d}$ and $\bar{C}^{d} := \mathbb{E}[C^{d}]$ denotes the expected stationary one-stage cost under policy $d$. We therefore identify the modified fixed point before transferring the tracking and local stability arguments. The transition condition is unchanged, whereas the projection bound and the local candidate must be checked for the new objective.

In this case, the augmented Q sample is $\widetilde C(v_n,s,a)+\lambda C(s,a)$. Its expected operator and target are $(\mathcal T_{v,\lambda}Q)(s,a)=(\mathcal T_vQ)(s,a)+\lambda\E[C(s,a)],~Q_\lambda^d:=Q_\lambda(\operatorname{VaR}^d)$, where $Q_\lambda(v)$ is the unique fixed point of $\mathcal T_{v,\lambda}$. The extra expected cost is independent of $Q$, so the contraction coefficient and weights are unchanged. It is also independent of $v$, so the threshold Lipschitz bound remains valid. For the projection bound, write $M_1:=\max_{s\in\mathcal S,\,a\in\mathcal A(s)}\E[|C(s,a)|]$ and $V:=\max\{|v_{\min}-r_v|,|v_{\max}+r_v|\}$. A sufficient radius is
\begin{align}\label{Eq.mean_projection_bound}
Q_{\max}&>1+\frac{w_{\max}\left(\left[1+\lambda(1-\phi)\right]M_1+(2-\phi)V\right)}{(1-\rho_Q)(1-\phi)}.
\end{align}
This is Assumption~\ref{A.Qprojection} with an additional term for the mean cost. The initialization still satisfies $|Q_0(s,a)|\le Q_{\max}$.
For local convergence, let $d_\lambda^*(s)=\delta(a_\lambda^*(s))$ be a deterministic candidate and require
\begin{align}\label{Eq.mean_strict_gap}
\Delta_{Q,\lambda}
&:=\min_{\substack{s\in\mathcal S,~a\in\mathcal A(s)\\a\ne a_\lambda^*(s)}}
\left[Q_\lambda^{d_\lambda^*}(s,a)-Q_\lambda^{d_\lambda^*}(s,a_\lambda^*(s))\right]>0.
\end{align}

The following result states the conclusions with the new fixed point and its own local entry condition.

\begin{corollary}\label{cor_mean_cvar}
Under Assumptions~\ref{A.cost}--\ref{A.mixing}, \ref{A.Qcontraction}, \ref{A.Qgain}, and condition~\eqref{Eq.mean_projection_bound}, we have $|v_n-\operatorname{VaR}^{d_n}|\to0$ and $\|Q_n-Q_\lambda^{d_n}\|_w\to0$ a.s.

If $d_\lambda^*$ satisfies condition~\eqref{Eq.mean_strict_gap}, it is a strict local minimizer of $J_\lambda$ among stationary randomized policies. Moreover, there exists a constant $r_\lambda>0$ independent of the sample path such that, almost surely, there is a finite $N_{0,\lambda}(\omega)$ such that, for every $N\ge N_{0,\lambda}(\omega)$, if $\|d_N(\omega)-d_\lambda^*\|\le r_\lambda/4$, then $d_n(\omega)$ converges to $d_\lambda^*$.
\end{corollary}

\section{Numerical Experiments}\label{sec_experi}

In this section, we first give a small model for illustration in Section \ref{sec_theory_example}. We then compare Algorithm~\ref{GOODALGORITHM}, denoted by CRL, with RED on a machine replacement problem under two cost distributions in Section \ref{sec_machine}, and consider a renewable energy storage problem with discrete costs in Section \ref{sec_renewable}. The small model illustrates the theoretical conditions, while the two applications illustrate performance beyond some of these conditions. The main comparison concerns methods that optimize CVaR. Section~\ref{sec_numerical_details} of the Supplementary material gives the implementation details and compares policies obtained by minimizing CVaR, mean cost, and their weighted sum.

\subsection{An Illustrative Model}\label{sec_theory_example}
We first consider a two state and two action model with reference state $s^{(1)}$. Both actions are available at both states. Table~\ref{tab:covered_model} gives the transition probabilities and the means of the Gaussian cost distributions. Every cost has standard deviation $0.25$, and the CVaR level is $\phi=0.9$.

\begin{table}[H]
\centering
\caption{Transition probabilities and mean costs.}
\label{tab:covered_model}
\renewcommand{\arraystretch}{1.0}
\begin{tabular}{ccccc}
\toprule
State & Action & $p(s^{(1)}\mid s,a)$ & $p(s^{(2)}\mid s,a)$ & Mean cost \\
\midrule
$s^{(1)}$ & $a^{(1)}$ & 0.80 & 0.20 & 0.00 \\
$s^{(1)}$ & $a^{(2)}$ & 0.85 & 0.15 & 1.00 \\
$s^{(2)}$ & $a^{(1)}$ & 0.75 & 0.25 & 2.00 \\
$s^{(2)}$ & $a^{(2)}$ & 0.70 & 0.30 & 0.50 \\
\bottomrule
\end{tabular}
\end{table}

The model satisfies Assumption~\ref{A.mixing} with $m_p=1$, $\nu=(0.5,0.5)$, and $\delta_p=0.3$. On the finite state space, the same condition also gives the required Lipschitz continuity of the stationary distribution. For the Q contraction, we take $w(s^{(1)})=w(s^{(2)})=1$. The left hand side of the certificate in Proposition~\ref{prop_weight_certificate} is $2[1-p(s^{(1)}\mid s,a)]$ and has maximum value $0.60$. Thus, Assumption~\ref{A.Qcontraction} holds with $\rho_Q=0.60$. For the bound in Assumption~\ref{A.Qprojection}, take $r_v=1$, $M_1=2[2\Phi(8)-1]+0.5\varphi(8)$, and $V=3+0.25\Phi^{-1}(0.9)$, where $\Phi$ and $\varphi$ are the standard normal distribution and density functions. Here $M_1=\max_{s,a}\E|C(s,a)|$ and $V\approx3.3204$. The resulting bound is $1+25(M_1+1.1V)\approx142.3107$, and we set $Q_{\max}=150$.

We use the following input sequences: $\alpha_n = 0.2/(n+1)^{0.75}$, $\beta_n(s,a) = 0.5/(N_n(s,a)+1)^{0.55}$, $\gamma_n = 0.5/(n+1)^{0.90}$, $\epsilon_n = 0.05/(n+1)^{0.10}$. These powers satisfy $0.55+0.10<0.75<0.90$, so Proposition~\ref{prop_gain_certificate} applies. Enumerating the four deterministic policies gives $d^*(s^{(1)})=\delta(a^{(1)})$ and $d^*(s^{(2)})=\delta(a^{(2)})$. This policy has stationary distribution $(7/9,2/9)$, VaR $0.5595$, CVaR $0.7298$, and Q gaps $4.3931$ and $13.6370$. Hence, Assumption~\ref{A.policylimit} is also satisfied.

\begin{figure}[htbp]
\centering
\includegraphics[width=\linewidth, trim=0cm 0cm 0cm 1.4cm, clip]{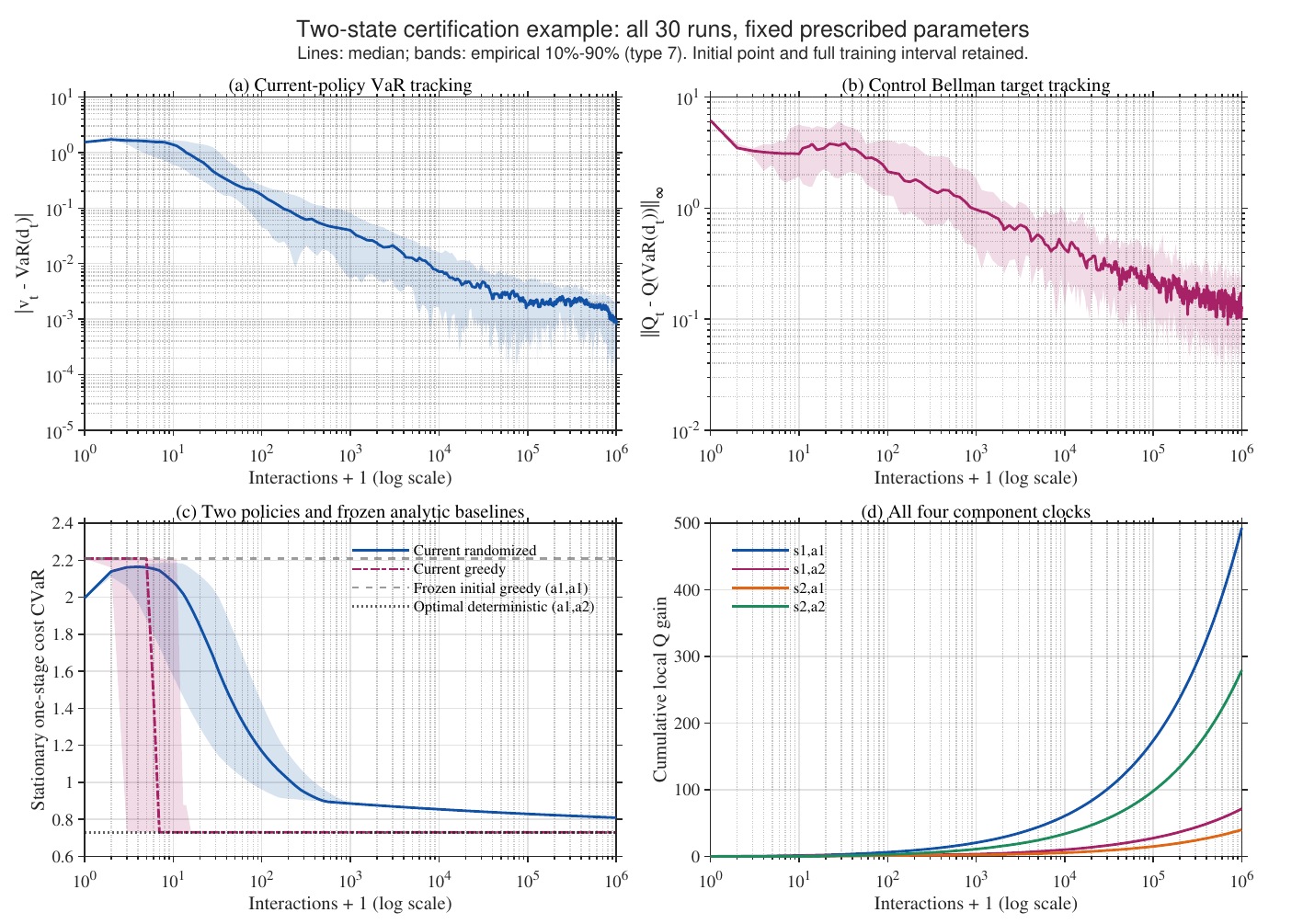}
\caption{Tracking and policy performance.}
\label{fig:covered_learning}
\end{figure}

We perform 30 independent replications of $10^6$ interactions, with $s_0=s^{(1)}$, $v_0=0$, $Q_0=0$, and a uniform initial policy. Conditional on each state and action, costs and next states are sampled independently. Visit counts include the current observation when evaluating $\beta_n(s,a)$. The final current randomized policies have mean VaR $0.5844$ and mean CVaR $0.8094$. The initial greedy policy selects $(a^{(1)},a^{(1)})$. Holding this policy fixed gives CVaR $2.2096$. All 30 final greedy policies select $(a^{(1)},a^{(2)})$, with CVaR $0.7298$. The remaining exploration explains the difference from the final greedy policy. Figure~\ref{fig:covered_learning} shows the errors $|v_n-\operatorname{VaR}^{d_n}|$ and $\|Q_n-Q^{d_n}\|_\infty$, policy CVaRs, and accumulated Q gains. Here $Q^{d_n}$ is the fixed threshold control solution defined in Section~\ref{sec_alg}, computed only for evaluation after training. The final mean errors are $0.0010$ and $0.1395$, respectively. No Q projection occurred in any run. Curves show medians and empirical $10\%$ to $90\%$ quantiles over all 30 runs.

This example shows that our Q-tracking result allows some actions to be sampled increasingly less often than others. On sample paths where the minimizing actions are eventually identified, Theorem~\ref{prop_post_rate} gives $d_n(s,a)=O(\epsilon_n)$ for each nonminimizing action, while the probability of selecting the minimizing action tends to one. In contrast, this sampling pattern does not satisfy the bounded update-frequency ratios required by \citet[Theorem C.1.1]{rojas2025burning}.

\FloatBarrier
\subsection{Machine Replacement}\label{sec_machine}
We next consider a machine replacement problem with state space $\mathcal{S}=\{s^{(1)},\ldots,s^{(6)}\}$, where the state represents the accumulated utilization of the machine. At each time epoch, the decision is to retain the current machine, denoted by $a^{(1)}=0$, or replace it with a new one, denoted by $a^{(2)}=1$. Replacement is the only feasible action at state $s^{(6)}$. Tables~\ref{tab:kernel1} and~\ref{tab:kernel2} give the transition matrices. Conditional on a feasible state and action pair, the cost follows a Gaussian distribution with standard deviation $0.5$ and mean given in Table~\ref{tab:mean_cost}. For the Student $t$ case, we take $C(s,a)=m(s,a)+0.5\sqrt{3/5}\,Z$, where $Z\sim t_5$ and $m(s,a)$ is the corresponding entry of Table~\ref{tab:mean_cost}. Thus both cost models have the same conditional means and standard deviation $0.5$. The objective is to minimize the long-run CVaR$_{0.9}$ of the one stage cost. The state dependent action restriction is represented by $\mathcal A(s)$. The always-replace policy has maximum stationary state probability $0.4960<1/2$, so Proposition~\ref{prop_weight_certificate} rules out Assumption~\ref{A.Qcontraction} for every reference state. The numerical powers below also fail the sufficient condition in Proposition~\ref{prop_gain_certificate}. This experiment is an empirical robustness study rather than a direct verification of Theorem~\ref{th_Q}.

\begin{table}[htbp]
\centering
\caption{Transition matrix for $a^{(1)}=0$.}
\renewcommand{\arraystretch}{1.0} 
\begin{tabular*}{\linewidth}{@{\extracolsep{\fill}}c*{6}{c}@{}}
\toprule
$a^{(1)}$ & $s^{(1)}$ & $s^{(2)}$ & $s^{(3)}$ & $s^{(4)}$ & $s^{(5)}$ & $s^{(6)}$ \\
\midrule
$s^{(1)}$ & 0.496 & 0.254 & 0.131 & 0.067 & 0.034 & 0.018 \\
$s^{(2)}$ & 0.000 & 0.505 & 0.259 & 0.133 & 0.068 & 0.035 \\
$s^{(3)}$ & 0.000 & 0.000 & 0.523 & 0.268 & 0.138 & 0.071 \\
$s^{(4)}$ & 0.000 & 0.000 & 0.000 & 0.563 & 0.289 & 0.148 \\
$s^{(5)}$ & 0.000 & 0.000 & 0.000 & 0.000 & 0.661 & 0.339 \\
$s^{(6)}$ & - & - & - & - & - & - \\
\bottomrule
\end{tabular*}
\label{tab:kernel1}
\end{table}
\begin{table}[htbp]
\centering
\caption{Transition matrix for $a^{(2)}=1$.}
\renewcommand{\arraystretch}{1.0} 
\begin{tabular*}{\linewidth}{@{\extracolsep{\fill}}c*{6}{c}@{}}
\toprule
$a^{(2)}$ & $s^{(1)}$ & $s^{(2)}$ & $s^{(3)}$ & $s^{(4)}$ & $s^{(5)}$ & $s^{(6)}$ \\
\midrule
$s^{(1)}$ & 0.496 & 0.254 & 0.131 & 0.067 & 0.034 & 0.018 \\
$s^{(2)}$ & 0.496 & 0.254 & 0.131 & 0.067 & 0.034 & 0.018 \\
$s^{(3)}$ & 0.496 & 0.254 & 0.131 & 0.067 & 0.034 & 0.018 \\
$s^{(4)}$ & 0.496 & 0.254 & 0.131 & 0.067 & 0.034 & 0.018 \\
$s^{(5)}$ & 0.496 & 0.254 & 0.131 & 0.067 & 0.034 & 0.018 \\
$s^{(6)}$ & 0.496 & 0.254 & 0.131 & 0.067 & 0.034 & 0.018 \\
\bottomrule
\end{tabular*}
\label{tab:kernel2}
\end{table}
\begin{table}[htbp]
\centering
\caption{Mean costs for each state-action pair.}
\renewcommand{\arraystretch}{1.0} 
\begin{tabular*}{\linewidth}{@{\extracolsep{\fill}}c*{6}{c}@{}}
\toprule
 & $s^{(1)}$ & $s^{(2)}$ & $s^{(3)}$ & $s^{(4)}$ & $s^{(5)}$ & $s^{(6)}$ \\
\midrule
$a^{(1)}$ & 0 & 3 & 6 & 9 & 12 & 15 \\
$a^{(2)}$ & 15 & 15 & 15 & 15 & 15 & 15 \\
\bottomrule
\end{tabular*}
\label{tab:mean_cost}
\end{table}

We compare CRL with the tabular RED CVaR Q learning method of \citet{rojas2025burning}. For RED, we use the reward $R=-C$, lower tail probability $0.1$, total random exploration probability $0.1$, Q step size $0.02$, and VaR and CVaR step sizes $0.002$. CRL uses $\alpha_n=10/(n+1)^{0.9}$, $\beta_n(s,a)=1/(N_n(s,a)+1)^{0.8}$, $\gamma_n=1/(n+1)^{0.99}$, and $\epsilon_n=1/[2(n+1)^{0.999}]$, with reference state $s^{(1)}$. Each algorithm is independently repeated 30 times with a budget of $N=10^6$ interactions, including 1000 warm up interactions. For evaluation, we enumerate all 32 feasible deterministic policies and compute their stationary cost distributions. Let $\bar d$ denote the policy with the smallest enumerated CVaR. It retains the machine in states $s^{(1)},\ldots,s^{(5)}$ and replaces it in state $s^{(6)}$ under both cost models. The zero initial Q tables and tie rule already select $\bar d$, so the saved paths describe departures from and returns to that benchmark.

\begin{figure}[htbp]
\centering
\includegraphics[width=\linewidth]{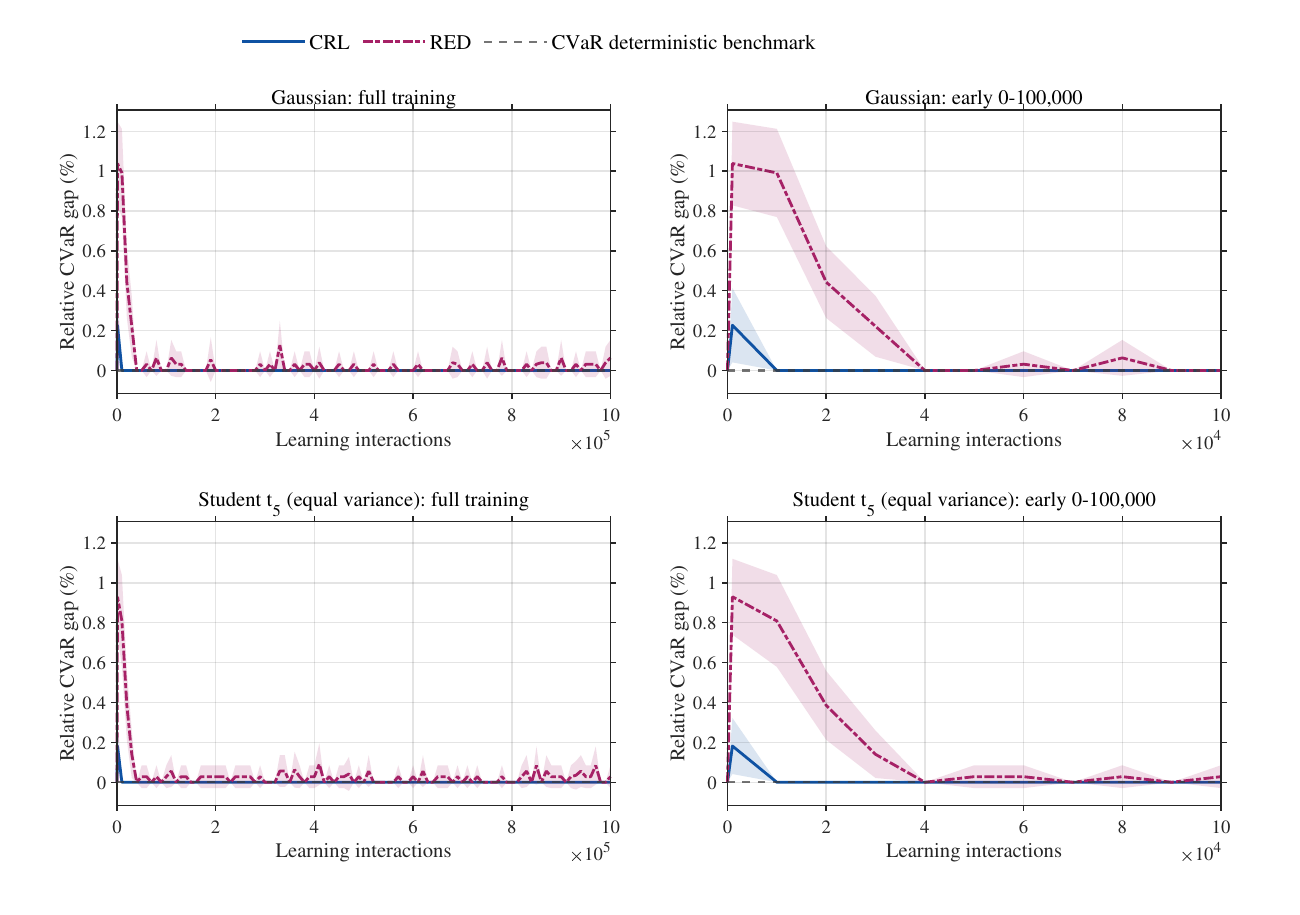}
\caption{Relative CVaR gaps of greedy policies.}
\label{fig: gap}
\end{figure}

The results are shown in Table~\ref{tab:performance} and Figure~\ref{fig: gap}. Parentheses in Table~\ref{tab:performance} give 95\% Student intervals based on the 30 replications; the deterministic benchmark is evaluated directly. The reported performance is the stationary CVaR of each saved deterministic greedy policy, evaluated with the same model equations as the benchmark. The stationary probabilities determine the cost mixture, its quantile is found from the mixture distribution, and its CVaR is computed from the positive part representation. To compare the saved paths, we use the gap $(\operatorname{CVaR}^{d}/\operatorname{CVaR}^{\bar d}-1)$, expressed as a percentage, for an evaluated policy $d$. The average gap in Table~\ref{tab:performance} is the trapezoidal integral through all 102 saved points, at iterations $0$, $1000$, and $10000,20000,\ldots,10^6$, divided by $N$. The number of observed policy changes compares adjacent saved action vectors. Figure~\ref{fig: gap} shows the complete budget and the same early window for each distribution. Left panels show the full budget and right panels show the common early window. Curves average all 30 runs, and shaded bands are pointwise $95\%$ Student intervals. Although the final values are close, CRL has a smaller average gap and fewer observed policy changes in both distributions under the stated configurations. The paired intervals in Table~\ref{tab:paired_supp} do not establish a final CVaR advantage.

\begin{table}[!t]
\centering
\caption{Greedy policy performance and saved paths.}
\label{tab:performance}
\small
\setlength{\tabcolsep}{3pt}
\renewcommand{\arraystretch}{1.15}
\begin{tabular}{@{}lccccc@{}}
\toprule
Method & Final CVaR & Mean cost & Average gap (\%) & Observed changes & Batch time (s) \\
\midrule
\multicolumn{6}{@{}l}{\textbf{Gaussian costs}} \\
Benchmark & 15.2239 & 8.1252 & 0 & --- & --- \\
CRL & \shortstack{15.2239 \\ {\footnotesize (15.2239, 15.2239)}} & \shortstack{8.1252 \\ {\footnotesize (8.1252, 8.1252)}} & \shortstack{0.00114 \\ {\footnotesize (0.00020, 0.00207)}} & \shortstack{0.40 \\ {\footnotesize (0.10, 0.70)}} & 26.309 \\
RED & \shortstack{15.2336 \\ {\footnotesize (15.2199, 15.2473)}} & \shortstack{8.0318 \\ {\footnotesize (7.8991, 8.1645)}} & \shortstack{0.03492 \\ {\footnotesize (0.02876, 0.04107)}} & \shortstack{4.80 \\ {\footnotesize (3.88, 5.72)}} & 15.827 \\
\midrule
\multicolumn{6}{@{}l}{\textbf{Student $t_5$ costs with matched conditional variance}} \\
Benchmark & 15.2120 & 8.1252 & 0 & --- & --- \\
CRL & \shortstack{15.2120 \\ {\footnotesize (15.2120, 15.2120)}} & \shortstack{8.1252 \\ {\footnotesize (8.1252, 8.1252)}} & \shortstack{0.00091 \\ {\footnotesize (0.00021, 0.00161)}} & \shortstack{0.40 \\ {\footnotesize (0.10, 0.70)}} & 27.956 \\
RED & \shortstack{15.2162 \\ {\footnotesize (15.2075, 15.2249)}} & \shortstack{8.0785 \\ {\footnotesize (7.9830, 8.1740)}} & \shortstack{0.03741 \\ {\footnotesize (0.03101, 0.04380)}} & \shortstack{6.67 \\ {\footnotesize (5.39, 7.94)}} & 18.722 \\
\bottomrule
\end{tabular}
\end{table}

Section~\ref{sec_numerical_details} gives further details of the experimental settings and additional comparisons between CRL and RED. CRL attains the best enumerated deterministic CVaR in all 30 final greedy policies under each cost distribution, while RED attains the same benchmark in 28 Gaussian runs and 29 Student $t_5$ runs. Table~\ref{tab:fixed_budget_supp} in the Supplementary material reports the relative CVaR gaps at iterations $10^4$, $5\times10^4$, $10^5$, $2\times10^5$, $5\times10^5$, and $10^6$. Under both cost distributions, the reported gaps for CRL are zero at these iterations, whereas those for RED are positive at some of them. In an additional comparison using the configurations in Table~\ref{tab:recorded_machine_inputs} and three different Q initializations, both methods attain the optimal deterministic CVaR in all 30 runs for each initialization and cost distribution.

\FloatBarrier
\subsection{Renewable Energy Storage System Scheduling}\label{sec_renewable}

We next consider a renewable energy storage scheduling problem. The state $s_n=B_n$ is the storage level, with $B_n\in\{0.4,1.0,1.6,2.2,2.8,3.4\}$. Renewable generation $G_n$ and demand $D_n$ are independent random variables with distributions given in Tables~\ref{tab: new energy} and~\ref{tab: demand}. The action $a_n\in\{-2.4,-1.2,0.6,1.2\}$ represents charging or discharging, where $a_n>0$ denotes discharging and $a_n<0$ denotes charging. Feasibility requires $-C_{min}\le a_n\le C_{max}$ and $B_n-B_{max}\le a_n\le B_n-B_{min}$. Thus, $C_{min}$ bounds the charging magnitude, $C_{max}$ bounds the discharging level, and the storage dynamics are $B_{n+1}=B_n-a_n$. This is a synthetic storage example with a prescribed action set that excludes idling. The action is selected from the current storage level before the generation and demand for that period are observed.

Let $W_n=D_n-G_n-a_n$ be the power exchanged with the backbone network, where $W_n>0$ indicates a shortage and $W_n<0$ a surplus. The one stage cost is $C(s_n,a_n)=p_b[W_n]^+-p_s[W_n]^-+c|a_n|+h(B_n-a_n)$, where $[x]^-=-\min\{x,0\}$, $p_b$ and $p_s$ are the buying and selling prices, $c$ is the battery throughput cost, and $h$ is the holding cost. The term $c|a_n|$ assigns a nonnegative utilization cost to both charging and discharging. Table~\ref{tab:prapmeter} gives the remaining parameters, and the objective is to minimize the CVaR$_{0.9}$ of the one stage cost under the long-run occupation distribution from $B_0=0.4$. The state dependent action sets are covered by $\mathcal A(s)$. The discrete costs and periodic policies fall outside the assumptions of Lemma~\ref{lemma_exchange_VaR}. Instead, we use the occupation distribution criterion as an extension. The example is treated as an empirical robustness study.

\begin{table}[htbp]
\centering
\caption{Distribution of $G_n$.}
\begin{tabular}{ccccccc}
\toprule
\textbf{Value} & 0.0 & 0.6 & 1.2 & 1.8 & 2.4 & 3.0 \\
\midrule
\textbf{Probability} & 0.10 & 0.30 & 0.20 & 0.10 & 0.15 & 0.15 \\
\bottomrule
\end{tabular}
\label{tab: new energy}
\end{table}

\begin{table}[htbp]
\centering
\caption{Distribution of $D_n$.}
\begin{tabular}{ccccccc}
\toprule
\textbf{Value} & 0.6 & 1.2 & 1.8 & 2.4 & 3.0 & 3.6 \\
\midrule
\textbf{Probability} & 0.05 & 0.25 & 0.15 & 0.25 & 0.2 & 0.1 \\
\bottomrule
\end{tabular}
\label{tab: demand}
\end{table}

\begin{table}[htbp]
\centering
\caption{Problem parameters.}
\renewcommand{\arraystretch}{1.0} 
\begin{tabular*}{\linewidth}{@{\extracolsep{\fill}}*{8}{c}@{}}
\toprule
\textbf{$C_{min}$} & \textbf{$C_{max}$} & \textbf{$B_{min}$} & \textbf{$B_{max}$} & \textbf{$p_b$} & \textbf{$p_s$} & \textbf{$c$} & \textbf{$h$} \\
\midrule
2.4 & 1.2 & 0.4 & 3.4 & 3 & 1.5 & 4 & 2 \\
\bottomrule
\end{tabular*}
\label{tab:prapmeter}
\end{table}

For this application, CRL is run for $6\times10^5$ interactions in 30 independent replications. The warm up period is included in the budget, and the input sequences are the machine replacement sequences except for $\epsilon_n=1/[4(n+1)^{0.999}]$. The Q projection bound is $Q_{\max}=500$. We compare warm up lengths $W$ of $2000$, $5000$, and $10000$.
The initial storage level for training and evaluation and the reference state are all $0.4$. During warm up, Q, VaR, and the visit counts continue to update. The first policy update occurs after observation $W+1$, so the first $W+1$ action samples are uniform.

To obtain a reference value, we evaluate all 216 feasible deterministic policies from $B_0=0.4$. One optimal policy makes the storage level repeat the sequence $0.4,1.6,1.0$. The cost distribution is obtained by assigning equal weight to the three positions in this sequence. This policy has VaR $17$, CVaR $18.41$, and mean cost $7.9180$, and minimizes both CVaR and mean cost among the policies considered. Here CVaR is computed after mixing, while averaging the three position-specific CVaRs instead gives $13.66$. Policies that differ at transient states may nevertheless enter this same cycle.

For each completed run, we select the action with the smallest final Q value at each storage level and evaluate the resulting policy. Table~\ref{tab: energy converge} gives the final CVaR values, averaged over the 30 runs, for both the greedy and current randomized policies. Each policy is evaluated from $B_0=0.4$ using its long-run state and action frequencies before the results are averaged. A run is successful if its final greedy CVaR is within $1\%$ of $18.41$. The reported intervals are Wilson $95\%$ intervals for the success probability.

\begin{table}[htbp]
\centering\small
\caption{CRL results for the storage example.}\label{tab: energy converge}
\setlength{\tabcolsep}{8pt}
\renewcommand{\arraystretch}{1.12}
\begin{tabular}{@{}rcccc@{}}\toprule
Warm up length & Greedy CVaR & Current CVaR & Successful runs & Wilson interval($\%$) \\\midrule
2000 & 18.4335 & 18.4751 & 29/30 & (83.3, 99.4) \\
5000 & 18.4100 & 18.5211 & 30/30 & (88.6, 100.0) \\
10000 & 18.4100 & 18.6341 & 30/30 & (88.6, 100.0) \\
\bottomrule\end{tabular}
\end{table}

With $5000$ or $10000$ warm up interactions, all 30 final greedy policies attain the benchmark values. With $2000$, 29 runs attain the benchmark and one has CVaR $19.1150$. Thus, CRL recovers an optimal deterministic policy in every run at the two longer warm up lengths. At $W=5000$, the final randomized policies have average CVaR $18.5211$ and mean cost $7.9549$, compared with $18.4100$ and $7.9180$ for the greedy policies. The randomized policies still select other actions, which explains their higher CVaR.

\FloatBarrier
\section{Conclusions}\label{sec_conclude}
This paper develops an RL method for finite state and finite action Markov decision processes under a stationary one stage cost CVaR criterion. The criterion measures the upper tail of a typical period in steady state and is different from the CVaR of a cumulative return. The method uses a single sample trajectory and separate but coupled recursions for the VaR, $Q$-function, and policy estimates. 
We first establish tracking of the policy dependent VaR and then use a component update argument to establish tracking of the \(Q\)-function under decreasing exploration.
The main theorem is a local stability result for the policy recursion. On sample paths for which a sufficiently late iterate enters a strict action gap neighborhood, the policy remains in the neighborhood and converges to a deterministic policy. After the minimizing actions have been identified, we also obtain the order $O(\epsilon_n)$ for the policy recursion.

\FloatBarrier
\section*{Supplementary material}
Supplementary material: Implementation details and additional numerical comparisons for the machine replacement experiment in Section~\ref{sec_machine}, including results for alternative objectives and final policies (Section~\ref{sec_numerical_details}).

\bibliographystyle{informs2014}
\bibliography{references}

@article{rockafellar2000optimization,
  title={Optimization of conditional value-at-risk},
  author={Rockafellar, R Tyrrell and Uryasev, Stanislav},
  journal={Journal of Risk},
  volume={2},
  pages={21--42},
  year={2000},
  publisher={Citeseer}
}

@article{hong2009simulating,
  title={Simulating sensitivities of conditional value at risk},
  author={Hong, L Jeff and Liu, Guangwu},
  journal={Management Science},
  volume={55},
  number={2},
  pages={281--293},
  year={2009},
  publisher={INFORMS}
}

@book{kushner2003sa,
	author={Kushner, Harold J. and  Yin, G. George},
	Pages = {166--169},
	Title = {Stochastic Approximation and Recursive Algorithms and Applications},
	publisher={Springer},
	Year = {2003}
}

@article{Hu2022,
Author = {Jiaqiao Hu and Yijie Peng and Gongbo Zhang and Qi Zhang},
Title = {A Stochastic Approximation Method for Simulation-Based Quantile Optimization},
Journal = {INFORMS Journal on Computing},
Volume = {34},
Number = {6},
Year = {2022},
Pages = {2889--2907},
DOI = {10.1287/ijoc.2022.1214}}

@inproceedings{tamar2015optimizing,
  title={Optimizing the {CVaR} via sampling},
  author={Tamar, Aviv and Glassner, Yonatan and Mannor, Shie},
  booktitle={Proceedings of the AAAI Conference on Artificial Intelligence},
  volume={29},
  year={2015}
}

@article{hu2024quantile,
  title={Quantile Optimization via Multiple-Timescale Local Search for Black-Box Functions},
  author={Hu, Jiaqiao and Song, Meichen and Fu, Michael C},
  journal={Operations Research},
  year={2025},
  volume={73},
  number={3},
  pages={1535--1557},
  publisher={INFORMS}
}

@article{xia2023risk,
  title={Risk-sensitive {Markov} decision processes with long-run {CVaR} criterion},
  author={Xia, Li and Zhang, Luyao and Glynn, Peter W},
  journal={Production and Operations Management},
  volume={32},
  number={12},
  pages={4049--4067},
  year={2023},
  publisher={Wiley Online Library}
}

@article{abounadi2001learning,
  title={Learning algorithms for {Markov} decision processes with average cost},
  author={Abounadi, Jinane and Bertsekas, Dimitri and Borkar, Vivek S},
  journal={SIAM Journal on Control and Optimization},
  volume={40},
  number={3},
  pages={681--698},
  year={2001},
  publisher={SIAM},
  doi={10.1137/S0363012999361974}
}

@inproceedings{chen2023provably,
  title={Provably efficient iterated {CVaR} reinforcement learning with function approximation and human feedback},
  author={Chen, Yu and Du, Yihan and Hu, Pihe and Wang, Siwei and Wu, Desheng and Huang, Longbo},
  booktitle={The Twelfth International Conference on Learning Representations},
  year={2024},
  url={https://openreview.net/forum?id=vW1SkPl4kp}
}

@article{yang2024relative,
  title={Relative {Q}-learning for average-reward {Markov} decision processes with continuous states},
  author={Yang, Xiangyu and Hu, Jiaqiao and Hu, Jian-Qiang},
  journal={IEEE Transactions on Automatic Control},
  year={2024},
  pages={6546--6560},
  volume={69},
  number={10},
  publisher={IEEE},
  doi={10.1109/TAC.2024.3371380}
}

@article{bonetti2023risk,
  title={Risk-averse optimization of reward-based coherent risk measures},
  author={Bonetti, Massimiliano and Bisi, Lorenzo and Restelli, Marcello},
  journal={Artificial Intelligence},
  volume={316},
  eid={103845},
  year={2023},
  publisher={Elsevier}
}

@article{chow2018risk,
  title={Risk-constrained reinforcement learning with percentile risk criteria},
  author={Chow, Yinlam and Ghavamzadeh, Mohammad and Janson, Lucas and Pavone, Marco},
  journal={Journal of Machine Learning Research},
  volume={18},
  number={167},
  pages={1--51},
  year={2018}
}

@article{jiang2018risk,
  title={Risk-averse approximate dynamic programming with quantile-based risk measures},
  author={Jiang, Daniel R and Powell, Warren B},
  journal={Mathematics of Operations Research},
  volume={43},
  number={2},
  pages={554--579},
  year={2018},
  publisher={INFORMS}
}

@book{prashanth2022risk,
  title={Risk-Sensitive Reinforcement Learning via Policy Gradient Search},
  author={Prashanth, L. A. and Fu, Michael C},
  journal={Foundations and Trends{\textregistered} in Machine Learning},
  year={2022},
  publisher={Now Publishers, Inc.}
}

@book{borkar2008stochastic,
  title={Stochastic Approximation: A Dynamical Systems Viewpoint},
  author={Borkar, Vivek S},
  volume={9},
  year={2008},
  publisher={Springer}
}

@article{wang2025,
  title={{Statistical inference in conditional value-at-risk optimization}},
  author={Wang, Qixin and Cao, Hao and Hu, Jian-Qiang},
  publisher={Available at SSRN: https://ssrn.com/abstract=5168942},
  year={2025},
  journal = {SSRN},
url = {https://ssrn.com/abstract=5168942}
}

@article{pflug2016time,
  title={Time-consistent decisions and temporal decomposition of coherent risk functionals},
  author={Pflug, Georg Ch and Pichler, Alois},
  journal={Mathematics of Operations Research},
  volume={41},
  number={2},
  pages={682--699},
  year={2016},
  publisher={INFORMS}
}

@inproceedings{chow2015risk,
  title={Risk-sensitive and robust decision-making: A {CVaR} optimization approach},
  author={Chow, Yinlam and Tamar, Aviv and Mannor, Shie and Pavone, Marco},
  booktitle={Proceedings of the 28th International Conference on Neural Information Processing Systems-Volume 1},
  pages={1522--1530},
  year={2015}
}

@inproceedings{borkar2010learning,
  title={Learning algorithms for risk-sensitive control},
  author={Borkar, Vivek S},
  booktitle={Proceedings of the 19th International Symposium on Mathematical Theory of Networks and Systems--MTNS},
  volume={5},
  year={2010}
}

@inproceedings{stanko2021cvar,
  title={{CVaR} {Q}-learning},
  author={Stanko, Silvestr and Macek, Karel},
  booktitle={Computational Intelligence: 11th International Joint Conference, {IJCCI} 2019, {Vienna}, {Austria}, {September} 17--19, 2019, Revised Selected Papers},
  pages={333--358},
  year={2021},
  organization={Springer},
  publisher={Springer},
  editor={Merelo, Juan Juli{\'a}n and Garibaldi, Jonathan and Linares-Barranco, Alejandro and Warwick, Kevin and Madani, Kurosh},
  doi={10.1007/978-3-030-70594-7_14}
}

@article{hu2024simulation,
  title={Simulation optimization of conditional value-at-risk},
  author={Hu, Jiaqiao and Song, Meichen and Fu, Michael C and Peng, Yijie},
  journal={IISE Transactions},
  pages={1167--1181},
  volume={57},
  number={10},
  year={2025},
  publisher={Taylor \& Francis}
}

@incollection{uryasev2001conditional,
  title={Conditional value-at-risk: Optimization approach},
  author={Uryasev, Stanislav and Rockafellar, R Tyrrell},
  booktitle={Stochastic Optimization: Algorithms and Applications},
  pages={411--435},
  year={2001},
  publisher={Springer}
}

@article{brown2007large,
  title={Large deviations bounds for estimating conditional value-at-risk},
  author={Brown, David B},
  journal={Operations Research Letters},
  volume={35},
  number={6},
  pages={722--730},
  year={2007},
  publisher={Elsevier}
}

@article{rojas2025burning,
    title={Burning {RED}: Unlocking Subtask-Driven Reinforcement Learning and Risk-Awareness in Average-Reward {Markov} Decision Processes},
    author={Rojas, Juan Sebastian and Lee, Chi-Guhn},
    journal={Reinforcement Learning Journal},
    volume={6},
    pages={431--477},
    year={2025}
}

\FloatBarrier
\begin{APPENDICES}
\numberwithin{equation}{section}
\numberwithin{table}{section}
\numberwithin{figure}{section}
\numberwithin{lemma}{section}
\numberwithin{theorem}{section}
\numberwithin{proposition}{section}
\numberwithin{corollary}{section}
\numberwithin{assumption}{section}
\numberwithin{remark}{section}

\section{Auxiliary Results}\label{sec_auxiliary_app}
\begin{lemma}\label{lemma_continuity_app}
Under Assumptions~\ref{A.cost}, \ref{A.function}, and~\ref{A.mixing}, the stationary quantile is unique, and $d\mapsto\operatorname{VaR}^d$ is Lipschitz continuous with constant $(L_\mu+1)/\varepsilon_f$.
\end{lemma}
\begin{proof}[Proof of Lemma~\ref{lemma_continuity_app}]
For a stationary policy $d$, define $g^d(v,s):=\sum_{a\in\mathcal A(s)}d(s,a)F(v;s,a)$, so that $F^d(v)=\sum_s\mu^d(s)g^d(v,s)$. Let $d$ and $d'$ be two stationary policies. Adding and subtracting $\sum_s\mu^{d'}(s)g^d(v,s)$ gives
\begin{align*}
F^d(v)-F^{d'}(v)
&=\sum_{s\in\mathcal S}\left(\mu^d(s)-\mu^{d'}(s)\right)g^d(v,s)
 +\sum_{s\in\mathcal S}\mu^{d'}(s)\left(g^d(v,s)-g^{d'}(v,s)\right).
\end{align*}
Taking absolute values and applying the triangle inequality yields
\begin{align*}
\left|F^d(v)-F^{d'}(v)\right| &\le\sum_{s\in\mathcal S}\left|\mu^d(s)-\mu^{d'}(s)\right|\left|g^d(v,s)\right| +\sum_{s\in\mathcal S}\mu^{d'}(s)\left|g^d(v,s)-g^{d'}(v,s)\right|\\
&\le\left\|\mu^d-\mu^{d'}\right\|_1 +\max_{s\in\mathcal S}\left\|d(s)-d'(s)\right\|_1 \le\left(L_\mu+1\right)\left\|d-d'\right\|.
\end{align*}
The second inequality follows from $0\le g^d(v,s)\le1$ and
\begin{align*}
\left|g^d(v,s)-g^{d'}(v,s)\right| &=\left|\sum_{a\in\mathcal A(s)}\left(d(s,a)-d'(s,a)\right)F(v;s,a)\right| \le\sum_{a\in\mathcal A(s)}\left|d(s,a)-d'(s,a)\right|,
\end{align*}
where the last inequality uses $0\le F(v;s,a)\le1$. The third inequality follows from Lemma~\ref{lemma_stationary_lipschitz}. Every stationary VaR lies in $[v_{\min},v_{\max}]$. Indeed, $F^d(v_{\min})\le\phi$ and $F^d(v_{\max})\ge\phi$. By Assumption~\ref{A.function}, every stationary mixture density is no smaller than $\varepsilon_f$ on this interval, and the stationary VaR is unique.

Let $q=\operatorname{VaR}^d$ and $q'=\operatorname{VaR}^{d'}$. Suppose first that $q\ge q'$. By the mean value theorem, there is a point $\bar q$ between $q'$ and $q$ such that
\begin{align*}
F^d(q)-F^d(q') &=\left(q-q'\right)\sum_{s\in\mathcal S}\mu^d(s)\sum_{a\in\mathcal A(s)}d(s,a)f(\bar q;s,a) \ge\varepsilon_f\left(q-q'\right).
\end{align*}
The inequality follows because both quantiles and $\bar q$ lie in $[v_{\min},v_{\max}]$. Using $F^d(q)=F^{d'}(q')=\phi$ gives
\begin{align*}
\varepsilon_f\left(q-q'\right) &\le F^d(q)-F^d(q') =F^{d'}(q')-F^d(q') \le\left(L_\mu+1\right)\left\|d-d'\right\|.
\end{align*}
The first inequality is the density lower bound, the equality uses the two quantile equations, and the last inequality is the distribution bound established above. Interchanging $d$ and $d'$ covers the other ordering. Therefore,
\begin{align}\label{Eq.VaR_Lipschitz_app}
\left|\operatorname{VaR}^d-\operatorname{VaR}^{d'}\right|
&\le\frac{L_\mu+1}{\varepsilon_f}\left\|d-d'\right\|.
\end{align}

This proves the quantile bound without using Assumption~\ref{A.Qlearning}.
\end{proof}

\begin{lemma}\label{lemma_poisson_app}
Under Assumptions~\ref{A.cost}, \ref{A.function}, and~\ref{A.mixing}, for $v\in\mathbb R$ and $d\in\mathcal D$, define $h(v,d,s) =\sum_{a\in\mathcal A(s)}d(s,a)F(v;s,a)-F^d(v)$.
There is a solution $u(v,d,\cdot)$ of
\begin{align}\label{Eq.poisson_app}
u(v,d,s)-\sum_{s'\in\mathcal S}P^d(s,s')u(v,d,s')
&=h(v,d,s),
\end{align}
normalized by $\sum_s\mu^d(s)u(v,d,s)=0$. Moreover, there are constants $C_u,L_u<\infty$ such that $\max_s|u(v,d,s)| \le C_u$, $\max_s|u(v,d,s)-u(v',d',s)| \le L_u(|v-v'|+\|d-d'\|)$.
\end{lemma}
\begin{proof}[Proof of Lemma~\ref{lemma_poisson_app}]
Assumption~\ref{A.mixing} implies a common Doeblin condition over $m_p$ steps. Therefore, there are constants $K_p<\infty$ and $\rho_p\in(0,1)$ such that, for every stationary policy $d$,
\begin{align}\label{Eq.geometric_mixing_app}
\max_{s\in\mathcal S}\left\|\left(P^d\right)^k(s,\cdot)-\mu^d\right\|_1
&\le K_p\rho_p^k,
\qquad k\ge0.
\end{align}
To see this, apply the minorization after every block of $m_p$ transitions. At the end of each block, the common part $\delta_p\nu$ can be coupled, so the remaining total variation distance is multiplied by at most $1-\delta_p$. To include the boundary case $\delta_p=1$, use the smaller valid minorization constant $\delta_p/2$ and take $\rho_p=(1-\delta_p/2)^{1/m_p}\in(0,1)$. Enlarging the constant for the incomplete final block then gives~\eqref{Eq.geometric_mixing_app}.

For fixed $v$ and $d$, define $u(v,d,s) =\sum_{j=0}^{\infty}\sum_{s'\in\mathcal S}[(P^d)^j(s,s')-\mu^d(s')]h(v,d,s')$.
Since $|h(v,d,s)|\le1$, Equation~\eqref{Eq.geometric_mixing_app} gives
\begin{align*}
\left|\sum_{s'\in\mathcal S}\left[\left(P^d\right)^j(s,s')-\mu^d(s')\right]h(v,d,s')\right| &\le\left\|\left(P^d\right)^j(s,\cdot)-\mu^d\right\|_1 \le K_p\rho_p^j.
\end{align*}
The first inequality follows from $|h|\le1$, and the second is~\eqref{Eq.geometric_mixing_app}. Because the geometric series is summable, the series defining $u$ converges absolutely and uniformly, and $\max_s|u(v,d,s)| \le K_p/(1-\rho_p)$.
Applying $P^d$ term by term and shifting the summation index gives $\sum_{s'}P^d(s,s')u(v,d,s') =\sum_{j=1}^{\infty}\sum_{s'}[(P^d)^j(s,s')-\mu^d(s')]h(v,d,s')$.
Subtracting this expression from the series for $u$ leaves only the term with $j=0$, which is $h(v,d,s)$. Thus,~\eqref{Eq.poisson_app} holds. The normalization follows because $\mu^d(P^d)^j=\mu^d$ for every $j$, so every term in the stationary average of the series is zero.

It remains to prove the Lipschitz bound. Let \(A(d)=I-P^d+\mathbf{1}\mu^d\), where \(\mu^d\) is viewed as a row vector and \(\mathbf{1}\) is the column vector of ones. On the finite state space, the normalized Poisson equation can be written as $A(d)u(v,d)=h(v,d)$. The matrix $A(d)$ is nonsingular because the eigenvalue one of $P^d$ is simple under the uniform minorization. The policy space is compact and $A(d)$ depends continuously on $d$, so $\sup_{d\in\mathcal D}\|A(d)^{-1}\| <\infty$.
The transition matrix $P^d$ is linear in $d$, and Lemma~\ref{lemma_stationary_lipschitz} makes $\mu^d$ Lipschitz continuous. Hence, for a finite constant $C_A$, $\|A(d)-A(d')\| \le C_A\|d-d'\|$.
The identity $A(d)^{-1}-A(d')^{-1}=A(d)^{-1}[A(d')-A(d)]A(d')^{-1}$ then shows that $A(d)^{-1}$ is Lipschitz continuous in $d$.

We also have $\max_s|h(v,d,s)-h(v',d',s)| \le 2M_f|v-v'|+(L_\mu+2)\|d-d'\|$.
The term involving $v-v'$ follows from the density upper bound. The policy term follows by adding and subtracting the same conditional and stationary distribution terms, as in the proof of Lemma~\ref{lemma_continuity_app}. Finally,
\begin{align*}
u(v,d)-u(v',d')
&=A(d)^{-1}\left[h(v,d)-h(v',d')\right]
 +\left[A(d)^{-1}-A(d')^{-1}\right]h(v',d').
\end{align*}
Taking norms, using the uniform inverse bound, the Lipschitz bounds just obtained, and $\|h\|_\infty\le1$ proves the stated inequality for a finite constant $L_u$. This completes the proof.
\end{proof}

\begin{lemma}\label{lemma_scalar_tracking_app}
Let $x_{n+1}=(1-a_nc_n)x_n+a_n\xi_{n+1}+r_n$, where $a_n>0$, $\sum_na_n=\infty$, $a_n\to0$, and $0<c_0\le c_n\le c_1<\infty$. If $\sum_na_n\xi_{n+1}$ converges and $r_n/a_n\to0$, then $x_n\to0$.
\end{lemma}
\begin{proof}[Proof of Lemma~\ref{lemma_scalar_tracking_app}]
Fix an integer $N$ large enough that $0<a_nc_n<1$ for all $n\ge N$. For $N\le j\le n$, define $p_{j,n} =\prod_{k=j}^n(1-a_kc_k)$,
with $p_{n+1,n}=1$. Repeated substitution in the recursion gives
\begin{align}\label{Eq.scalar_variation_app}
x_{n+1}
&=p_{N,n}x_N+\sum_{j=N}^np_{j+1,n}a_j\xi_{j+1}+\sum_{j=N}^np_{j+1,n}r_j.
\end{align}
This identity follows by induction. It is true for $n=N$, and multiplying the representation at time $n$ by $1-a_{n+1}c_{n+1}$ and adding the new two terms gives the representation at time $n+1$.

Because $c_n\ge c_0$ and $1-y\le e^{-y}$ for $y\ge0$, $0\le p_{N,n} \le\exp\{-c_0\sum_{k=N}^na_k\}\longrightarrow0$.
Thus, the first term on the right hand side of~\eqref{Eq.scalar_variation_app} tends to zero.

Let $S_n=\sum_{k=0}^{n-1}a_k\xi_{k+1}$ and let $S_\infty$ be its limit. The second term in~\eqref{Eq.scalar_variation_app} can be written as $\sum_{j=N}^np_{j+1,n}(S_{j+1}-S_j)$.
Applying discrete summation by parts after subtracting $S_\infty$ gives the exact identity
\begin{align*}
\sum_{j=N}^np_{j+1,n}\left(S_{j+1}-S_j\right) &=\left(S_{n+1}-S_\infty\right)-p_{N+1,n}\left(S_N-S_\infty\right) -\sum_{j=N+1}^n\left(p_{j+1,n}-p_{j,n}\right)\left(S_j-S_\infty\right).
\end{align*}
The identity follows by writing the two sums separately, shifting the index in the terms containing $S_{j+1}$, and collecting the coefficients of each $S_j-S_\infty$. The sequence $p_{j+1,n}$ is nondecreasing in $j$, its final value is one, and $\sum_{j=N+1}^n(p_{j+1,n}-p_{j,n}) =1-p_{N+1,n}$.
Taking absolute values in the preceding identity therefore yields $|\sum_{j=N}^np_{j+1,n}(S_{j+1}-S_j)| \le2\sup_{j\ge N}|S_j-S_\infty|$.
Since $S_j\to S_\infty$, the right hand side tends to zero as $N$ tends to infinity, uniformly in $n\ge N$. Therefore, the noise term tends to zero.

Write $r_j=a_j\eta_j$, where $\eta_j\to0$. Since $c_j\ge c_0$, $\sum_{j=N}^np_{j+1,n}a_j \le (1/c_0) \sum_{j=N}^np_{j+1,n}a_jc_j = (1/c_0) (1-p_{N,n}) \le 1/c_0$.
The equality follows from $p_{j+1,n}-p_{j,n}=p_{j+1,n}a_jc_j$ and telescoping. Consequently, $|\sum_{j=N}^np_{j+1,n}r_j| \le (1/c_0) \sup_{j\ge N}|\eta_j|$.
The right hand side tends to zero as $N$ tends to infinity. All three terms in~\eqref{Eq.scalar_variation_app} therefore tend to zero, which proves $x_n\to0$.
\end{proof}

\begin{lemma}\label{lemma_Q_continuity_app}
Under Assumptions~\ref{A.cost} and~\ref{A.Qcontraction}, $Q(v)$ exists uniquely and is Lipschitz continuous in $v$ on $[v_{\min}-r_v,v_{\max}+r_v]$. If Assumptions~\ref{A.function} and~\ref{A.mixing} also hold, $d\mapsto Q^d$ is Lipschitz continuous in the weighted and ordinary maximum norms.
\end{lemma}
\begin{proof}[Proof of Lemma~\ref{lemma_Q_continuity_app}]
Finite first moments make $\mathcal T_vQ$ finite for every finite $v$ and every $Q$-function. The space of feasible Q entries is finite dimensional and complete in $\|\cdot\|_w$. Condition~\ref{A.Qcontraction} therefore gives a unique fixed point by the contraction mapping theorem. Although the condition is stated on the interval needed for tracking, its difference in $Q$ is independent of $v$, so fixed point existence also holds for each other fixed finite threshold. The positive part function is Lipschitz continuous with constant one. Hence, for every $v,v'\in[v_{\min}-r_v,v_{\max}+r_v]$ and every state and action pair,
\begin{align*}
\left|\tilde c(v,s,a)-\tilde c(v',s,a)\right| &\le\left|v-v'\right|+\frac{1}{1-\phi}\E\left[\left|\left(C(s,a)-v\right)^+-\left(C(s,a)-v'\right)^+\right|\right] \le\frac{2-\phi}{1-\phi}\left|v-v'\right|.
\end{align*}
The first inequality follows from the definition of $\tilde c$, and the second follows from $|(x-v)^+-(x-v')^+|\le|v-v'|$. Because $w(s,a)\ge1$, the same bound holds in the weighted maximum norm: $\|\mathcal T_vQ-\mathcal T_{v'}Q\|_w \le ((2-\phi)/(1-\phi)) |v-v'|$.
Because $Q(v)=\mathcal T_vQ(v)$ and $Q(v')=\mathcal T_{v'}Q(v')$, adding and subtracting $\mathcal T_vQ(v')$ gives
\begin{align*}
\left\|Q(v)-Q(v')\right\|_w
&\le\left\|\mathcal T_vQ(v)-\mathcal T_vQ(v')\right\|_w
 +\left\|\mathcal T_vQ(v')-\mathcal T_{v'}Q(v')\right\|_w\\
&\le\rho_Q\left\|Q(v)-Q(v')\right\|_w
 +\frac{2-\phi}{1-\phi}|v-v'|.
\end{align*}
The first inequality is the triangle inequality. The second uses Assumption~\ref{A.Qcontraction} and the preceding cost bound. Moving the contraction term to the left gives
\begin{align}\label{Eq.Q_Lipschitz_app}
\left\|Q(v)-Q(v')\right\|_w
&\le\frac{2-\phi}{(1-\rho_Q)(1-\phi)}|v-v'|.
\end{align}
Since $\|Q\|_\infty\le w_{\max}\|Q\|_w$, the same mapping is Lipschitz continuous in the ordinary maximum norm. Finally, apply Lemma~\ref{lemma_continuity_app} to $Q^d=Q(\operatorname{VaR}^d)$ and combine~\eqref{Eq.VaR_Lipschitz_app} with~\eqref{Eq.Q_Lipschitz_app}. This proves the policy Lipschitz property with the additional model conditions stated in the lemma. This completes the proof.
\end{proof}

For a state and action pair $(s,a)$, let $I_n(s,a):=\I\{s_n=s,a_n=a\}$. On the event $I_n(s,a)=1$, define
\begin{align*}
\zeta_{n+1}(s,a) &:=\widetilde C(v_n,s,a)-\tilde c(v_n,s,a) +\min_{a'\in\mathcal A(s_{n+1})}Q_n(s_{n+1},a') -\sum_{s'\in\mathcal S}p(s'\mid s,a)\min_{a'\in\mathcal A(s')}Q_n(s',a'),
\end{align*}
and set $\zeta_{n+1}(s,a)=0$ when $I_n(s,a)=0$. Define
\begin{align}\label{Eq.Q_martingale_definition_app}
M_m(s,a)
&:=\sum_{n=0}^{m-1}\beta_n(s,a)I_n(s,a)\zeta_{n+1}(s,a),
\qquad m\ge1.
\end{align}

\begin{lemma}\label{lemma_projected_cycle_app}
Assume Assumption~\ref{A.Qcontraction} and the local gain bounds $0<\beta_k(s,a)<1$ from Assumption~\ref{A.beta}. Fix integers $p<q$ and a constant $T>0$. Suppose that $v_p\in[v_{\min}-r_v,v_{\max}+r_v]$, $|Q(v_p)(s,a)|\le Q_{\max}-1$ for every pair, and $\sum_{k=p}^{q-1}\beta_k(s,a)I_k(s,a) \ge T,~ s\in\mathcal S,~a\in\mathcal A(s)$.
Let
\begin{align*}
\Delta_M(p,q) &:=\max_{s\in\mathcal S,~a\in\mathcal A(s)}\max_{p\le m\le q}|M_m(s,a)-M_p(s,a)|, \qquad \Delta_v(p,q) :=\max_{p\le m\le q}|v_m-v_p|.
\end{align*}
If $\Delta_M(p,q)<1$, then, with $L_c=(2-\phi)/(1-\phi)$ and $\kappa_T=\rho_Q+(1-\rho_Q)e^{-T}$,
\begin{align}
\left\|Q_q-Q(v_p)\right\|_w &\le\kappa_T\left\|Q_p-Q(v_p)\right\|_w+\Delta_M(p,q) +\frac{(1+\rho_Q)\Delta_M(p,q)+L_c\Delta_v(p,q)}{1-\rho_Q}.\label{Eq.projected_cycle_endpoint_app}
\end{align}
and, for every $p\le m\le q$,
\begin{align}
\left\|Q_m-Q(v_p)\right\|_w &\le\left\|Q_p-Q(v_p)\right\|_w+\Delta_M(p,q) +\frac{(1+\rho_Q)\Delta_M(p,q)+L_c\Delta_v(p,q)}{1-\rho_Q}.\label{Eq.projected_cycle_uniform_app}
\end{align}
\end{lemma}

\begin{proof}[Proof of Lemma~\ref{lemma_projected_cycle_app}]
Write $\bar Q=Q(v_p)$ and, for $p\le k\le q$, define $m_k(s,a):=M_k(s,a)-M_p(s,a)$ and $Y_k(s,a):=Q_k(s,a)-m_k(s,a)$. Thus, $m_p=0$ and $Y_p=Q_p$. We first compare $Y_k$ with the fixed point $\bar Q$.

Consider an iteration $k$ at which component $i=(s,a)$ is visited, and write $b_k=\beta_k(s,a)$. By~\eqref{Eq.Q_martingale_definition_app}, $m_{k+1}(i)-m_k(i) =b_k\zeta_{k+1}(i)$.
Since $|\bar Q(i)|\le Q_{\max}-1$ and $|m_{k+1}(i)|\le\Delta_M(p,q)<1$, the point $\bar Q(i)+m_{k+1}(i)$ belongs to $[-Q_{\max},Q_{\max}]$. Projection onto this interval is nonexpansive. Therefore,
\begin{align*}
&|Y_{k+1}(i)-\bar Q(i)| =\left|\Pi_{Q_{\max}}\left((1-b_k)Q_k(i)+b_k(\mathcal T_{v_k}Q_k)(i)+b_k\zeta_{k+1}(i)\right)-m_{k+1}(i)-\bar Q(i)\right|\\
&\quad\le\left|(1-b_k)Q_k(i)+b_k(\mathcal T_{v_k}Q_k)(i)+b_k\zeta_{k+1}(i)-m_{k+1}(i)-\bar Q(i)\right|\\
&\quad=\left|(1-b_k)\left[Y_k(i)-\bar Q(i)\right]+b_k\left[(\mathcal T_{v_k}Q_k)(i)-(\mathcal T_{v_p}\bar Q)(i)\right]-b_km_k(i)\right|.
\end{align*}
The inequality follows by comparing the projections of the first argument and $\bar Q(i)+m_{k+1}(i)$, whose projection is itself. The equality follows from $Q_k=Y_k+m_k$, the martingale increment above, and $\mathcal T_{v_p}\bar Q=\bar Q$.

Add and subtract $(\mathcal T_{v_p}Q_k)(i)$ in the last bracket. Assumption~\ref{A.Qcontraction}, the weighted norm definition, and the Lipschitz property of the transformed one stage cost give
\begin{align*}
\frac{|(\mathcal T_{v_k}Q_k)(i)-(\mathcal T_{v_p}\bar Q)(i)|}{w(i)} &\le\frac{|(\mathcal T_{v_p}Q_k)(i)-(\mathcal T_{v_p}\bar Q)(i)|}{w(i)} +\frac{|(\mathcal T_{v_k}Q_k)(i)-(\mathcal T_{v_p}Q_k)(i)|}{w(i)}\\
&\le\rho_Q\|Q_k-\bar Q\|_w+L_c|v_k-v_p| \le\rho_Q\|Y_k-\bar Q\|_w+\rho_Q\|m_k\|_w+L_c\Delta_v(p,q)\\
&\le\rho_Q\|Y_k-\bar Q\|_w+\rho_Q\Delta_M(p,q)+L_c\Delta_v(p,q).
\end{align*}
The first inequality is the triangle inequality. The second uses Assumption~\ref{A.Qcontraction} and the bound $\|\mathcal T_vQ-\mathcal T_{v'}Q\|_w\le L_c|v-v'|$. The third uses $Q_k=Y_k+m_k$. The last inequality follows from $w(i)\ge1$, which gives $\|m_k\|_w\le\Delta_M(p,q)$.

Set $D_k :=\|Y_k-\bar Q\|_w, ~ \eta_{p,q}:=(1+\rho_Q)\Delta_M(p,q)+L_c\Delta_v(p,q)$.
Dividing the preceding projection inequality by $w(i)$ and using $|m_k(i)|/w(i)\le\Delta_M(p,q)$ yields
\begin{align}\label{Eq.projected_component_recursion_app}
\frac{|Y_{k+1}(i)-\bar Q(i)|}{w(i)}
&\le(1-b_k)\frac{|Y_k(i)-\bar Q(i)|}{w(i)}+b_k\rho_QD_k+b_k\eta_{p,q}.
\end{align}
Every unvisited component of $Y_k$ remains unchanged. Let $B_{p,q} :=\max\{D_p, (\eta_{p,q}/(1-\rho_Q)) \}$.
If $D_k\le B_{p,q}$, then the right hand side of~\eqref{Eq.projected_component_recursion_app} is no larger than
\begin{align*}
(1-b_k)B_{p,q}+b_k\rho_QB_{p,q}+b_k\eta_{p,q} &=B_{p,q}-b_k\left[(1-\rho_Q)B_{p,q}-\eta_{p,q}\right] \le B_{p,q}.
\end{align*}
The last inequality follows from the definition of $B_{p,q}$. Since $D_p\le B_{p,q}$, induction over the calendar iterations gives $D_k\le B_{p,q}$ for every $p\le k\le q$. Hence, $\|Q_k-\bar Q\|_w \le\|Y_k-\bar Q\|_w+\|m_k\|_w \le D_p+ \eta_{p,q}/(1-\rho_Q) +\Delta_M(p,q)$,
which is~\eqref{Eq.projected_cycle_uniform_app}.

It remains to obtain the strict reduction at the end of the cycle. Fix a component $i$, and list its visit times in $[p,q)$ in increasing order. Between two visits its value is unchanged. Repeated substitution in~\eqref{Eq.projected_component_recursion_app}, using $D_k\le B_{p,q}$ at every visit, gives $(|Y_q(i)-\bar Q(i)|)/(w(i)) \le\vartheta_{p,q}(i)D_p+[1-\vartheta_{p,q}(i)](\rho_QB_{p,q}+\eta_{p,q})$,
where $\vartheta_{p,q}(i) :=\prod_{k=p}^{q-1}(1-\beta_k(i)I_k(i))$.
The formula follows by induction over the visits: after each visit, the previous component error keeps weight $1-b_k$, and the common upper bound $\rho_QB_{p,q}+\eta_{p,q}$ receives weight $b_k$. Since the accumulated local gain is at least $T$ and $1-x\le e^{-x}$, $\vartheta_{p,q}(i) \le\exp\{-\sum_{k=p}^{q-1}\beta_k(i)I_k(i)\} \le e^{-T}$.
Also, $B_{p,q}\le D_p+\eta_{p,q}/(1-\rho_Q)$, and hence $\rho_QB_{p,q}+\eta_{p,q} \le\rho_QD_p+ \eta_{p,q}/(1-\rho_Q)$.
Substituting this bound into the repeated substitution formula gives
\begin{align*}
\frac{|Y_q(i)-\bar Q(i)|}{w(i)} &\le\vartheta_{p,q}(i)D_p+\left[1-\vartheta_{p,q}(i)\right] \left[\rho_QD_p+\frac{\eta_{p,q}}{1-\rho_Q}\right]\\
&=\left[\rho_Q+(1-\rho_Q)\vartheta_{p,q}(i)\right]D_p +\frac{1-\vartheta_{p,q}(i)}{1-\rho_Q}\eta_{p,q} \le\left[\rho_Q+(1-\rho_Q)e^{-T}\right]D_p +\frac{\eta_{p,q}}{1-\rho_Q}.
\end{align*}
The first inequality uses the upper bound on $\rho_QB_{p,q}+\eta_{p,q}$. The equality collects the coefficients of $D_p$. The last inequality uses $0\le\vartheta_{p,q}(i)\le e^{-T}$ for the first coefficient and $0\le1-\vartheta_{p,q}(i)\le1$ for the perturbation coefficient. Taking the maximum over the components and then adding $\|m_q\|_w\le\Delta_M(p,q)$ proves~\eqref{Eq.projected_cycle_endpoint_app}. This completes the proof.
\end{proof}

\section{Proofs for Section~\ref{sec_prob_alg}}\label{sec_preliminary_app}
\begin{proof}[Proof of Lemma~\ref{lemma_exchange_VaR}]
Fix a stationary policy $d$ and an initial state $s_0=s$. Let $F_n^d(x) =\PP^d(C(s_n,a_n)\le x\mid s_0=s)$.
Conditioning on the state and action at time $n$ gives $F_n^d(x) =\sum_{s'\in\mathcal S}\sum_{a\in\mathcal A(s')}\PP^d(s_n=s',a_n=a\mid s_0=s)F(x;s',a)$.
Under Assumption~\ref{A.mixing}, the state and action probabilities converge to $\pi^d(s',a)$. Since the number of state and action pairs is finite and $0\le F(x;s',a)\le1$,
\begin{align*}
\sup_{x\in\mathbb R}\left|F_n^d(x)-\PP(C^d\le x)\right|
&\le\sum_{s'\in\mathcal S}\sum_{a\in\mathcal A(s')}\left|\PP^d\left(s_n=s',a_n=a\mid s_0=s\right)-\pi^d(s',a)\right|\longrightarrow0.
\end{align*}
The inequality follows by subtracting the two mixture representations and applying the triangle inequality. Thus, the time $n$ cost distribution converges uniformly to the stationary cost distribution.

Let $q=\operatorname{VaR}[C^d]$ and let $q_n$ be the VaR of $C(s_n,a_n)$ conditional on $s_0=s$. At $v_{\min}$, every component cumulative distribution function is no larger than $\phi$, whereas at $v_{\max}$ every component cumulative distribution function is no smaller than $\phi$. The same is therefore true for every finite mixture $F_n^d$ and for the stationary mixture. Consequently, $q_n,q\in[v_{\min},v_{\max}]$. By Assumption~\ref{A.function}, every finite mixture density and the stationary mixture density are therefore no smaller than $\varepsilon_f$ on $[v_{\min},v_{\max}]$, so all these quantiles are unique. In particular, the stationary density is positive at $q$. Therefore, for every sufficiently small $\eta>0$, $\PP(C^d\le q-\eta) <\phi <\PP(C^d\le q+\eta)$.
Uniform convergence of the distribution functions implies that the same two strict inequalities hold with $F_n^d$ for all sufficiently large $n$. By the definition of a quantile, $q-\eta<q_n<q+\eta$. Since $\eta$ is arbitrary, $q_n\to q$. The Cesaro average of a convergent sequence has the same limit, so $\operatorname{VaR}^d(s) =\lim_{N\to\infty} (1/N) \sum_{n=0}^{N-1}q_n = q =\operatorname{VaR}[C^d]$.
For continuous distributions, the CVaR representation is $\operatorname{CVaR}_n^d =q_n+ (1/(1-\phi)) \E^d[(C(s_n,a_n)-q_n)^+\mid s_0=s]$.
Let $p_n(s',a)=\PP^d(s_n=s',a_n=a\mid s_0=s)$. By conditioning on the state and action,
\begin{align*}
&\E^d\left[\left(C(s_n,a_n)-q_n\right)^+\mid s_0=s\right] =\sum_{s'\in\mathcal S}\sum_{a\in\mathcal A(s')}p_n(s',a)\E\left[\left(C(s',a)-q_n\right)^+\right].
\end{align*}
Add and subtract the same expression with $p_n$ replaced by $\pi^d$ and $q_n$ replaced by $q$. The absolute difference is bounded by
\begin{align*}
&\sum_{s'\in\mathcal S}\sum_{a\in\mathcal A(s')}\left|p_n(s',a)-\pi^d(s',a)\right|\E\left[\left|C(s',a)\right|+|q_n|+|q|\right]\\
&\quad+\sum_{s'\in\mathcal S}\sum_{a\in\mathcal A(s')}\pi^d(s',a)\E\left[\left|\left(C(s',a)-q_n\right)^+-\left(C(s',a)-q\right)^+\right|\right].
\end{align*}
The first term tends to zero because $p_n(s',a)\to\pi^d(s',a)$, the state and action spaces are finite, the first moments are finite, and the convergent sequence $\{q_n\}$ is bounded. The second term is no larger than $|q_n-q|$ because the positive part function is Lipschitz continuous with constant one. Hence, $\operatorname{CVaR}_n^d \longrightarrow q+ (1/(1-\phi)) \E[(C^d-q)^+] =\operatorname{CVaR}[C^d]$.
Taking Cesaro averages gives $\operatorname{CVaR}^d(s)=\operatorname{CVaR}[C^d]$. This completes the proof.
\end{proof}

\section{Proofs for Section~\ref{sec_converge}}\label{sec_convergence_app}

\begin{proof}[Proof of Lemma~\ref{lemma_stationary_lipschitz}]
Fix two stationary policies $d$ and $d'$. Under Assumption~\ref{A.mixing}, the stationary transition matrix of policy $d$ satisfies $(P^d)^{m_p}(s,y) \ge\delta_p\nu(y),~ s,y\in\mathcal S$.
If $\delta_p<1$, write $(P^d)^{m_p}(s,y) =\delta_p\nu(y)+(1-\delta_p)R^d(s,y)$,
where every row of $R^d$ is a probability vector. If $\delta_p=1$, the same argument below applies with the second term equal to zero. Let $z$ be a row vector satisfying $\sum_sz(s)=0$. The common part cancels, and hence $z(P^d)^{m_p} =(1-\delta_p)zR^d$.
Since multiplication by a stochastic matrix does not increase the one norm of a row vector, $\|zR^d\|_1 =\sum_{y\in\mathcal S}|\sum_{s\in\mathcal S}z(s)R^d(s,y)| \le\sum_{s\in\mathcal S}|z(s)|\sum_{y\in\mathcal S}R^d(s,y) =\|z\|_1$.
The first inequality is the triangle inequality, and the last equality uses that every row of $R^d$ sums to one. Consequently, $\|z(P^d)^{m_p}\|_1 \le(1-\delta_p)\|z\|_1$.
The equality preceding this bound uses $\sum_sz(s)=0$, which removes the common term $\delta_p\nu$.

We next compare the transition matrices generated by $d$ and $d'$. For every state $s$,
\begin{align*}
\sum_{y\in\mathcal S}\left|P^d\left(s,y\right)-P^{d'}\left(s,y\right)\right| &=\sum_{y\in\mathcal S}\left|\sum_{a\in\mathcal A(s)}\left[d(s,a)-d'(s,a)\right]p\left(y\mid s,a\right)\right|\\
&\le\sum_{a\in\mathcal A(s)}\left|d(s,a)-d'(s,a)\right| \sum_{y\in\mathcal S}p\left(y\mid s,a\right) =\left\|d(s)-d'(s)\right\|_1 \le\left\|d-d'\right\|.
\end{align*}
The first inequality is the triangle inequality, the second equality uses that each transition row sums to one, and the last inequality follows from the definition of the policy norm.

The difference of the $m_p$ step transition matrices can be expanded as $(P^d)^{m_p}-(P^{d'})^{m_p} =\sum_{\ell=0}^{m_p-1}(P^d)^\ell (P^d-P^{d'}) (P^{d'})^{m_p-1-\ell}$.
This identity follows by adding and subtracting the intermediate products in which the factors are changed from $P^d$ to $P^{d'}$ one at a time. Left multiplication by a stochastic matrix forms convex combinations of rows, and right multiplication by a stochastic matrix does not increase the one norm of a row. Therefore, $\max_{s\in\mathcal S}\sum_{y\in\mathcal S} |(P^d)^{m_p}(s,y) -(P^{d'})^{m_p}(s,y)| \le m_p\|d-d'\|$.
A stationary distribution exists because the state space is finite. Let $\mu^d$ and $\mu^{d'}$ be stationary distributions of the two transition matrices. Using stationarity and adding and subtracting $\mu^{d'}(P^d)^{m_p}$ gives $\mu^d-\mu^{d'} =(\mu^d-\mu^{d'})(P^d)^{m_p} +\mu^{d'}[(P^d)^{m_p} -(P^{d'})^{m_p}]$.
The vector $\mu^d-\mu^{d'}$ has components summing to zero. Applying the contraction to the first term and the row bound to the second term yields $\|\mu^d-\mu^{d'}\|_1 \le(1-\delta_p)\|\mu^d-\mu^{d'}\|_1 +m_p\|d-d'\|$.
Moving the first term on the right to the left and dividing by $\delta_p$ gives $\|\mu^d-\mu^{d'}\|_1 \le (m_p/\delta_p) \|d-d'\|$.
Finally, applying the same contraction to two stationary distributions of the same policy gives their equality. Hence, the stationary distribution is unique, and the stated Lipschitz bound follows. This completes the proof.
\end{proof}

The VaR recursion is not projected. The proof first uses the deterministic tail points $\underline v$ and $\overline v$ to show that the mean drift points toward a fixed interval. It then removes the cost observation martingale and applies a last crossing argument. This gives an integrable bound on the entire sample path rather than only a bound at each fixed iteration.

\begin{proof}[Proof of Lemma~\ref{lemma_v_bound}]
Conditional on $\mathscr G_n$, define $\zeta_{n+1} :=F(v_n;s_n,a_n)-\I\{C(s_n,a_n)\le v_n\}$.
The conditional cost model gives $\E[\zeta_{n+1}\mid\mathscr G_n]=0$ and $|\zeta_{n+1}|\le1$. Procedure~\eqref{Eq.v} can therefore be written as
\begin{align}\label{Eq.v_drift_noise_app}
v_{n+1}
&=v_n+\alpha_n\left[\phi-F(v_n;s_n,a_n)\right]+\alpha_n\zeta_{n+1}.
\end{align}
Let $S_n :=\sum_{k=0}^{n-1}\alpha_k\zeta_{k+1},~ S_0:=0$.
The sequence $\{S_n\}$ is a square integrable martingale. Orthogonality of the martingale differences gives $\sup_m\E[S_m^2] =\sup_m\sum_{k=0}^{m-1}\alpha_k^2\E[\zeta_{k+1}^2] \le\sum_{k=0}^{\infty}\alpha_k^2<\infty$.
The equality follows because the cross products of distinct martingale differences have zero expectation, and the inequality follows from $|\zeta_{k+1}|\le1$. Doob's $L^2$ maximal inequality now yields
\begin{align}
\E\left[\sup_{n\ge0}|S_n|^2\right] &\le4\sup_m\E[S_m^2] \le4\sum_{k=0}^{\infty}\alpha_k^2<\infty.\label{Eq.v_martingale_max_app}
\end{align}
In particular, $S_n$ converges and $\sup_n|S_n|<\infty$ almost surely.

We next use the direction of the mean drift. By the definition of $\underline v$, $\overline v$, and $\eta_v$,
\begin{align*}
v_n\ge\overline v &\quad\Longrightarrow\quad \phi-F(v_n;s_n,a_n)\le-\eta_v, \qquad v_n\le\underline v \quad\Longrightarrow\quad \phi-F(v_n;s_n,a_n)\ge\eta_v.
\end{align*}
The first implication follows because every cumulative distribution function is nondecreasing and $F(\overline v;s,a)\ge\phi+\eta_v$; the second follows in the same way from $F(\underline v;s,a)\le\phi-\eta_v$.

Fix a sample path and an integer $n$. We first obtain an upper bound. If $v_k>\overline v$ for all $0\le k\le n$, then every mean drift in~\eqref{Eq.v_drift_noise_app} is nonpositive, and hence $v_n =v_0+\sum_{k=0}^{n-1}\alpha_k[\phi-F(v_k;s_k,a_k)]+S_n \le v_0+|S_n|$.
If the path is at or below $\overline v$ at some time before $n$, let $\tau<n$ be the last index such that $v_\tau\le\overline v$. Since the absolute value of both terms inside the VaR increment is no larger than one, $v_{\tau+1} \le\overline v+2\alpha_\tau \le\overline v+2\alpha_0$,
where the last inequality uses the monotonicity of $\alpha_n$. By the definition of the last crossing, $v_k>\overline v$ for $\tau+1\le k<n$, so the corresponding mean drifts are nonpositive. Summing~\eqref{Eq.v_drift_noise_app} from $\tau+1$ to $n-1$ gives
\begin{align*}
v_n &=v_{\tau+1}+\sum_{k=\tau+1}^{n-1}\alpha_k\left[\phi-F(v_k;s_k,a_k)\right]+S_n-S_{\tau+1} \le\overline v+2\alpha_0+|S_n-S_{\tau+1}| \le\overline v+2\alpha_0+2\sup_{m\ge0}|S_m|.
\end{align*}
The second inequality uses the nonpositive drift after the last crossing, and the third is the triangle inequality. Combining the two cases gives
\begin{align}\label{Eq.v_upper_path_app}
\sup_{n\ge0}v_n
&\le\max\left\{v_0,\overline v+2\alpha_0\right\}+2\sup_{m\ge0}|S_m|.
\end{align}

The lower bound is obtained in the same way. If the path stays below $\underline v$, every mean drift is nonnegative. Otherwise, let $\sigma<n$ be the last index such that $v_\sigma\ge\underline v$. The one step change is no smaller than $-2\alpha_\sigma$, and every later mean drift before time $n$ is nonnegative. Therefore,
\begin{align}\label{Eq.v_lower_path_app}
\inf_{n\ge0}v_n
&\ge\min\left\{v_0,\underline v-2\alpha_0\right\}-2\sup_{m\ge0}|S_m|.
\end{align}
The derivation is identical to the upper bound with all inequalities reversed.

Let $B_v$ be the maximum of the absolute values of the two deterministic constants in~\eqref{Eq.v_upper_path_app} and~\eqref{Eq.v_lower_path_app}. These two displays imply $\sup_{n\ge0}|v_n| \le B_v+2\sup_{m\ge0}|S_m|$.
Using $(x+y)^2\le2x^2+2y^2$ and~\eqref{Eq.v_martingale_max_app}, we obtain $\E[\sup_{n\ge0}v_n^2] \le2B_v^2+8\E[\sup_{m\ge0}|S_m|^2] <\infty$.
This proves the stated result. The two consequences follow immediately because $v_n^2\le\sup_kv_k^2$ for every $n$, and an integrable random variable is finite almost surely. This completes the proof.
\end{proof}

\begin{proof}[Proof of Theorem~\ref{th_v}]
Because $\epsilon_n$ is decreasing, $d_n(s)\in\mathcal D_{\epsilon_{n-1}}\subseteq\mathcal D_{\epsilon_n}$. The projection $\Pi_n$ is nonexpansive, and $d_n(s)$ is feasible for the same projection. Hence, $\|d_{n+1}(s)-d_n(s)\|_2 \le\gamma_n\|\delta(\arg\min_aQ_{n+1}(s,a))-d_n(s)\|_2 \le2\gamma_n$.
The first inequality follows by comparing the projections of the unprojected policy update and $d_n(s)$. The second follows because both vectors inside the norm are probability vectors. Since all norms are equivalent in the finite dimensional policy space, there is a constant $C_d$ such that
\begin{align}\label{Eq.policy_increment_app}
\left\|d_{n+1}-d_n\right\|
&\le C_d\gamma_n.
\end{align}
Lemma~\ref{lemma_continuity_app} therefore gives
\begin{align}\label{Eq.VaR_target_increment_app}
\left|\operatorname{VaR}^{d_{n+1}}-\operatorname{VaR}^{d_n}\right|
&\le C_q\gamma_n
\end{align}
for a finite constant $C_q$.

We next identify the noise in Procedure~\eqref{Eq.v}. Let $g_n(s) =\sum_{a\in\mathcal A(s)}d_n(s,a)F(v_n;s,a)$.
Conditional on $\mathscr F_n$, the action is sampled from $d_n(s_n)$ and the cost has cumulative distribution function $F(\cdot;s_n,a_n)$. Consequently, $\E[\I\{C(s_n,a_n)\le v_n\}\mid\mathscr F_n] =g_n(s_n)$.
Add and subtract $g_n(s_n)$ and $F^{d_n}(v_n)$ in the VaR increment to obtain
\begin{align}\label{Eq.VaR_noise_decomposition_app}
\phi-\I\left\{C(s_n,a_n)\le v_n\right\}
&=\phi-F^{d_n}(v_n)+\xi_{n+1}-h(v_n,d_n,s_n),
\end{align}
where $\xi_{n+1} =g_n(s_n)-\I\{C(s_n,a_n)\le v_n\}$.
The sign of the final term follows from $h(v_n,d_n,s_n)=g_n(s_n)-F^{d_n}(v_n)$. The random variable $\xi_{n+1}$ is a martingale difference with respect to $\mathscr F_n$, and $|\xi_{n+1}|\le1$. Therefore, $\sum_{n=0}^{\infty}\alpha_n\xi_{n+1}$
converges almost surely, because its quadratic variation is bounded by $\sum_n\alpha_n^2<\infty$.

It remains to treat the state term. Let $u_n(s)=u(v_n,d_n,s)$ be the Poisson solution in Lemma~\ref{lemma_poisson_app}. Equation~\eqref{Eq.poisson_app} gives $h(v_n,d_n,s_n) =u_n(s_n)-\sum_{s'\in\mathcal S}P^{d_n}(s_n,s')u_n(s')$.
Adding and subtracting $u_n(s_{n+1})$ yields
\begin{align}\label{Eq.Poisson_split_app}
\alpha_nh(v_n,d_n,s_n)
&=\alpha_n\left[u_n(s_n)-u_n(s_{n+1})\right]
 +\alpha_n\left[u_n(s_{n+1})-\sum_{s'}P^{d_n}(s_n,s')u_n(s')\right].
\end{align}
Conditional on $\mathscr F_n$, the next state has transition matrix $P^{d_n}$ because the action is sampled from $d_n(s_n)$. Thus, the second bracket in~\eqref{Eq.Poisson_split_app} is a martingale difference. It is uniformly bounded by $2C_u$, so its weighted series converges almost surely by $\sum_n\alpha_n^2<\infty$.

We now show that the series formed by the first bracket also converges. For integers $N\le m$, discrete summation by parts gives
\begin{align*}
&\sum_{n=N}^m\alpha_n\left[u_n(s_n)-u_n(s_{n+1})\right] =\alpha_Nu_N(s_N)-\alpha_mu_m(s_{m+1}) +\sum_{n=N+1}^m\left[\alpha_nu_n(s_n)-\alpha_{n-1}u_{n-1}(s_n)\right].
\end{align*}
To verify the identity, write the two sums separately, shift the index of the terms containing $s_{n+1}$, and cancel the common interior terms. For each term in the final sum, add and subtract $\alpha_{n-1}u_n(s_n)$ to obtain
\begin{align*}
\left|\alpha_nu_n(s_n)-\alpha_{n-1}u_{n-1}(s_n)\right|
&\le C_u\left|\alpha_n-\alpha_{n-1}\right|
 +\alpha_{n-1}\max_s\left|u_n(s)-u_{n-1}(s)\right|\\
&\le C_u\left|\alpha_n-\alpha_{n-1}\right|
 +L_u\alpha_{n-1}\left(\left|v_n-v_{n-1}\right|+\left\|d_n-d_{n-1}\right\|\right)\\
&\le C_u\left|\alpha_n-\alpha_{n-1}\right|+L_u\alpha_{n-1}\left(\alpha_{n-1}+C_d\gamma_{n-1}\right).
\end{align*}
The first inequality is the triangle inequality. The second uses the Lipschitz bound in Lemma~\ref{lemma_poisson_app}. The third uses $|v_n-v_{n-1}|\le\alpha_{n-1}$ from Procedure~\eqref{Eq.v} and~\eqref{Eq.policy_increment_app} for the policy increment. Since $\alpha_n$ is nonincreasing, $\sum_n|\alpha_n-\alpha_{n-1}|<\infty$. Also, $\sum_n\alpha_n^2<\infty$, and $\gamma_n=o(\alpha_n)$ implies $\sum_n\alpha_n\gamma_n<\infty$. Hence, the last upper bound is summable. The two boundary terms converge because $\alpha_n\to0$ and $u_n$ is bounded. Therefore, $\sum_{n=0}^{\infty}\alpha_nh(v_n,d_n,s_n)$
converges almost surely.

Combining the two convergent series with~\eqref{Eq.VaR_noise_decomposition_app}, we can write Procedure~\eqref{Eq.v} as
\begin{align}\label{Eq.VaR_tracking_recursion_app}
v_{n+1}
&=v_n+\alpha_n\left[\phi-F^{d_n}(v_n)\right]+\alpha_n\eta_{n+1},
\end{align}
where the partial sums $\sum_n\alpha_n\eta_{n+1}$ converge almost surely.

Let $q_n=\operatorname{VaR}^{d_n}$ and $e_n=v_n-q_n$. Since $F^{d_n}(q_n)=\phi$, the mean value theorem gives a point $\widetilde v_n$ between $v_n$ and $q_n$ such that $\phi-F^{d_n}(v_n) =-c_ne_n$,
where $c_n =\sum_{s\in\mathcal S}\mu^{d_n}(s)\sum_{a\in\mathcal A(s)}d_n(s,a)f(\widetilde v_n;s,a)$.
We next show that the deterministic density interval in Assumption~\ref{A.function} contains the VaR iterates after a finite random time. Let $R_n =\sum_{k=0}^{n-1}\alpha_k\eta_{k+1},~R_0=0$.
The sequence $\{R_n\}$ converges almost surely. Fix a sample path on which it converges. Since $\alpha_n\to0$ and a convergent sequence is Cauchy, there is a finite random index $N_1$ such that, for all $m\ge n\ge N_1$, $\alpha_n \le r_v/8 ,~ |R_{m+1}-R_n| \le r_v/8$.
Suppose first that $v_n>\overline v$ for every $n\ge N_1$. The upper tail drift then satisfies $\phi-F^{d_n}(v_n)\le-\eta_v$. Summing~\eqref{Eq.VaR_tracking_recursion_app} gives
\begin{align*}
v_{m+1} &=v_{N_1}+\sum_{n=N_1}^{m}\alpha_n\left[\phi-F^{d_n}(v_n)\right] +R_{m+1}-R_{N_1} \le v_{N_1}-\eta_v\sum_{n=N_1}^{m}\alpha_n+\frac{r_v}{8}.
\end{align*}
The right hand side tends to minus infinity because $\sum_n\alpha_n=\infty$, which contradicts $v_{m+1}>\overline v$. Hence, there is an index $\tau\ge N_1$ such that $v_\tau\le\overline v$.

Fix $m\ge\tau$. If $v_m\le\overline v$, the desired upper bound is immediate. Otherwise, let $\ell<m$ be the last index in $[\tau,m]$ satisfying $v_\ell\le\overline v$. At the crossing step, the absolute value of the mean drift is at most one, and therefore
\begin{align*}
v_{\ell+1} &=v_\ell+\alpha_\ell\left[\phi-F^{d_\ell}(v_\ell)\right] +R_{\ell+1}-R_\ell \le\overline v+\alpha_\ell+\left|R_{\ell+1}-R_\ell\right| \le\overline v+\frac{r_v}{4}.
\end{align*}
By the definition of $\ell$, every later iterate up to time $m$ is above $\overline v$, so every corresponding mean drift is nonpositive. Summing from $\ell+1$ to $m-1$ yields
\begin{align*}
v_m &=v_{\ell+1}+\sum_{k=\ell+1}^{m-1}\alpha_k\left[\phi-F^{d_k}(v_k)\right] +R_m-R_{\ell+1} \le\overline v+\frac{r_v}{4}+\frac{r_v}{8} <\overline v+r_v.
\end{align*}
Thus, $v_m\le\overline v+r_v$ for every $m\ge\tau$. Reversing the inequalities and using the lower tail drift gives another finite random index after which $v_m\ge\underline v-r_v$. Consequently, there is a finite random $N_2$ such that $\underline v-r_v \le v_n\le\overline v+r_v,~ n\ge N_2$.
Every $q_n=\operatorname{VaR}^{d_n}$ lies in $[v_{\min},v_{\max}]$, which is contained in $[\underline v,\overline v]$. Therefore, for $n\ge N_2$, the point $\widetilde v_n$ between $v_n$ and $q_n$ lies in $[\underline v-r_v,\overline v+r_v]$. Assumption~\ref{A.function} gives the deterministic bound $\varepsilon_f\le c_n\le M_f$ for all such $n$. Subtracting $q_{n+1}$ from both sides of~\eqref{Eq.VaR_tracking_recursion_app} gives $e_{n+1} =(1-\alpha_nc_n)e_n+\alpha_n\eta_{n+1}-(q_{n+1}-q_n)$.
By~\eqref{Eq.VaR_target_increment_app} and $\gamma_n=o(\alpha_n)$, $(|q_{n+1}-q_n|)/\alpha_n \le C_q (\gamma_n/\alpha_n) \longrightarrow0$.
All conditions of Lemma~\ref{lemma_scalar_tracking_app} are now satisfied with $a_n=\alpha_n$ and $r_n=-(q_{n+1}-q_n)$. Therefore, $e_n\to0$ almost surely, which is the desired result.
\end{proof}

\begin{proof}[Proof of Proposition~\ref{prop_weight_certificate}]
For a $Q$-function $Q$, write $h_Q(s)=\min_{a\in\mathcal A(s)}Q(s,a)$. We first compare the two minimum functions. Let $a'$ be a minimizing action of $Q'$ at state $s$. Then $h_Q(s)-h_{Q'}(s) \le Q(s,a')-Q'(s,a') \le w(s)\|Q-Q'\|_w$.
The first inequality follows because $h_Q(s)$ is no larger than $Q(s,a')$. The second follows from $w(s,a')=w(s)$ and the definition of the weighted maximum norm. Interchanging $Q$ and $Q'$ gives the reverse bound, and hence $|h_Q(s)-h_{Q'}(s)| \le w(s)\|Q-Q'\|_w$.
The transformed one stage cost is the same in $\mathcal T_vQ$ and $\mathcal T_vQ'$, so it cancels. For every feasible pair $(s,a)$, separate the reference state term to obtain
\begin{align*}
&\left(\mathcal T_vQ-\mathcal T_vQ'\right)(s,a)\\
&\quad=\sum_{y\ne s^{(0)}}p(y\mid s,a)\left[h_Q(y)-h_{Q'}(y)\right]
-\left[1-p\left(s^{(0)}\mid s,a\right)\right]
\left[h_Q\left(s^{(0)}\right)-h_{Q'}\left(s^{(0)}\right)\right].
\end{align*}
The equality follows by combining the term $p(s^{(0)}\mid s,a)[h_Q(s^{(0)})-h_{Q'}(s^{(0)})]$ with the normalization term $-[h_Q(s^{(0)})-h_{Q'}(s^{(0)})]$. Applying the triangle inequality and the preceding bound gives
\begin{align*}
&\left|\left(\mathcal T_vQ-\mathcal T_vQ'\right)(s,a)\right| \le\left\{\sum_{y\ne s^{(0)}}p(y\mid s,a)w(y) +\left[1-p\left(s^{(0)}\mid s,a\right)\right]w\left(s^{(0)}\right)\right\} \left\|Q-Q'\right\|_w\\
&\quad\le\rho_Qw(s)\left\|Q-Q'\right\|_w.
\end{align*}
The last inequality is the certificate in the proposition. Since $w(s,a)=w(s)$, division by $w(s,a)$ and maximization over all feasible pairs yield $\|\mathcal T_vQ-\mathcal T_vQ'\|_w \le\rho_Q\|Q-Q'\|_w$.
The argument does not depend on $v$, because the transformed cost cancels from the difference. This proves Assumption~\ref{A.Qcontraction}. For a fixed value of $\rho_Q$, the certificate inequalities are linear and homogeneous in the positive weights. Feasibility can therefore be checked by a linear program with $w(s)\ge1$ for every state $s$, for example by minimizing $\sum_{s\in\mathcal S}w(s)$, and the smallest certified value can be found by bisection on $\rho_Q$. If a positive certificate is obtained under another temporary scale, define $c_w:=\min_{s\in\mathcal S}w(s)>0$ and replace every weight by $w(s)/c_w$. The certificate is unchanged because both sides of every inequality are divided by the same positive number. The rescaled weights satisfy $\min_{s\in\mathcal S}w(s)=1$ and are the weights used to compute $w_{\max}$ and the projection bound $Q_{\max}$. We now prove the necessary occupancy restriction using the global contraction condition, without assuming that the sufficient state weight certificate was used. For the positive feasible pair weights in Assumption~\ref{A.Qcontraction}, set $\bar w(s):=\max_{a\in\mathcal A(s)}w(s,a)$ and select an action $\bar a(s)$ attaining this maximum. Define $Q'(s,\bar a(s))=0$ and $Q'(s,a)=2w_{\max}$ for every other feasible action. Obtain $Q$ by replacing the selected entry by $\bar w(s)$ when $s\ne s^{(0)}$ and by $-\bar w(s^{(0)})$ when $s=s^{(0)}$, leaving the other entries unchanged. The selected action is the unique minimum of both functions, because $2w_{\max}>\bar w(s)$. Also, $\|Q-Q'\|_w=1$, since the only nonzero difference at each state has magnitude $w(s,\bar a(s))$.

For this particular pair of $Q$-functions, the changes in the state minima have opposite signs at the reference state and elsewhere. Substitution into the operator gives, for every feasible pair,
\begin{align*}
\left(\mathcal T_vQ-\mathcal T_vQ'\right)(s,a)
&=\sum_{y\ne s^{(0)}}p(y\mid s,a)\bar w(y)
 +\left[1-p\left(s^{(0)}\mid s,a\right)\right]\bar w\left(s^{(0)}\right)\le\rho_Qw(s,a).
\end{align*}
The equality follows by collecting the reference transition and the subtracted reference minimum. Every term on its right is nonnegative. The inequality is therefore the component bound obtained directly from global contraction and $\|Q-Q'\|_w=1$.

Fix any stationary randomized policy $d$. Multiply this component inequality by $\mu^d(s)d(s,a)$ and sum. Stationarity gives $\sum_s\mu^d(s)\sum_a d(s,a)p(y\mid s,a)=\mu^d(y)$ for every $y$, and thus
\begin{align*}
\sum_{y\ne s^{(0)}}\mu^d(y)\bar w(y) +\left[1-\mu^d\left(s^{(0)}\right)\right]\bar w\left(s^{(0)}\right) &\le\rho_Q\sum_s\mu^d(s)\sum_{a\in\mathcal A(s)}d(s,a)w(s,a) \le\rho_Q\sum_s\mu^d(s)\bar w(s).
\end{align*}
The first inequality averages the preceding component bounds. The second uses $w(s,a)\le\bar w(s)$ and $\sum_a d(s,a)=1$. The left hand side equals $\sum_s\mu^d(s)\bar w(s)+[1-2\mu^d(s^{(0)})]\bar w(s^{(0)})$. Moving the full weighted sum to the right and changing signs yields $[2\mu^d(s^{(0)})-1]\bar w(s^{(0)}) \ge(1-\rho_Q)\sum_s\mu^d(s)\bar w(s)>0$.
The strict inequality uses $\rho_Q<1$, positive weights, and the fact that $\mu^d$ is a probability distribution. Since $\bar w(s^{(0)})>0$, this proves~\eqref{Eq.reference_occupancy_necessary} for every stationary randomized policy. The preceding component inequality also gives a state weight certificate: replace its upper bound $\rho_Qw(s,a)$ by $\rho_Q\bar w(s)$. Thus, allowing weights to depend on the action does not enlarge the class of models for which some global certificate exists.

It remains to prove the converse under Assumption~\ref{A.mixing}. Assume~\eqref{Eq.reference_occupancy_necessary} for every stationary policy. When there is only one state, the continuation minimum and the reference minimum cancel, so the operator is independent of $Q$ and the assertion is immediate. Suppose there are at least two states. Starting at $s\ne s^{(0)}$, let $\tau:=\inf\{n\ge0:s_n=s^{(0)}\}$ and define $H_j(s)$ as the largest expected value of $\min\{\tau,j\}$ over admissible policies. Set $H_j(s^{(0)})=0$. The conditional minorization gives, for every admissible policy, $\PP(\tau>km_p\mid s_0=s) \le[1-\delta_p\nu(s^{(0)})]^k$, $\E[\tau\mid s_0=s] =\sum_{n=0}^{\infty}\PP(\tau>n\mid s_0=s) \le m_p/(\delta_p\nu(s^{(0)}))$.
Here the first inequality follows by conditioning successively at the block endpoints: failure to visit the reference state during a block implies failure to be there at its endpoint. The second equality is the tail sum formula for a nonnegative integer valued random variable, and the last bound groups the summands into blocks of length $m_p$. These bounds hold uniformly over the policies.

Finite horizon conditioning on the first transition gives $H_0(s)=0$ and, for $s\ne s^{(0)}$, $H_{j+1}(s) =1+\max_{a\in\mathcal A(s)}\sum_{y\ne s^{(0)}}p(y\mid s,a)H_j(y)$.
The sequence increases with $j$ and is bounded by the preceding mean hitting time bound. Hence it converges to a finite $H(s)\ge1$ at each nonreference state. Since the action sets are finite, passage to the limit through the maximum is valid and yields
\begin{align}\label{Eq.hitting_certificate_app}
H(s)&=1+\max_{a\in\mathcal A(s)}\sum_{y\ne s^{(0)}}p(y\mid s,a)H(y),\qquad s\ne s^{(0)}.
\end{align}
Choose a deterministic stationary policy $\widehat d$ attaining these maxima, and choose its reference action to maximize $\sum_{y\ne s^{(0)}}p(y\mid s^{(0)},a)H(y)$. Denote this last maximum by $R$. The stationary equations for $\mu^{\widehat d}$ imply
\begin{align*}
\sum_{y\ne s^{(0)}}\mu^{\widehat d}(y)H(y)
&=\mu^{\widehat d}\left(s^{(0)}\right)R
 +\sum_{s\ne s^{(0)}}\mu^{\widehat d}(s)\sum_{y\ne s^{(0)}}P^{\widehat d}(s,y)H(y)\\
&=\mu^{\widehat d}\left(s^{(0)}\right)R
 +\sum_{s\ne s^{(0)}}\mu^{\widehat d}(s)\left[H(s)-1\right].
\end{align*}
The first equality substitutes stationarity and separates the reference state. The second uses the maximizing actions and~\eqref{Eq.hitting_certificate_app}. Cancelling the common weighted sums gives $1-\mu^{\widehat d}(s^{(0)})=\mu^{\widehat d}(s^{(0)})R$. Therefore $\mu^{\widehat d}(s^{(0)})=1/(1+R)$, and the assumed occupancy condition gives $R<1$.

Choose $\eta>0$ so small that $(1+\eta)R<1$. Define $w(s^{(0)})=1$ and $w(s)=(1+\eta)H(s)-1$ for $s\ne s^{(0)}$. These weights are positive because $H(s)\ge1$. For a nonreference state and every feasible action, the certificate expression satisfies
\begin{align*}
&\sum_{y\ne s^{(0)}}p(y\mid s,a)w(y) +\left[1-p\left(s^{(0)}\mid s,a\right)\right]w\left(s^{(0)}\right) =(1+\eta)\sum_{y\ne s^{(0)}}p(y\mid s,a)H(y)\\
&\quad\le(1+\eta)\left[H(s)-1\right]=w(s)-\eta<w(s).
\end{align*}
The equality cancels $-\sum_{y\ne s^{(0)}}p(y\mid s,a)$ with $1-p(s^{(0)}\mid s,a)$. The inequality follows from~\eqref{Eq.hitting_certificate_app}. At the reference state the same expression is at most $(1+\eta)R<1=w(s^{(0)})$. Divide each strict inequality by its positive state weight. There are finitely many feasible pairs, so the largest ratio is strictly below one and can be bounded by a common $\rho_Q\in(0,1)$. The required certificate follows. Finally, divide every weight by its minimum before computing $w_{\max}$ or $Q_{\max}$; this common rescaling does not change any certificate ratio. This proves the converse and completes the characterization.
\end{proof}

\begin{proof}[Proof of Proposition~\ref{prop_gain_certificate}]
Fix a feasible pair $(s,a)$ and write $I_k(s,a)=\I\{s_k=s,a_k=a\}$. We first show that every sufficiently long calendar interval contains enough visits to this pair. Let $K>0$ be a constant to be selected below and define $L_n:=\lceil K(n+1)^{p_\beta+p_\epsilon}\rceil$. Because $p_\beta+p_\epsilon<p_\alpha<1$, we have $L_n/n\to0$. Hence, for all sufficiently large $n$, $L_n\le n$ and every iteration between $n$ and $n+L_n$ is no larger than $2n$.

Consider the block times $t_j=n+j(m_p+1)$, where $j=1,\ldots,J_n$ and $J_n=\lfloor L_n/(m_p+1)\rfloor$. Also set $r_j:=t_j-m_p=t_{j-1}+1$, with $t_0:=n$. The observations at $t_1,\ldots,t_{j-1}$ are contained in $\mathscr F_{r_j}$. By Assumption~\ref{A.mixing}, $\PP(s_{t_j}=s\mid\mathscr F_{r_j}) =\PP(s_{r_j+m_p}=s\mid\mathscr F_{r_j}) \ge\delta_p\nu(s)$.
Conditional on $\mathscr F_{t_j}$, the action is sampled from $d_{t_j}(s_{t_j})$. The statewise policy projection in Algorithm~\ref{GOODALGORITHM} gives $d_{t_j}(s,a)\ge\epsilon_{t_j-1}$. Applying the tower property first at $\mathscr F_{t_j}$ and then at $\mathscr F_{r_j}$ yields
\begin{align*}
\E\left[I_{t_j}(s,a)\mid\mathscr F_{r_j}\right] &=\E\left[\E\left[I_{t_j}(s,a)\mid\mathscr F_{t_j}\right]\middle|\mathscr F_{r_j}\right] =\E\left[\I\{s_{t_j}=s\}d_{t_j}(s,a)\mid\mathscr F_{r_j}\right]\\
&\ge\epsilon_{t_j-1}\PP\left(s_{t_j}=s\mid\mathscr F_{r_j}\right) \ge\delta_p\nu(s)\epsilon_{t_j-1}.
\end{align*}
The first equality is the tower property, the second follows from the action sampling rule, the first inequality uses the exploration lower bound, and the last inequality uses the common minorization. Since $t_j\le n+L_n\le2n$ for all sufficiently large $n$ and $\epsilon_k=\epsilon_0/(k+1)^{p_\epsilon}$ is decreasing, $\E[I_{t_j}(s,a)\mid\mathscr F_{r_j}] \ge q_n, ~ q_n:= (\delta_p\nu(s)\epsilon_0)/(2n+1)^{p_\epsilon}$.
We next obtain a probability bound that is uniform in the dependence among the block observations. Let $\lambda>0$. Since $I_{t_j}(s,a)$ takes values zero and one,
\begin{align*}
\E\left[\exp\{-\lambda I_{t_j}(s,a)\}\mid\mathscr F_{r_j}\right] &=1-\left(1-e^{-\lambda}\right) \E\left[I_{t_j}(s,a)\mid\mathscr F_{r_j}\right] \le\exp\left\{-\left(1-e^{-\lambda}\right)q_n\right\}.
\end{align*}
The equality follows by considering the two values of the indicator. The inequality uses $1-x\le e^{-x}$ and the preceding conditional probability lower bound. To justify the iteration, observe that $I_{t_1}(s,a),\ldots,I_{t_{j-1}}(s,a)$ are $\mathscr F_{r_j}$ measurable. Hence,
\begin{align*}
&\E\left[\exp\left\{-\lambda\sum_{\ell=1}^{j}I_{t_\ell}(s,a)\right\}\middle|\mathscr F_n\right] =\E\left[\exp\left\{-\lambda\sum_{\ell=1}^{j-1}I_{t_\ell}(s,a)\right\} \E\left[\exp\{-\lambda I_{t_j}(s,a)\}\mid\mathscr F_{r_j}\right]\middle|\mathscr F_n\right]\\
&\quad\le\exp\left\{-\left(1-e^{-\lambda}\right)q_n\right\} \E\left[\exp\left\{-\lambda\sum_{\ell=1}^{j-1}I_{t_\ell}(s,a)\right\}\middle|\mathscr F_n\right].
\end{align*}
Applying this inequality successively for $j=J_n,J_n-1,\ldots,1$ gives $\E[\exp\{-\lambda\sum_{j=1}^{J_n}I_{t_j}(s,a)\}\mid \mathscr F_n] \le\exp\{-(1-e^{-\lambda})J_nq_n\}$.
By Markov's inequality, $\PP(\sum_{j=1}^{J_n}I_{t_j}(s,a)\le ((J_nq_n)/2) \mid \mathscr F_n) \le\exp\{ (\lambda J_nq_n)/2 -(1-e^{-\lambda})J_nq_n\}$.
Take $\lambda=\log2$. Then $1-e^{-\lambda}=1/2$, and the exponent is $-(1-\log2)J_nq_n/2$. For all sufficiently large $n$, $J_n\ge L_n/[2(m_p+1)]$, so there is a constant $c_1>0$, independent of $n$ and $K$, such that $J_nq_n \ge c_1K(n+1)^{p_\beta}$.
The power on the right follows by subtracting $p_\epsilon$ from the calendar length power $p_\beta+p_\epsilon$. Consequently, $\sum_{n=1}^{\infty} \PP(\sum_{j=1}^{J_n}I_{t_j}(s,a)\le (J_nq_n)/2 ) <\infty$,
because $\sum_n\exp\{-cK n^{p_\beta}\}<\infty$ for every $c,K>0$. The Borel Cantelli lemma now implies that, on an event of probability one, there is a finite random index after which $\sum_{j=1}^{J_n}I_{t_j}(s,a) \ge (J_nq_n)/2 \ge c_2K(n+1)^{p_\beta}$
for a constant $c_2>0$.

Every block visit counted above occurs no later than $n+L_n\le2n$. Since $N_k(s,a)\le k+1$ and $f_\beta$ is decreasing, the Q gain used at each such visit satisfies $\beta_k(s,a) = \beta_0/(N_k(s,a)+1)^{p_\beta} \ge \beta_0/(2n+2)^{p_\beta}$.
Multiplying this lower bound by the number of block visits gives a constant $c_3>0$ such that, eventually almost surely, $\sum_{k=n}^{n+L_n}\beta_k(s,a)I_k(s,a) \ge c_3K$.
Choose $K$ so that $c_3K>T$. Then every sufficiently large interval $[n,n+L_n]$ accumulates more than $T$ units of the local Q gain sequence. Therefore, if $m\ge n$ and $\sum_{k=n}^{m}\beta_k(s,a)I_k(s,a) \le T$,
we must have $m<n+L_n$ for all sufficiently large $n$. Since $\alpha_k$ is decreasing,
\begin{align*}
\sup\left\{\sum_{k=n}^{m}\alpha_k:m\ge n,~\sum_{k=n}^{m}\beta_k(s,a)I_k(s,a)\le T\right\} &\le\left(L_n+1\right)\alpha_n =O\left(n^{p_\beta+p_\epsilon-p_\alpha}\right) \longrightarrow0.
\end{align*}
The first inequality uses the calendar length bound, the equality follows from the polynomial forms of $L_n$ and $\alpha_n$, and the convergence follows from $p_\beta+p_\epsilon<p_\alpha$.

Finally, take $T=1$. On the probability one event obtained above, choose a sufficiently large integer $n_1$ and define recursively $n_{j+1}:=n_j+L_{n_j}+1$. The intervals $[n_j,n_j+L_{n_j}]$ are disjoint, and every one of them contributes more than one unit of local Q gain. Hence, for every $J\ge1$, $\sum_{k=n_1}^{n_J+L_{n_J}}\beta_k(s,a)I_k(s,a)>J$, which implies $\sum_n\beta_n(s,a)I_n(s,a)=\infty$ almost surely. First carry out the argument for every positive integer $T$ and take the countable intersection of the resulting probability one events. For any real $T>0$, the supremum in~\ref{A.Qgain} is bounded by the supremum for the integer $\lceil T\rceil$, so the conclusion holds on the same event. Since the set of feasible state and action pairs is finite, this event can also be intersected over all pairs. Both parts of Assumption~\ref{A.Qgain} follow. This completes the proof.
\end{proof}

\begin{proof}[Proof of Theorem~\ref{th_Q}]
We first verify that the partial sums in~\eqref{Eq.Q_martingale_definition_app} converge. When $I_n(s,a)=1$, the conditional expectation of $\zeta_{n+1}(s,a)$ given $\mathscr G_n$ is zero. The cost term is centered because $\E[\widetilde C(v_n,s,a)\mid\mathscr G_n]=\tilde c(v_n,s,a)$, and the next state term is centered by the conditional transition distribution. Since $\beta_n(s,a)I_n(s,a)$ is $\mathscr G_n$ measurable, the tower property gives $\E[\beta_n(s,a)I_n(s,a)\zeta_{n+1}(s,a)\mid\mathscr F_n] =0$.
Thus, $\{M_m(s,a)\}$ is a martingale for every feasible state and action pair.

We next bound its second moment without treating the random local gain as a deterministic quantity. The Q projection gives $|Q_n(s,a)|\le Q_{\max}$. Moreover, $|\widetilde C(v_n,s,a)| \le |v_n|+ (|C(s,a)|+|v_n|)/(1-\phi) \le C_0(1+|v_n|+|C(s,a)|)$
for a finite constant $C_0$. Jensen's inequality gives $|\tilde c(v_n,s,a)|^2\le\E[\widetilde C^2(v_n,s,a)\mid\mathscr G_n]$, and the centered next state term is bounded by $2Q_{\max}$. Hence, for a finite deterministic constant $C_1$,
\begin{align}\label{Eq.Q_noise_conditional_bound_app}
\E\left[\zeta_{n+1}^2(s,a)\mid\mathscr G_n\right]I_n(s,a)
&\le C_1I_n(s,a)\left(1+v_n^2+\E[C^2(s,a)]\right).
\end{align}
The inequality follows from $(x_1+\cdots+x_4)^2\le4(x_1^2+\cdots+x_4^2)$, the preceding bound on $\widetilde C$, Jensen's inequality, and the bounded Q terms.

Orthogonality of the martingale differences and~\eqref{Eq.Q_noise_conditional_bound_app} yield
\begin{align*}
\E[M_m^2(s,a)] &=\sum_{n=0}^{m-1}\E\left[\beta_n^2(s,a)I_n(s,a)\zeta_{n+1}^2(s,a)\right]\\
&\le C_1\E\left[\left(1+\sup_{k\ge0}v_k^2\right)\sum_{n=0}^{m-1}\beta_n^2(s,a)I_n(s,a)\right] +C_1\E[C^2(s,a)]\E\left[\sum_{n=0}^{m-1}\beta_n^2(s,a)I_n(s,a)\right].
\end{align*}
The first equality follows because cross products of distinct martingale increments have zero expectation. The first inequality uses $v_n^2\le\sup_kv_k^2$ and separates the deterministic second moment of the component cost. Along every sample path, the gains used at successive visits to $(s,a)$ are distinct members of the sequence $\{f_\beta(j)\}$. Therefore,
\begin{align}\label{Eq.random_local_gain_sum_app}
\sum_{n=0}^{\infty}\beta_n^2(s,a)I_n(s,a)
&\le\sum_{j=0}^{\infty}f_\beta^2(j)
<\infty
\qquad\text{a.s.}
\end{align}
Substituting this deterministic upper bound and using Lemma~\ref{lemma_v_bound} gives
\begin{align*}
\sup_m\E[M_m^2(s,a)] &\le C_2\left(1+\E\left[\sup_{k\ge0}v_k^2\right]+\E[C^2(s,a)]\right)\sum_{j=0}^{\infty}f_\beta^2(j) <\infty.
\end{align*}
Thus, every $M_m(s,a)$ converges almost surely and in $L^2$. Because the number of state and action pairs is finite,
\begin{align}\label{Eq.Q_martingale_tail_app}
\max_{s\in\mathcal S,~a\in\mathcal A(s)}\sup_{m\ge n}|M_m(s,a)-M_n(s,a)|
&\longrightarrow0
\qquad\text{a.s.}
\end{align}
The convergence follows from the Cauchy property of each convergent martingale and finiteness of the component set.

We next verify the fixed point bound required by Lemma~\ref{lemma_projected_cycle_app}. For $v\in[v_{\min}-r_v,v_{\max}+r_v]$, the fixed point equation and Assumption~\ref{A.Qcontraction} give
\begin{align*}
\|Q(v)\|_w &=\|\mathcal T_vQ(v)\|_w \le\|\mathcal T_vQ(v)-\mathcal T_v0\|_w+\|\mathcal T_v0\|_w\\
&\le\rho_Q\|Q(v)\|_w+\|\mathcal T_v0\|_w.
\end{align*}
The first inequality is the triangle inequality, and the second is the weighted contraction. Moving the first term on the last line to the left gives $\|Q(v)\|_w \le (\|\mathcal T_v0\|_w)/(1-\rho_Q)$.
Since $\mathcal T_v0(s,a)=\tilde c(v,s,a)$, $w(s,a)\ge1$, and $|\tilde c(v,s,a)| \le |v|+ (\E[|C(s,a)|]+|v|)/(1-\phi) = (\E[|C(s,a)|]+(2-\phi)|v|)/(1-\phi)$,
Assumption~\ref{A.Qprojection} implies
\begin{align}\label{Eq.Q_target_margin_app}
\max_{s\in\mathcal S,~a\in\mathcal A(s)}|Q(v)(s,a)|
&\le w_{\max}\|Q(v)\|_w
<Q_{\max}-1.
\end{align}
Thus, every fixed point used below has a uniform distance from the projection boundary.

Theorem~\ref{th_v} and $\operatorname{VaR}^{d_n}\in[v_{\min},v_{\max}]$ imply that, on an event of probability one, there is a finite $N_1$ such that $v_n\in[v_{\min}-r_v,v_{\max}+r_v]$ for every $n\ge N_1$. Fix a sample path in the intersection of this event and the event in~\eqref{Eq.Q_martingale_tail_app}.

Fix $T>0$. Starting from $n_0\ge N_1$, define $n_{r+1}$ to be the first integer larger than $n_r$ such that $\sum_{k=n_r}^{n_{r+1}-1}\beta_k(s,a)I_k(s,a) \ge T, ~ s\in\mathcal S,~a\in\mathcal A(s)$.
Assumption~\ref{A.Qgain} makes every endpoint finite. Since $n_{r+1}$ is the first common endpoint, the pair updated at the final iteration has accumulated a local gain no larger than $T+\max_{k\ge n_r,s,a}\beta_k(s,a)$. The maximum tends to zero. If the following conclusion failed, a subsequence of cycles would have a further subsequence with the same final pair because there are finitely many pairs. The second part of Assumption~\ref{A.Qgain}, applied to this pair with upper bound $T+1$, would then give a contradiction. Therefore, $\sum_{k=n_r}^{n_{r+1}-1}\alpha_k \longrightarrow0$.
Since $|v_{k+1}-v_k|\le\alpha_k$, it follows that
\begin{align}\label{Eq.v_cycle_app}
\Delta_v(n_r,n_{r+1})
&\le\sum_{k=n_r}^{n_{r+1}-1}\alpha_k
\longrightarrow0.
\end{align}
Equation~\eqref{Eq.Q_martingale_tail_app} also gives
\begin{align}\label{Eq.M_cycle_app}
\Delta_M(n_r,n_{r+1})
&\longrightarrow0.
\end{align}
In particular, $\Delta_M(n_r,n_{r+1})<1$ for all sufficiently large $r$, so Lemma~\ref{lemma_projected_cycle_app} applies.

Let $R_r :=\|Q_{n_r}-Q(v_{n_r})\|_w$.
The endpoint bound~\eqref{Eq.projected_cycle_endpoint_app} and the Lipschitz bound~\eqref{Eq.Q_Lipschitz_app} give $R_{r+1} \le\|Q_{n_{r+1}}-Q(v_{n_r})\|_w+\|Q(v_{n_r})-Q(v_{n_{r+1}})\|_w \le\kappa_TR_r+\varepsilon_r$,
where
\begin{align*}
\varepsilon_r &:=\Delta_M(n_r,n_{r+1}) +\frac{(1+\rho_Q)\Delta_M(n_r,n_{r+1})+L_c\Delta_v(n_r,n_{r+1})}{1-\rho_Q} +\frac{L_c}{1-\rho_Q}\Delta_v(n_r,n_{r+1}).
\end{align*}
Equations~\eqref{Eq.v_cycle_app} and~\eqref{Eq.M_cycle_app} imply $\varepsilon_r\to0$, whereas $0<\kappa_T<1$. Repeated substitution gives, for $r>r_0$, $R_r \le\kappa_T^{r-r_0}R_{r_0}+\sum_{j=r_0}^{r-1}\kappa_T^{r-1-j}\varepsilon_j$.
The first term tends to zero. For the second term, fix $\eta>0$ and choose $r_0$ so large that $\varepsilon_j\le\eta(1-\kappa_T)$ for all $j\ge r_0$. Then $\sum_{j=r_0}^{r-1}\kappa_T^{r-1-j}\varepsilon_j \le\eta(1-\kappa_T)\sum_{j=r_0}^{r-1}\kappa_T^{r-1-j} \le\eta$.
Because $\eta$ is arbitrary, $R_r\to0$.

We now pass from the cycle endpoints to every calendar iteration. For $n_r\le n\le n_{r+1}$, the uniform cycle bound~\eqref{Eq.projected_cycle_uniform_app} and~\eqref{Eq.Q_Lipschitz_app} give
\begin{align*}
\|Q_n-Q(v_n)\|_w
&\le\|Q_n-Q(v_{n_r})\|_w+\|Q(v_{n_r})-Q(v_n)\|_w\\
&\le R_r+\Delta_M(n_r,n_{r+1})\\
&\quad+\frac{(1+\rho_Q)\Delta_M(n_r,n_{r+1})+L_c\Delta_v(n_r,n_{r+1})}{1-\rho_Q}\\
&\quad+\frac{L_c}{1-\rho_Q}\Delta_v(n_r,n_{r+1}).
\end{align*}
Every term on the last three lines tends to zero. Hence, $\|Q_n-Q(v_n)\|_w\to0$ almost surely. Finally, Theorem~\ref{th_v} and~\eqref{Eq.Q_Lipschitz_app} imply
\begin{align*}
\|Q(v_n)-Q^{d_n}\|_w &=\left\|Q(v_n)-Q\left(\operatorname{VaR}^{d_n}\right)\right\|_w \le\frac{L_c}{1-\rho_Q}\left|v_n-\operatorname{VaR}^{d_n}\right| \longrightarrow0.
\end{align*}
The triangle inequality proves $\|Q_n-Q^{d_n}\|_w\to0$. Since $\|Q\|_\infty\le w_{\max}\|Q\|_w$, the ordinary maximum norm conclusion follows as well.

For completeness, we verify that the Q projection is eventually inactive. By~\eqref{Eq.Q_target_margin_app} and the convergence just proved, every Q component is eventually at least three quarters of a unit from the boundary. The Q terms in one unprojected increment are bounded and their local gain tends to zero. For the random cost term, Markov's inequality gives, for every fixed $\eta>0$,
\begin{align*}
&\sum_n\PP\left(\beta_n(s,a)I_n(s,a)|\widetilde C(v_n,s,a)|>\eta\right) \le\frac{1}{\eta^2}\E\left[\sum_n\beta_n^2(s,a)I_n(s,a)\widetilde C^2(v_n,s,a)\right]\\
&\quad\le C_3\left(1+\E\left[\sup_n v_n^2\right]+\E[C^2(s,a)]\right)\sum_{j=0}^{\infty}f_\beta^2(j)\\
&\quad<\infty.
\end{align*}
The first inequality is Markov's inequality applied term by term. The second uses the same pathwise local gain bound as~\eqref{Eq.random_local_gain_sum_app}, and the last line follows from Lemma~\ref{lemma_v_bound}. The Borel--Cantelli lemma therefore gives $\beta_n(s,a)I_n(s,a)|\widetilde C(v_n,s,a)|\to0$ almost surely. Since there are finitely many components, the complete unprojected Q increment tends to zero. It is eventually smaller than one half, while the current Q value is at least three quarters of a unit from the boundary. Hence, the projection does not activate after an almost surely finite iteration. This final observation is not needed for the convergence argument, but it confirms that the projection does not alter the limiting recursion. This completes the proof.
\end{proof}

\begin{proof}[Proof of Proposition~\ref{prop_performance_difference}]
Let $q^*=\operatorname{VaR}^{d^*}$, $Q^*=Q(q^*)$, and $h^*(s):=\min_{b\in\mathcal A(s)}Q^*(s,b)$. The fixed point equation $Q^*=\mathcal T_{q^*}Q^*$ can be written as
\begin{align}\label{Eq.performance_poisson_app}
Q^*(s,a)+g^*
&=\tilde c(q^*,s,a)+\sum_{s'\in\mathcal S}p(s'\mid s,a)h^*(s'),
\end{align}
where $g^*:=\min_bQ^*(s^{(0)},b)$. Since Assumption~\ref{A.policylimit} makes $a^*(s)$ the unique minimizer, $Q^*(s,a^*(s))=h^*(s)$.

We first identify $g^*$. Multiply~\eqref{Eq.performance_poisson_app} by $\mu^{d^*}(s)d^*(s,a)$ and sum over the state and action pairs. The left hand side becomes $\sum_s\mu^{d^*}(s)h^*(s)+g^*$,
because $d^*(s)$ selects $a^*(s)$. The continuation term on the right hand side is $\sum_{s\in\mathcal S}\mu^{d^*}(s)\sum_{a\in\mathcal A(s)}d^*(s,a)\sum_{s'\in\mathcal S}p(s'\mid s,a)h^*(s') =\sum_{s'\in\mathcal S}\mu^{d^*}(s')h^*(s')$.
The equality follows from stationarity of $\mu^{d^*}$ under $P^{d^*}$. The two terms containing $h^*$ therefore cancel, which gives
\begin{align*}
g^* &=\sum_{s\in\mathcal S}\mu^{d^*}(s)\sum_{a\in\mathcal A(s)}d^*(s,a)\tilde c(q^*,s,a) =q^*+\frac{1}{1-\phi}\E\left[(C^{d^*}-q^*)^+\right] =\operatorname{CVaR}^{d^*}.
\end{align*}
The second equality is the definition of the stationary cost mixture, and the last equality follows because $q^*$ is the VaR of $C^{d^*}$.

Now fix an arbitrary stationary policy $d$. Multiply~\eqref{Eq.performance_poisson_app} by $\mu^d(s)d(s,a)$ and sum. Stationarity of $\mu^d$ under $P^d$ again cancels the continuation term with $\sum_s\mu^d(s)h^*(s)$. Consequently,
\begin{align}
&\sum_{s\in\mathcal S}\mu^d(s)\sum_{a\in\mathcal A(s)}d(s,a)\tilde c(q^*,s,a)-\operatorname{CVaR}^{d^*} =\sum_{s\in\mathcal S}\mu^d(s)\sum_{a\in\mathcal A(s)}d(s,a)\left[Q^*(s,a)-h^*(s)\right].\label{Eq.fixed_var_performance_app}
\end{align}
This is an exact identity. No limiting argument or first order expansion has been used.

Recall the function $\Psi_d$ defined before~\eqref{Eq.threshold_control_bridge}. By the stationary cost mixture, it also satisfies $\Psi_d(y)=y+\E[(C^d-y)^+]/(1-\phi)$.
Then the left hand side of~\eqref{Eq.fixed_var_performance_app} is $\Psi_d(q^*)-\operatorname{CVaR}^{d^*}$, whereas $\operatorname{CVaR}^d=\min_y\Psi_d(y)=\Psi_d(q_d)$ with $q_d=\operatorname{VaR}^d$. Thus,
\begin{align}
\operatorname{CVaR}^d-\operatorname{CVaR}^{d^*} &=\sum_{s\in\mathcal S}\mu^d(s)\sum_{a\in\mathcal A(s)}d(s,a)\left[Q^*(s,a)-h^*(s)\right] -\left[\Psi_d(q^*)-\Psi_d(q_d)\right].\label{Eq.CVaR_performance_identity_app}
\end{align}
The identity follows by adding and subtracting $\Psi_d(q^*)$ and then using~\eqref{Eq.fixed_var_performance_app}.

We next bound the final bracket. Since the stationary cost distribution is continuous, $\Psi_d'(y) =1- (\PP(C^d>y))/(1-\phi) = (F^d(y)-\phi)/(1-\phi)$.
The first equality follows by differentiating the positive part expectation, and the second uses $\PP(C^d>y)=1-F^d(y)$. The density of $F^d$ is a convex combination of the component densities and is therefore no larger than $M_f$. Hence, $|\Psi_d'(y)-\Psi_d'(z)| \le (M_f/(1-\phi)) |y-z|$.
Also, $\Psi_d'(q_d)=0$. If $q^*\ge q_d$, the fundamental theorem of calculus gives
\begin{align*}
0\le\Psi_d(q^*)-\Psi_d(q_d) &=\int_{q_d}^{q^*}\left[\Psi_d'(y)-\Psi_d'(q_d)\right]dy \le\frac{M_f}{1-\phi}\int_{q_d}^{q^*}(y-q_d)dy =\frac{M_f}{2(1-\phi)}(q^*-q_d)^2.
\end{align*}
The first inequality follows because $q_d$ minimizes $\Psi_d$, the second follows from the Lipschitz bound on the derivative, and the last equality evaluates the elementary integral. The same bound follows when $q^*<q_d$ by integrating from $q^*$ to $q_d$. Therefore, $0\le\Psi_d(q^*)-\Psi_d(q_d) \le (M_f/(2(1-\phi))) |q^*-q_d|^2 \le C_{\rm p}\|d-d^*\|^2$,
where the last inequality follows from the VaR Lipschitz bound~\eqref{Eq.VaR_Lipschitz_app}, and one may take $C_{\rm p} = (M_f/(2(1-\phi))) ( (L_\mu+1)/\varepsilon_f )^2$.
Substituting this bound into~\eqref{Eq.CVaR_performance_identity_app} proves the inequality stated in the proposition.

It remains to show local optimality. Assumption~\ref{A.mixing}, applied to a stationary policy $d$, gives $(P^d)^{m_p}(s,s')\ge\delta_p\nu(s')$. Multiplying by $\mu^d(s)$ and summing over $s$ yields $\mu^d(s') =\sum_s\mu^d(s)(P^d)^{m_p}(s,s') \ge\delta_p\nu(s')$.
The equality uses stationarity, and the inequality uses the common minorization. Thus, with $\mu_{\min}:=\delta_p\min_s\nu(s)>0$, every stationary policy satisfies $\mu^d(s)\ge\mu_{\min}$.

For every state, $h^*(s)=Q^*(s,a^*(s))$. The strict action gap therefore gives
\begin{align*}
&\sum_{s\in\mathcal S}\mu^d(s)\sum_{a\in\mathcal A(s)}d(s,a)\left[Q^*(s,a)-h^*(s)\right] =\sum_{s\in\mathcal S}\mu^d(s)\sum_{a\in\mathcal A(s),~a\ne a^*(s)}d(s,a)\left[Q^*(s,a)-Q^*(s,a^*(s))\right]\\
&\quad\ge\mu_{\min}\Delta_Q\max_{s\in\mathcal S}\sum_{a\in\mathcal A(s),~a\ne a^*(s)}d(s,a) =\frac{\mu_{\min}\Delta_Q}{2}\|d-d^*\|.
\end{align*}
The equality removes the minimizing action, whose advantage is zero. The inequality uses the strict gap and the positive stationary probabilities. The final equality follows from $d^*(s)=\delta(a^*(s))$, which gives $\|d(s)-d^*(s)\|_1=2\sum_{a\ne a^*(s)}d(s,a)$.

Combining the last display with the proposition bound gives $\operatorname{CVaR}^d-\operatorname{CVaR}^{d^*} \ge ((\mu_{\min}\Delta_Q)/2) \|d-d^*\|-C_{\rm p}\|d-d^*\|^2$.
The right hand side is strictly positive whenever $0<\|d-d^*\|<\mu_{\min}\Delta_Q/(2C_{\rm p})$. Hence, every sufficiently close stationary randomized policy different from $d^*$ has a larger stationary one stage cost CVaR. This is the usual neighborhood definition of local optimality and covers all nearby policies directly, without restricting them to a prescribed radial representation. This completes the proof.
\end{proof}

\begin{proof}[Proof of Theorem~\ref{th_opt}]
Let $\Delta_Q$ be the strict action gap in Assumption~\ref{A.policylimit}. By Lemma~\ref{lemma_Q_continuity_app}, $d\mapsto Q^d$ is continuous in the ordinary maximum norm. Therefore, there is $r>0$ such that
\begin{align}\label{Eq.local_gap_app}
\|d-d^*\|\le r
&\quad\Longrightarrow\quad
Q^d(s,a)-Q^d(s,a^*(s))\ge\frac{\Delta_Q}{2}
\end{align}
for every state $s$ and every action $a\ne a^*(s)$. To verify this, choose $r$ so small that $\|Q^d-Q^{d^*}\|_\infty\le\Delta_Q/4$ whenever $\|d-d^*\|\le r$. Then
\begin{align*}
Q^d(s,a)-Q^d(s,a^*(s)) &\ge Q^{d^*}(s,a)-Q^{d^*}(s,a^*(s))-2\|Q^d-Q^{d^*}\|_\infty \ge\Delta_Q-\frac{\Delta_Q}{2} =\frac{\Delta_Q}{2}.
\end{align*}
The first inequality adds and subtracts the two corresponding values of $Q^{d^*}$, and the second uses the strict gap and the choice of $r$.

Let $\Omega_0$ be the event on which Theorem~\ref{th_Q} holds, the policy increment bound~\eqref{Eq.policy_increment_app} holds for all iterations, and all deterministic sequence limits used below hold. This event has probability one. Fix a sample path $\omega\in\Omega_0$.

Theorem~\ref{th_Q} gives $\|Q_{n+1}-Q^{d_{n+1}}\|_\infty\to0$. On the other hand,~\eqref{Eq.policy_increment_app} and Lemma~\ref{lemma_Q_continuity_app} imply $\|Q^{d_{n+1}}-Q^{d_n}\|_\infty \le L_d\|d_{n+1}-d_n\| \le L_dC_d\gamma_n \longrightarrow0$,
where $L_d$ is the Lipschitz constant of $d\mapsto Q^d$. The first inequality follows from Lemma~\ref{lemma_Q_continuity_app}, the second follows from~\eqref{Eq.policy_increment_app}, and the last relation follows from $\gamma_n\to0$. Therefore, $\|Q_{n+1}-Q^{d_n}\|_\infty \le\|Q_{n+1}-Q^{d_{n+1}}\|_\infty+\|Q^{d_{n+1}}-Q^{d_n}\|_\infty \longrightarrow0$.
The inequality is the triangle inequality. Consequently, there is a finite integer $N_Q(\omega)$ such that
\begin{align}\label{Eq.Q_identification_app}
\|Q_{n+1}-Q^{d_n}\|_\infty
&\le\frac{\Delta_Q}{8},
\qquad n\ge N_Q(\omega).
\end{align}
If $n\ge N_Q(\omega)$ and $\|d_n-d^*\|\le r$, then, for $a\ne a^*(s)$,
\begin{align*}
Q_{n+1}(s,a)-Q_{n+1}(s,a^*(s)) &\ge Q^{d_n}(s,a)-Q^{d_n}(s,a^*(s))-2\|Q_{n+1}-Q^{d_n}\|_\infty \ge\frac{\Delta_Q}{2}-\frac{\Delta_Q}{4}\\
&=\frac{\Delta_Q}{4}>0.
\end{align*}
The first inequality adds and subtracts the two values of $Q^{d_n}$, and the second uses~\eqref{Eq.local_gap_app} and~\eqref{Eq.Q_identification_app}. Thus, Procedure~\eqref{Eq.d} selects $a^*(s)$ at every state whenever the policy lies in the radius $r$ neighborhood after time $N_Q(\omega)$.

We next analyze the policy projection. Fix a state $s$, let $m=|\mathcal A(s)|$, and write $R_n(s):=\sum_{a\ne a^*(s)}d_n(s,a)$. Since $d^*(s)$ is a vertex of the simplex,
\begin{align}\label{Eq.policy_distance_mass_app}
\|d_n(s)-d^*(s)\|_1
&=2R_n(s).
\end{align}
When the selected action is $a^*(s)$, let $x_n(s):=d_n(s)+\gamma_n[\delta(a^*(s))-d_n(s)]$ be the vector before projection. Thus, $x_n(s,a)=(1-\gamma_n)d_n(s,a)$ for $a\ne a^*(s)$. Since $d_n(s,a)\ge\epsilon_{n-1}$ and $\epsilon_n\le\epsilon_{n-1}$, $[\epsilon_n-x_n(s,a)]^+ \le[\epsilon_{n-1}-(1-\gamma_n)\epsilon_{n-1}]^+ =\gamma_n\epsilon_{n-1}$.
The inequality uses the lower bound on $d_n(s,a)$ and monotonicity of $\epsilon_n$, and the equality is direct simplification. Hence, the sum of the deficits over all nonselected components is no larger than $(m-1)\gamma_n\epsilon_{n-1}$.

Raise every deficient nonselected component to $\epsilon_n$ and remove the same total mass from the selected component. Whenever $\|d_n-d^*\|\le r<1$, the selected component before projection is at least $(1-\gamma_n)(1-r/2)+\gamma_n$. Since $\epsilon_n\to0$ and $(m-1)\gamma_n\epsilon_{n-1}\to0$, there is a deterministic integer $N_2$ such that this constructed vector is feasible for every $n\ge N_2$. Its positive changes have total mass at most $(m-1)\gamma_n\epsilon_{n-1}$, and the selected component changes by the negative of the same amount. Its Euclidean distance from $x_n(s)$ is therefore no larger than $\sqrt2(m-1)\gamma_n\epsilon_{n-1}$.

Because $\Pi_n[x_n(s)]$ is the Euclidean projection, it cannot be farther from $x_n(s)$ than the feasible vector just constructed. The relation $\|z\|_1\le\sqrt m\|z\|_2$ gives
\begin{align}\label{Eq.policy_projection_bound_app}
\|d_{n+1}(s)-x_n(s)\|_1
&\le\sqrt{2m}(m-1)\gamma_n\epsilon_{n-1},
\qquad n\ge N_2.
\end{align}
The factor $\sqrt2$ comes from moving the same total mass into the deficient components and out of the selected component, and $\sqrt m$ converts the Euclidean norm to the one norm.

The nonselected mass of $x_n(s)$ is $(1-\gamma_n)R_n(s)$. Combining this identity with~\eqref{Eq.policy_projection_bound_app} yields
\begin{align}\label{Eq.policy_mass_recursion_app}
R_{n+1}(s)
&\le(1-\gamma_n)R_n(s)+C\gamma_n\epsilon_{n-1},
\end{align}
where $C=\sqrt{2m}(m-1)$. This inequality separates the contraction produced by the original policy update from the correction produced by the changing exploration set.

Choose a deterministic $N_1$ so large that $C\epsilon_{n-1}\le r/8$ for every $n\ge N_1$, and set $N_0(\omega) :=\max\{N_Q(\omega),N_1,N_2\}$.
Now take any realized index $N\ge N_0(\omega)$ satisfying $\|d_N(\omega)-d^*\|\le r/4$. Then $R_N(s)\le r/8$ for every state by~\eqref{Eq.policy_distance_mass_app}. If $R_n(s)\le r/4$,~\eqref{Eq.policy_mass_recursion_app} gives $R_{n+1}(s) \le(1-\gamma_n) r/4 +\gamma_n r/8 \le r/4$.
Thus, induction shows that $\|d_n-d^*\|\le r/2<r$ for all $n\ge N$. The strict gap neighborhood is never left, and the minimizing action remains $a^*(s)$ at every later iteration. Notice that the argument applies to every realized $N\ge N_0(\omega)$ satisfying the entry condition, which proves the probability statement in the theorem.

It remains to show convergence. Iterating~\eqref{Eq.policy_mass_recursion_app} from $N$ to $n$ gives $R_{n+1}(s) \le[\prod_{k=N}^n(1-\gamma_k)]R_N(s) +C\sum_{j=N}^n[\prod_{k=j+1}^n(1-\gamma_k)]\gamma_j\epsilon_{j-1}$.
The first term tends to zero because $\sum_n\gamma_n=\infty$. For the second term, fix $\eta>0$ and choose $N_\eta$ so large that $\epsilon_{j-1}\le\eta$ for every $j\ge N_\eta$. The finite part before $N_\eta$ tends to zero because it is multiplied by a product of the form $\prod(1-\gamma_k)$. The remaining part satisfies $C\eta\sum_{j=N_\eta}^n[\prod_{k=j+1}^n(1-\gamma_k)]\gamma_j =C\eta[1-\prod_{k=N_\eta}^n(1-\gamma_k)] \le C\eta$.
The equality follows by telescoping. Since $\eta$ is arbitrary, $R_n(s)\to0$. Equation~\eqref{Eq.policy_distance_mass_app} then gives $d_n(s)\to d^*(s)$ for every state, and hence $d_n\to d^*$. Proposition~\ref{prop_performance_difference} supplies the local optimality statement. This completes the proof.
\end{proof}

\begin{proof}[Proof of Theorem~\ref{prop_post_rate}]
Fix a state $s$, let $m=|\mathcal A(s)|$, set $C=\sqrt{2m}(m-1)$, and write $R_n(s)=\sum_{a\ne a^*(s)}d_n(s,a)$. After the minimizing actions have been identified, Equation~\eqref{Eq.policy_mass_recursion_app} gives
\begin{align}\label{Eq.rate_mass_app}
R_{n+1}(s)
&\le\left(1-\gamma_n\right)R_n(s)+C\gamma_n\epsilon_{n-1}.
\end{align}
For the power sequence $\epsilon_n=\epsilon_0/(n+1)^{p_\epsilon}$, $\epsilon_{n-1}/\epsilon_n =( (n+1)/n )^{p_\epsilon} =1+ p_\epsilon/n +O( 1/n^2 )$.
Because $p_\gamma<1$, $1/n=o(\gamma_n)$. Hence,
\begin{align}\label{Eq.epsilon_ratio_app}
\frac{\epsilon_{n-1}}{\epsilon_n}
&=1+o\left(\gamma_n\right).
\end{align}

Define $z_n=R_n(s)/\epsilon_{n-1}$. Dividing~\eqref{Eq.rate_mass_app} by $\epsilon_n$ and using~\eqref{Eq.epsilon_ratio_app} gives $z_{n+1} \le(1-\gamma_n) (\epsilon_{n-1}/\epsilon_n) z_n +C\gamma_n \epsilon_{n-1}/\epsilon_n \le(1- \gamma_n/2 )z_n+2C\gamma_n$
for all sufficiently large $n$. The last inequality follows because the ratio is $1+o(\gamma_n)$. Let $B_s=\max\{z_N,4C\}$ for an index $N$ after which the inequality holds. If $z_n\le B_s$, then $z_{n+1} \le(1- \gamma_n/2 )B_s+2C\gamma_n \le B_s$,
where the last inequality uses $B_s\ge4C$. Induction gives $z_n\le B_s$ for every $n\ge N$. Therefore, $R_n(s)=O(\epsilon_{n-1})$. Since $\epsilon_{n-1}/\epsilon_n\to1$, this is equivalent to $R_n(s)=O(\epsilon_n)$. By~\eqref{Eq.policy_distance_mass_app}, $\|d_n(s)-d^*(s)\|_1 =2R_n(s)=O(\epsilon_n)$ a.s.
The argument starts after the random action identification time. It does not bound the tail of that time and therefore does not give an expectation or high probability rate for the complete learning procedure. This completes the proof.
\end{proof}

\begin{proof}[Proof of Corollary~\ref{cor_mean_cvar}]
First, the added expected cost cancels from a difference of two $Q$-functions. Thus, $\|\mathcal T_{v,\lambda}Q-\mathcal T_{v,\lambda}Q'\|_w =\|\mathcal T_vQ-\mathcal T_vQ'\|_w \le\rho_Q\|Q-Q'\|_w$.
The equality holds because the same $\lambda\E[C(s,a)]$ is added to each operator value. The inequality is condition~\ref{A.Qcontraction}. Completeness of the finite dimensional Q space therefore gives a unique fixed point $Q_\lambda(v)$. Since the added term is also independent of $v$, the proof of Lemma~\ref{lemma_Q_continuity_app} gives $\|Q_\lambda(v)-Q_\lambda(v')\|_w \le ((2-\phi)/((1-\rho_Q)(1-\phi))) |v-v'|$.
Combining this bound with Lemma~\ref{lemma_continuity_app} gives continuity of $d\mapsto Q_\lambda^d$ with the same finite Lipschitz constant as before. The targets are generally different even though the constants in these difference bounds are unchanged.

To verify the radius, use the constants $M_1$ and $V$ defined in Section~\ref{sec.mean-cvar}. For $|v|\le V$ and any feasible pair, the triangle inequality and $|(C-v)^+|\le |C|+|v|$ give
\begin{align*}
\left|\tilde c(v,s,a)+\lambda\E[C(s,a)]\right| &\le |v|+\frac{\E[|C(s,a)|]+|v|}{1-\phi}+\lambda\E[|C(s,a)|] \le\frac{\left[1+\lambda(1-\phi)\right]M_1+(2-\phi)V}{1-\phi}.
\end{align*}
The first inequality bounds the positive part and uses $|\E C|\le\E|C|$; the second collects the two cost terms and the two threshold terms. By the fixed point equation and the contraction estimate,
\begin{align*}
\left\|Q_\lambda(v)\right\|_w
&\le\rho_Q\left\|Q_\lambda(v)\right\|_w+\left\|\mathcal T_{v,\lambda}0\right\|_w,\\
\left\|Q_\lambda(v)\right\|_\infty
&\le\frac{w_{\max}}{1-\rho_Q}\left\|\mathcal T_{v,\lambda}0\right\|_w\le\frac{w_{\max}\left(\left[1+\lambda(1-\phi)\right]M_1+(2-\phi)V\right)}{(1-\rho_Q)(1-\phi)}.
\end{align*}
The second line follows by moving the contraction term to the left and using $\|Q\|_\infty\le w_{\max}\|Q\|_w$. The last inequality uses $w(s,a)\ge1$ and the preceding cost bound. Equation~\eqref{Eq.mean_projection_bound} therefore leaves at least one unit between every required fixed point entry and the projection boundary.

We next retain the random local gain inside the moment calculation. Let $I_n(s,a)=\I\{s_n=s,a_n=a\}$, as in the Q proof, and let $\zeta_{n+1,\lambda}(s,a)$ denote the centered augmented Q observation on that event. Conditional on $\mathscr G_n$, the expected augmented observation is $(\mathcal T_{v_n,\lambda}Q_n)(s,a)$. The variance is bounded by the conditional second moment. Apply $(x+y)^2\le2x^2+2y^2$ to the augmented cost sample before taking its conditional expectation. On the event $I_n(s,a)=1$, its conditional second moment is at most $2\E[\widetilde C(v_n,s,a)^2\mid\mathscr G_n]+2\lambda^2\E[C(s,a)^2]$. The conditional law supplies the fixed component second moment in the last term. Using $|\widetilde C(v,s,a)|\le (|C(s,a)|+(2-\phi)|v|)/(1-\phi)$, and bounding the two Q minima by $Q_{\max}$, we obtain a finite constant $C_\lambda$ such that
\begin{align*}
\E\left[\zeta_{n+1,\lambda}(s,a)^2\mid\mathscr G_n\right]I_n(s,a)
&\le C_\lambda I_n(s,a)\left(1+v_n^2+\max_{s'\in\mathcal S,\,a'\in\mathcal A(s')}\E[C(s',a')^2]\right).
\end{align*}
This estimate does not require independence of the cost and next state, because the square of their sum was bounded before taking expectations. As in~\eqref{Eq.random_local_gain_sum_app}, distinct visits use distinct entries of the local gain sequence. Consequently,
\begin{align*}
\E\sum_n\beta_n(s,a)^2 I_n(s,a)\zeta_{n+1,\lambda}(s,a)^2
&\le C_\lambda\E\left[1+\sup_n v_n^2+\max_{s'\in\mathcal S,\,a'\in\mathcal A(s')}\E[C(s',a')^2]\right]\sum_j f_\beta(j)^2<\infty.
\end{align*}
The first inequality conditions at each iteration, applies the pathwise bound on the random gain sum, and then takes expectations. The last inequality uses Lemma~\ref{lemma_v_bound}, Assumption~\ref{A.cost}, and condition~\ref{A.beta}. Thus each augmented Q noise martingale converges. The direct cycle bound in Lemma~\ref{lemma_projected_cycle_app} applies with $Q_\lambda(v)$: its contraction, threshold sensitivity, noise tails, and interior projection comparison have all just been verified. The component step size condition is unchanged. The proof of Theorem~\ref{th_Q} therefore yields $\|Q_n-Q_\lambda(v_n)\|_w\to0$. Theorem~\ref{th_v} does not depend on the Q contraction or on the form of the chosen minimizing actions, beyond the policy increment bound, so it also applies to the augmented run. Together with the Lipschitz bound above, it yields both tracking statements in the corollary.

For the static performance statement, fix the candidate in~\eqref{Eq.mean_strict_gap}. Write $q_\lambda^*:=\operatorname{VaR}^{d_\lambda^*}$, $h_\lambda^*(s):=\min_{a\in\mathcal A(s)}Q_\lambda^{d_\lambda^*}(s,a)$, and $g_\lambda^*:=h_\lambda^*(s^{(0)})$. The augmented fixed point equation is $Q_\lambda^{d_\lambda^*}(s,a)+g_\lambda^* =\tilde c(q_\lambda^*,s,a)+\lambda\E[C(s,a)]+\sum_y p(y\mid s,a)h_\lambda^*(y)$.
Averaging under $\mu^{d_\lambda^*}(s)d_\lambda^*(s,a)$ cancels the continuation term with the state minimum, because the candidate is greedy for its own augmented $Q$-function. Hence $g_\lambda^*=\Psi_{d_\lambda^*}(q_\lambda^*)+\lambda\bar C^{d_\lambda^*}=J_\lambda(d_\lambda^*)$. Averaging the same equation under an arbitrary stationary policy $d$ instead gives
\begin{align*}
\Psi_d(q_\lambda^*)+\lambda\bar C^d-J_\lambda(d_\lambda^*)
&=\sum_s\mu^d(s)\sum_{a\in\mathcal A(s)}d(s,a)\left[Q_\lambda^{d_\lambda^*}(s,a)-h_\lambda^*(s)\right].
\end{align*}
Both cancellations use only $\mu^dP^d=\mu^d$, just as in~\eqref{Eq.threshold_control_bridge}. With $q_d:=\operatorname{VaR}^d$, subtraction of $\Psi_d(q_\lambda^*)-\Psi_d(q_d)$ now gives the exact identity
\begin{align*}
J_\lambda(d)-J_\lambda(d_\lambda^*)
&=\sum_s\mu^d(s)\sum_{a\in\mathcal A(s)}d(s,a)\left[Q_\lambda^{d_\lambda^*}(s,a)-h_\lambda^*(s)\right]-\left[\Psi_d(q_\lambda^*)-\Psi_d(q_d)\right].
\end{align*}
There is no extra approximation from the mean term: it enters linearly and is already included in the exact identity. The bracket is nonnegative because $q_d$ minimizes $\Psi_d$, and it is at most $C_{\rm p}\|d-d_\lambda^*\|^2$ by the derivative and quantile bounds in the proof of Proposition~\ref{prop_performance_difference}. Using the new strict gap in the other term gives $J_\lambda(d)-J_\lambda(d_\lambda^*) \ge ((\mu_{\min}\Delta_{Q,\lambda})/2) \|d-d_\lambda^*\|-C_{\rm p}\|d-d_\lambda^*\|^2$.
Here $\mu_{\min}=\delta_p\min_s\nu(s)>0$ follows from minorization. As before, the factor $1/2$ follows from $\|d(s)-\delta(a_\lambda^*(s))\|_1=2\sum_{a\ne a_\lambda^*(s)}d(s,a)$. The right hand side is strictly positive for every sufficiently close policy different from $d_\lambda^*$, proving strict local optimality for $J_\lambda$.

Finally, choose $r_\lambda>0$ small enough that $\|Q_\lambda^d-Q_\lambda^{d_\lambda^*}\|_\infty\le\Delta_{Q,\lambda}/4$ whenever $\|d-d_\lambda^*\|\le r_\lambda$. On a probability one event, Q tracking makes $\|Q_{n+1}-Q_\lambda^{d_n}\|_\infty$ smaller than $\Delta_{Q,\lambda}/8$ for all sufficiently large $n$. To see the use of $Q_{n+1}$, add and subtract $Q_\lambda^{d_{n+1}}$ and use the policy increment bound~\eqref{Eq.policy_increment_app}. Thus every minimizing action used by the update is $a_\lambda^*(s)$ while the policy remains in this neighborhood. The projection comparison in the proof of Theorem~\ref{th_opt} then gives, with $R_n(s):=\sum_{a\ne a_\lambda^*(s)}d_n(s,a)$, $R_{n+1}(s) \le(1-\gamma_n)R_n(s)+C\gamma_n\epsilon_{n-1}$.
Choose the finite $N_{0,\lambda}(\omega)$ so that the tracking bound holds and the forcing term is small enough for the neighborhood induction in that proof. Starting at any realized $N\ge N_{0,\lambda}(\omega)$ with $\|d_N-d_\lambda^*\|\le r_\lambda/4$, the same induction keeps the policy in the larger neighborhood. Expanding the last recursion gives a vanishing product multiplying $R_N(s)$ and a weighted sum of $\epsilon_{n-1}$ that tends to zero, because $\sum_n\gamma_n=\infty$ and $\epsilon_n\to0$. Hence $R_n(s)\to0$ for every state and $d_n\to d_\lambda^*$. The candidate, gap, radius, and entry time in this argument all belong to the augmented problem; no corresponding fact for $\lambda=0$ was assumed to carry over. This completes the proof.
\end{proof}

\end{APPENDICES}

\clearpage
\setcounter{nowappendix}{0}
\setcounter{page}{1}
\renewcommand{\thepage}{S\arabic{page}}
\makeatletter
\setlength{\@fptop}{0pt}
\setlength{\@fpsep}{11pt}
\setlength{\@fpbot}{0pt plus 1fil}
\makeatother
\setcounter{section}{0}
\renewcommand{\thesection}{S\arabic{section}}
\setcounter{table}{0}
\setcounter{figure}{0}
\setcounter{equation}{0}
\renewcommand{\thetable}{S\arabic{table}}
\renewcommand{\thefigure}{S\arabic{figure}}
\renewcommand{\theequation}{S\arabic{equation}}
\pdfbookmark[0]{Supplementary material}{supplementary-material}
\begin{center}
{\Large\bfseries Supplementary material\par}
\medskip
{\large Long-Run Conditional Value at Risk Reinforcement Learning\par}
\medskip
Qixin Wang, Hao Cao, Jian-Qiang Hu, Mingjie Hu, and Li Xia\par
\end{center}
\noindent Implementation details and additional numerical comparisons for the machine replacement experiment in Section~\ref{sec_machine} of the main manuscript, including alternative objectives and final policies. 

\section{Additional Numerical Results}\label{sec_numerical_details}
\noindent\textbf{Configurations and implementation.}
Table~\ref{tab:recorded_machine_inputs} gives the inputs used in the additional machine comparison in Section~\ref{sec_machine}. The observation index is $k=1,\ldots,10^6$, and $m$ is the number of visits to the updated component, including the current visit. Both the initial state and reference state are $s^{(1)}$. For CRL, exploration denotes the lower bound on each feasible action probability. For RED, it denotes the weight of the uniform distribution in its greedy behavior policy, and the VaR and CVaR scalar recursions use the same step size.

\begin{table}[H]
\centering\normalsize
\caption{Recorded inputs for the additional machine comparison.}\label{tab:recorded_machine_inputs}
\setlength{\tabcolsep}{6pt}
\renewcommand{\arraystretch}{1.12}
\begin{tabular}{@{}>{\raggedright\arraybackslash}p{3.2cm}ccccc@{}}
\toprule
Method and costs & VaR step & Q step & Policy step & Exploration & Q bound \\
\midrule
CRL, Gaussian & $10/k^{0.9}$ & $0.8/(m+1)^{0.8}$ & $0.5/k^{0.99}$ & $0.5/k^{0.999}$ & 1300 \\
CRL, Student $t_5$ & $20/k^{0.8}$ & $0.8/(m+1)^{0.65}$ & $0.8/k^{0.9}$ & $0.3/k^{0.3}$ & 1300 \\
RED, both costs & $0.2/k^{0.75}$ & $0.5/(m+1)^{0.55}$ & --- & $0.1$ & none \\
\bottomrule
\end{tabular}
\end{table}

The initialization names specify Q, rather than an independently sampled policy. For CRL, legacy sets every feasible Q entry to zero, replace sets the entries for replacement to $-0.5$ and the other feasible entries to zero, and random uses $0.1$ times independent standard normal entries. RED uses the opposite Q signs. Ties select the smallest feasible action index. In every case, CRL starts from the uniform feasible policy, whereas RED starts from the mixture of its greedy policy and the uniform feasible policy with weights $0.9$ and $0.1$. Both methods use $v_0=0$ and zero visit counts. RED also starts its scalar CVaR estimate at zero. Infeasible actions are excluded from all updates and comparisons. For each method and cost distribution, the recorded pilot protocol considers 12 configurations, each with $10^6$ observations. Selection first minimizes the mean final greedy CVaR. Values within $10^{-10}$ are compared using the time average CVaR gap computed from the recorded policy change events. The selected candidates are 1 and 10 for Gaussian and Student CRL, respectively, and 7 for RED under both distributions. 

The historical machine runs use $Q_{\max}=500$ for CRL, M-CRL, and MRL. Their archived implementations project only the updated Q component. The additional machine CRL runs use $Q_{\max}=1300$, with the executed update source identified by the hash in the saved records. RED does not project Q. All saved projection counts in the reported formal runs are zero. MRL uses the Q and policy updates of CRL with the sampled cost in place of the CVaR transform and omits the VaR recursion. Random Q initialization uses the separate stream described above. These runs use \mbox{MATLAB} R2023b Update 5 with a C update kernel on Windows, and the recorded processor identifier is Intel64 Family 6 Model 186 Stepping 2, with 16 logical processors. Each reported time is one elapsed time for a batch of 30 trajectories. 

\noindent\textbf{Budgets and comparisons.}
Table~\ref{tab:fixed_budget_supp} reports the same percentage gap at all six selected budgets. The intervals are untruncated Student intervals for the mean over replications, so a lower endpoint can be negative even though each deterministic policy gap is nonnegative. They describe sampling uncertainty in the repeated learning outcomes. 

\begin{table}[htbp]
\centering\normalsize\setlength{\tabcolsep}{3pt}
\caption{CVaR gaps at selected budgets.}\label{tab:fixed_budget_supp}
\begin{tabular}{@{}lcccccc@{}}\toprule
Method & $10^4$ & $5\times10^4$ & $10^5$ & $2\times10^5$ & $5\times10^5$ & $10^6$ \\\midrule
\multicolumn{7}{@{}l}{\textbf{Gaussian}} \\
CRL & \shortstack{0.0000 \\ {\scriptsize (0.0000, 0.0000)}} & \shortstack{0.0000 \\ {\scriptsize (0.0000, 0.0000)}} & \shortstack{0.0000 \\ {\scriptsize (0.0000, 0.0000)}} & \shortstack{0.0000 \\ {\scriptsize (0.0000, 0.0000)}} & \shortstack{0.0000 \\ {\scriptsize (0.0000, 0.0000)}} & \shortstack{0.0000 \\ {\scriptsize (0.0000, 0.0000)}} \\
RED & \shortstack{0.9909 \\ {\scriptsize (0.7696, 1.2121)}} & \shortstack{0.0000 \\ {\scriptsize (0.0000, 0.0000)}} & \shortstack{0.0000 \\ {\scriptsize (0.0000, 0.0000)}} & \shortstack{0.0000 \\ {\scriptsize (0.0000, 0.0000)}} & \shortstack{0.0000 \\ {\scriptsize (0.0000, 0.0000)}} & \shortstack{0.0633 \\ {\scriptsize (-0.0267, 0.1533)}} \\
\midrule
\multicolumn{7}{@{}l}{\textbf{Student $t_5$}} \\
CRL & \shortstack{0.0000 \\ {\scriptsize (0.0000, 0.0000)}} & \shortstack{0.0000 \\ {\scriptsize (0.0000, 0.0000)}} & \shortstack{0.0000 \\ {\scriptsize (0.0000, 0.0000)}} & \shortstack{0.0000 \\ {\scriptsize (0.0000, 0.0000)}} & \shortstack{0.0000 \\ {\scriptsize (0.0000, 0.0000)}} & \shortstack{0.0000 \\ {\scriptsize (0.0000, 0.0000)}} \\
RED & \shortstack{0.8089 \\ {\scriptsize (0.5785, 1.0394)}} & \shortstack{0.0280 \\ {\scriptsize (-0.0293, 0.0853)}} & \shortstack{0.0280 \\ {\scriptsize (-0.0293, 0.0853)}} & \shortstack{0.0280 \\ {\scriptsize (-0.0293, 0.0853)}} & \shortstack{0.0000 \\ {\scriptsize (0.0000, 0.0000)}} & \shortstack{0.0280 \\ {\scriptsize (-0.0293, 0.0853)}} \\
\bottomrule\end{tabular}\end{table}

\begin{table}[htbp]
\centering\normalsize
\caption{Differences CRL minus RED.}\label{tab:paired_supp}
\begin{tabular}{llrrr}\toprule
Distribution & Quantity & Mean difference & Lower 95\% & Upper 95\% \\\midrule
Gaussian & Final CVaR & -0.00964 & -0.02334 & 0.00406 \\
Gaussian & Final mean cost & 0.09338 & -0.03932 & 0.22608 \\
Gaussian & Average gap (\%) & -0.03378 & -0.04002 & -0.02754 \\
Gaussian & Observed changes & -4.400 & -5.344 & -3.456 \\
\midrule
Student $t_5$ & Final CVaR & -0.00426 & -0.01297 & 0.00445 \\
Student $t_5$ & Final mean cost & 0.04669 & -0.04880 & 0.14219 \\
Student $t_5$ & Average gap (\%) & -0.03650 & -0.04302 & -0.02997 \\
Student $t_5$ & Observed changes & -6.267 & -7.615 & -4.919 \\
\bottomrule\end{tabular}\end{table}

\noindent\textbf{Other objectives and final policies.}
For comparison with other operating objectives, MRL uses the mean cost criterion and M-CRL uses $J=\operatorname{CVaR}+0.3\,\mathrm{Mean}$. Table~\ref{tab:all_methods_supp} records the final greedy policy performance reported for all four methods. 

\begin{table}[htbp]
\centering\normalsize
\caption{Final greedy policy performance.}\label{tab:all_methods_supp}
\begin{tabular}{llrrrr}\toprule
Distribution & Method & VaR & CVaR & Mean cost & $J$ \\\midrule
Gaussian & CVaR benchmark & 14.6883 & 15.2239 & 8.1252 & 17.6615 \\
Gaussian & Mean / J benchmark & 15.1138 & 15.4741 & 6.0100 & 17.2771 \\
Gaussian & CRL & 14.6883 & 15.2239 & 8.1252 & 17.6615 \\
Gaussian & RED & 14.7059 & 15.2336 & 8.0318 & 17.6431 \\
Gaussian & MRL & 15.1533 & 15.5021 & 6.0109 & 17.3054 \\
Gaussian & M-CRL & 15.1138 & 15.4741 & 6.0100 & 17.2771 \\
Student $t_5$ & CVaR benchmark & 14.7420 & 15.2120 & 8.1252 & 17.6495 \\
Student $t_5$ & Mean / J benchmark & 15.0929 & 15.4381 & 6.0100 & 17.2411 \\
Student $t_5$ & CRL & 14.7420 & 15.2120 & 8.1252 & 17.6495 \\
Student $t_5$ & RED & 14.7493 & 15.2162 & 8.0785 & 17.6398 \\
Student $t_5$ & MRL & 15.1258 & 15.4659 & 6.0109 & 17.2692 \\
Student $t_5$ & M-CRL & 15.0929 & 15.4381 & 6.0100 & 17.2411 \\
\bottomrule\end{tabular}\end{table}

Table~\ref{tab:policy_frequency_supp} gives the policy frequencies behind these results. Here R denotes retain and X denotes replace, with actions listed in state order. Each method and distribution has 30 final policies. 
\begin{table}[htbp]
\centering\normalsize
\caption{Frequencies of the final deterministic policies.}\label{tab:policy_frequency_supp}
\begin{tabular}{llcc}\toprule
Distribution & Method & Policy & Number of runs \\\midrule
Gaussian & CRL & RRRRRX & 30 \\
Gaussian & RED & RRRRRX / RRRRXX & 28 / 2 \\
Gaussian & MRL & RRRXXX / RRXXXX & 22 / 8 \\
Gaussian & M-CRL & RRRXXX & 30 \\
Student $t_5$ & CRL & RRRRRX & 30 \\
Student $t_5$ & RED & RRRRRX / RRRRXX & 29 / 1 \\
Student $t_5$ & MRL & RRRXXX / RRXXXX & 22 / 8 \\
Student $t_5$ & M-CRL & RRRXXX & 30 \\
\bottomrule\end{tabular}\end{table}

Table~\ref{tab:current_supp} reports the stationary VaR, CVaR, and mean cost of the final randomized policies obtained by CRL, RED, MRL, and M-CRL under the Gaussian and Student $t_5$ cost models. Each policy is evaluated separately, and the reported values are averages over 30 independent replications.

\begin{table}[htbp]
\centering\normalsize
\caption{Final randomized policy performance.}\label{tab:current_supp}
\begin{tabular}{llrrr}\toprule
Distribution & Method & VaR & CVaR & Mean cost \\\midrule
Gaussian & CRL & 14.6897 & 15.2246 & 8.1238 \\
Gaussian & RED & 14.8828 & 15.3278 & 7.9453 \\
Gaussian & MRL & 15.1789 & 15.5195 & 6.0129 \\
Gaussian & M-CRL & 15.1141 & 15.4743 & 6.0114 \\
Student $t_5$ & CRL & 14.7431 & 15.2125 & 8.1238 \\
Student $t_5$ & RED & 14.9001 & 15.3006 & 7.9831 \\
Student $t_5$ & MRL & 15.1461 & 15.4818 & 6.0129 \\
Student $t_5$ & M-CRL & 15.0931 & 15.4383 & 6.0114 \\
\bottomrule\end{tabular}\end{table}

\end{document}